\documentclass[twoside]{amsart}
\usepackage{float}
\usepackage{graphicx}
\input{diagrams.tex}
\diagramstyle[Postscript=dvips]
\newarrow{Eqto}{=}{=}{=}{=}{=}
\newarrow{Line}{-}{-}{-}{-}{-}
\newarrow{Dashto}{}{dash}{}{dash}>
\newarrow{Dotsto}{.}{.}{.}{.}>
\newtheorem{theorem}{\sc Theorem}[section]
\newtheorem{proposition}[theorem]{\sc Proposition}
\newtheorem{lemma}[theorem]{\sc Lemma}
\newtheorem{corollary}[theorem]{\sc Corollary}
\theoremstyle{definition}
\newtheorem{definition}[theorem]{\sc Definition}
\newtheorem{definitions}[theorem]{\sc Definitions}
\newtheorem{example}[theorem]{\sc Example}
\newtheorem{examples}[theorem]{\sc Examples}
\newtheorem{problem}[theorem]{\sc Problem}
\theoremstyle{remark}
\newtheorem{remark}[theorem]{\sc Remark}
\newtheorem{remarks}[theorem]{\sc Remarks}
\newtheorem{claim}[theorem]{}

\newcommand{\Hb}[2]{\mathbf{#1}^{#2}}
\newcommand{\rest}[1]{\!\mid_{_#1}}
\newcommand{\ro}[1]{{}_{\left( #1\right)}}
\newcommand{\un}[1]{{}_{\langle #1\rangle}}
\newcommand{\ep}[1]{\varepsilon_{#1}}

\newcommand{\de}[1]{_{({#1})}}

\def\cmc{{}_C\frak{M}_C}
\def\Hom{\mathrm{Hom}}

\def\calM{\frak{M}}
\def\Hu{\mathbf{H}^\bullet}

\def\cmc{{}^C\!\calM^C}
\def\cmct{({}^C\!\calM^C,\otimes,\mathbf{1})}
\def\lin{--$\,$}
\def\Hu{\mathbf{H}^\bullet}
\def\ker{\mathrm{Ker\,}}
\def\bu{\bullet}

\def\calM{\mathcal{M}}

\def\mm{\mathfrak{M}}

\def\hmh{{}_{H}^{H}\mathfrak{M}_{H}^{H}}

\def\yd{{}_H^H\mathcal{YD}}

\def\ot{\otimes}

\def\ra{\rightarrow}

\def\hm{{}^H \mathfrak{M}}

\def\cotH{\square_H}

\def\td{\widetilde{\delta}}

\def\sm{R\#H}
\def\dsm{\Delta_{R\#H}}

\def\de{\delta}
\def\om{\omega}

\begin{document}
\pagestyle{headings}
\title{Hochschild cohomology of algebras in monoidal categories
and splitting morphisms of bilgebras}
\author{A. Ardizzoni, C. Menini}
\author{D. \c{S}tefan}
\subjclass{Primary 16; Secondary 16}
\date{January 6, 2003}
\begin{abstract}
The main goal of this paper is to investigate the structure of
Hopf algebras with the property that either its Jacobson radical
is a Hopf ideal or its coradical is a subalgebra. In order to do
that we define the Hochschild cohomology of an algebra in an
abelian monoidal category. Then we characterize those algebras
which have dimension less than or equal to 1 with respect to
Hochschild cohomology. Now let us consider a Hopf algebra $A$ such
that its Jacobson radical $J$ is a nilpotent Hopf ideal and
$H:=A/J$ is a semisimple algebra. By using our homological
results, we prove that the canonical projection of $A$ on $H$ has
a section which is an $H$\lin colinear algebra map. Furthermore,
if $H$ is cosemisimple too, then we can choose this section to be
an $(H,H)$\lin bicolinear algebra morphism. This fact allows us to
describe $A$ as a `generalized bosonization' of a certain algebra
$R$ in the category of Yetter\lin Drinfeld modules over $H$. As an
application we give a categorical proof of Radford's result about
Hopf algebras with projections. We also consider the dual
situation. In this case, many results that we obtain hold true for
a large enough class of $H$\lin module coalgebras, where $H$ is a
cosemisimple Hopf algebra.
\end{abstract}
\maketitle
\section*{Introduction} \markboth{\sc{A. ARDIZZONI, C. MENINI AND D.
\c{S}TEFAN}} {\sc{HOCHSCHILD COHOMOLOGY AND SPLITTING MORPHISMS OF
BIALGEBRAS}} Let $H$ be a Hopf algebra. The categories $\yd$ and
$\hmh$, of Yetter$\,$\lin Drinfeld modules and respectively Hopf
bimodules, appeared, in particular, as an attempt to construct new
solutions to the Yang$\,$\lin Baxter equation. Nowadays we can
recognize their most important properties into the definition of
braided categories, a very general and abstract setting useful,
not only for providing new solutions to the Yang$\,$\lin Baxter
equation, but also in many other areas of mathematics, like the
theory of quantum groups and low dimensional topology.

Partially motivated by these applications, the theory of Hopf
algebras knew in 80's an outstanding development. Besides many
striking results obtained since then, we would like to recall,
more or less chronologically, a few of them that will play a very
important role in our paper.\smallskip

$\bu$ The description of the coradical filtration of a pointed
coalgebra, result due to Taft and Wilson \cite{TW}, that is
crucial in the classification of  finite dimensional pointed Hopf
algebras.\smallskip

$\bu$ The characterization of bialgebras with projection
\cite{Rad} -- later Majid \cite{Maj1} showed that this result can
be interpreted in terms of bialgebras in a braided
category.\smallskip

$\bu$ The equivalence of braided categories $\yd\simeq\hmh$
\cite{Sch1} and \cite{AD}.\smallskip

$\bu$ The classification of certain classes of pointed Hopf
algebras of finite dimension. One of the used method is the
`lifting' method \cite{AS1}, \cite{AS2}, \cite{AS3}, \cite{AS4}.
Let $A$ be a Hopf algebra such that its coradical is a Hopf
subalgebra $H$. Then the coradical filtration of $A$ is a
filtration of Hopf algebras, so $\mathrm{gr}\,A$ is a graded Hopf
algebra. One of the main steps of the `lifting' method is to
describe $\mathrm{gr}\,A$, by using the second mentioned result,
as the `bosonization' of a certain Hopf algebra $R$ in $\yd$ by
$H$. The next step is to find all Hopf algebras $A$ having a given
graded Hopf algebra $\mathrm{gr}\,A$.\smallskip

$\bu$ Let $A$ be a finite dimensional Hopf algebra over a field
$k$ of characteristic zero whose coradical, say $H$, forms a Hopf
subalgebra. Then the left $H$-module coalgebra $A$ is a cosmash in
the sense that there exists a $H$-linear coalgebra map
$\gamma:A\rightarrow H$ such that $\gamma\rest{H}=\mathrm{Id_H}$,
see \cite{SvO}. In \cite{Mas} it is shown, with a different
method, that the above result still holds true without any
assumption on the dimension of $A$ and
$\mathrm{char}\,k$.\smallskip

$\bu$ For a Hopf algebra $A$ a conjectural formula for $A_1$, the
first component of the coradical filtration of $A$, is proposed in
\cite{AS5}. This formula is proved in the same paper in the case
when $A$ is a graded Hopf algebra such that its coradical is a
Hopf subalgebra of $A$. In \cite{CDMM} the conjecture is proved in
the ungraded case.\medskip

One of the main aims of this paper is to strengthen some of the
results that we mentioned above, by combining the formalism of
homological algebra and monoidal categories.

Let $(\calM,\ot,\mathbf{1})$ be a monoidal category. We start by
recalling some basic facts about algebras and $(A,A)$\lin
bimodules in $\calM$. Then, for an algebra $A$ and an $(A,A)$\lin
bimodule $M$ in $\calM$, we define the Hochschild cohomology
$\Hu(A,M)$ of $A$ with coefficients in $M$, and show that many
classical results can be extended to this more general context.
For example, we call an algebra $A$ separable if $A$ is
$%
\mathcal{E}$-projective as an $(A,A)$\lin bimodule in $\calM$.
Such an $A$ is characterized by the fact that its Hochschild
dimension is zero, that is $\Hb{H}{1}(A,M)=0$, for every bimodule
$M$. Also, by defining in an appropriate way Hochschild extensions
of $A$ with kernel $M$, we prove that the set of their equivalence
classes is in one$\,$\lin to$\,$\lin one correspondence with
$\Hb{H}{2}(A,M).$ An algebra $A$ will be called (non$\,$\lin
commutative) formally smooth if it has no non$\,$\lin trivial
extensions. In the particular case when $\calM$ is the category of
$K$\lin vector spaces, we recognize the definition of
quasi$\,$\lin free algebras introduced by J. Cuntz and D. Quillen
in \cite{CQ}. The first important results, Theorem \ref{X coro
UniE} and Theorem \ref{X formally smooth teo}, give different
equivalent characterizations of formally smooth algebras in a
monoidal category, and can be interpreted as generalizations of
Wedderburn$\,$\lin Malcev Theorem. We obtain immediately that
separable algebras are formally smooth. In view of these theorems,
roughly speaking, formally smoothness is useful to prove that
certain algebra morphisms have algebra morphism sections in the
category that we work in.

In the third section of the paper we apply this technique to
produce algebra sections in the case when $\calM$ is either
$\mm^H$ (the category of right $H$\lin comodules) or $^H\mm^H$
(the category of $(H,H)$\lin bicomodules), where $H$ is a
semisimple Hopf algebra. Let $A$ be a Hopf algebra such that its
Jacobson radical $J$ is a Hopf ideal. We denote the Hopf algebras
quotient $A/J$ by $H$. If $J$ is nilpotent and $H$ is semisimple,
in the main application of this section, Corollary
\ref{co:SectonBicolin}, we prove that the canonical projection
$\pi:A \ra H$ has a section which is a morphism of algebras in
$\mm^H$. This is a generalization of one of the main results of
\cite{SvO}. The second part of the corollary establishes that $\pi
$ has a section which is an algebra map in $^H\mm^H$, if we assume
in addition that $H$ is cosemisimple too (in fact we prove a more
general result, that we shall not mention here to simplify the
exposition). This corollary recall us \cite{Rad}, where it is
assumed that a Hopf algebra morphism $\pi:A\ra H$ has a section
$\sigma:H\ra A$ which is a morphism of Hopf algebras. In
\cite{Rad} it is shown that there is a bialgebra $R$ in $\yd$ such
that $A$ is the smash product algebra and the smash product
coalgebra of $R$ by $H$.

Taking into account Corollary \ref{co:SectonBicolin}, it is
natural to look for a similar description of a bialgebra $A$,
supposing that $\pi:A\ra H$ has a section $\sigma$ which is only a
morphism of  algebras in $^H\mm^H$. This will be done in the
fourth section of the paper.

The starting point is the  simple observation that $A$ becomes in
a natural way an object in $\hmh$. Of course the left and right
comodule structures are induced by $\pi$. Since $\sigma$ is a
morphism of algebras, $A$ is a bimodule over $H$, and the fact
that $\sigma$ is a morphism of bicomodules is enough to have the
required compatibility relations. By using the equivalence
$\hmh\simeq\yd$ we have  $A\simeq R\ot H$ (isomorphism in $\hmh$),
where $R=A^{co(H)}$. Moreover, the multiplication of $A$ is a
morphism in $\hmh$ and the unit of $A$ is in $R$. Therefore $R$
becomes an algebra in $\yd$ and $A$ can be identified as an
algebra with the smash product $R\# H$. We can not repeat this
argument for the coalgebra structure since $\Delta$ is only
$(H,H)$\lin colinear. Thus, by identifying $A$ and $R\# H$ as
algebras, the problem of describing all bialgebras $A$ as above is
equivalent to find all coalgebras structures on $R\#H$ such that
the comultiplication is a morphism of $(H,H)$\lin bicomodules. We
prove that $\dsm$ is uniquely determined by a pair of $K$\lin
linear maps $\de:R\ra R\ot R$ and $\om:H\ra R\ot R$. Let $\ep{}$
be the restriction of the counit of $A$ to $R$. The properties of
$\de$, $\om$ and $\ep{}$ necessary to get a bialgebra structure on
$R\#H$ are listed in Definition \ref{de:YDq}. The result that we
obtain is stated in Theorem \ref{pro Boss}.

Let us remark that $(R,\de,\ep{})$ is  not a coalgebra since $\de$
is not coassociative in general. In fact the coassociativity rule
is broken by the map $\om$, who is a normalized (non\lin
commutative) 1\lin cocycle of $H$. If $\om$ is the trivial cocycle
then $R$ is a bialgebra in $\yd$ and $A$ is isomorphic as an
algebra and coalgebra with the smash product, i.e. with the
`bosonization' of $R$ by $H$. The main results of this section are
Theorem \ref{th:BialgInBimod2} and Theorem \ref{te:A=YDq}.

In the last section of the paper we state the dual results and use
them to prove some applications. We start by defining the
Hochschild cohomology of a coalgebra in a monoidal category
$\calM$ and by giving homological characterizations of coseparable
and formally smooth coalgebras. As a consequence, by taking the
monoidal category $\mathcal M$ to be either $\mm_H$ or $_H\mm_H$,
we prove Theorem \ref{te:SectionBicomod}. According to this
theorem, under some assumptions on $H$, for every coalgebra $C$ in
$\calM$ such that $C_0=H$, there exists a morphism of coalgebras
$\pi:C\ra H$ in $\calM$ so that $\pi\rest{H}=\mathrm{Id}_H$. The
first assertion of Theorem \ref{te:SectionBicomod} had already
been proved by A. Masuoka in the case when $C$ is the underlying
coalgebra structure of a Hopf algebra $A$ with with the property
that its coradical is a subalgebra. Now we can describe the
coradical filtration of such a coalgebra $C$ as in Theorem
\ref{te:Corad_Filt}. Finally, we prove that a Hopf algebra $A$,
having the coradical a semisimple and cosemisimple Hopf
subalgebra, is as a Hopf algebra, not only as a coalgebra, a kind
of smash product, see Theorem \ref{te:Boson}. We expect that this
last result is strongly connected with the lifting method
introduced by N. Andruskiewitsch and H.J. Schneider. Probably
Theorem \ref{te:Boson} can be used to get direct information about
a Hopf algebra $A$ with the property that its coradical is a
subalgebra, skipping the step when the associated graded Hopf
algebra $gr A$ is investigated.

Parts of this paper were presented by the second author in the
talks she delivered at the meetings "2003 Spring Eastern Sectional
AMS Meeting" (Special Session on Hopf Algebras and Quantum
Groups), New York, NY (U.S.A.), April 12-13, 2003 and "Non
Commutative Geometry and Rings", Almeria (Spain), September 2 - 7,
2002.\medskip

\textbf{Notation.} \ In a category $\mathcal{M}$ the set of
morphisms from $X$ to $Y$ will be denoted by $\mathcal{M}(X,Y).$
If $X$ is an object in $\mathcal{M}$ then the functor
$\mathcal{M}(X,-)$ from $\mathcal{M}$ to $\mathfrak{Sets}$
associates to any morphism $u:U\rightarrow V$ in $\mathcal{M}$ the
function that will be denoted by $\mathcal{M}(X,u).$
\section{Hochschild cohomology in monoidal categories}
In this section we define  and study the Hochschild cohomology of
an algebra in a monoidal category. We start by recalling the
definitions of monoidal categories and of algebras in such
categories. In order to define Hochschild cohomology we will use
relative homological algebra, for details on this matter see
\cite[Chapter IX]{HiSt}.

\begin{claim}
\label{MonCat}A \emph{monoidal category} means  a category
$\mathcal{M}$ that is endowed with a functor $\otimes
:\mathcal{M}\times \mathcal{M}\rightarrow \mathcal{M}$, an object
$\mathbf{1}\in \mathcal{M}$ and functorial isomorphisms:
$a_{X,Y,Z}:(X\otimes Y)\otimes Z\rightarrow X\otimes (Y\otimes
Z),$ $l_{X}:\mathbf{1}\otimes X\rightarrow X$ and $r_{X}:X\otimes \mathbf{1}%
\rightarrow X.$ The functorial morphism $a$ is called the \emph{%
associativity constraint }and\emph{\ }satisfies the \emph{Pentagon Axiom, }%
that is the first diagram below is commutative, for every $U,\,V,$
$W,$ $X$ in $\mathcal{M.}$ The morphisms $l$ and $r$ are called
the \emph{unit constraints} and they are assumed to satisfy the
\emph{Triangle Axiom, }i.e. the second diagram is commutative. The
object $\mathbf{1}$ is called the \emph{unit} of $\mathcal{M}$.
\begin{equation*}
\begin{diagram}[small] (U\otimes (V\otimes W))\otimes X&\lTo^{a_{U, V,
W}\otimes X}& ((U\otimes V)\otimes W)\otimes X\\ & &
\dTo>{a_{U\otimes V, W,X}}\\ \dTo<{a_{U,V\otimes W,X}}& &
(U\otimes V)\otimes (W\otimes X) \\ & & \dTo>{a_{U, V, W\otimes
X}}\\ U\otimes ((V\otimes W)\otimes X)&\rTo_{U\otimes a_{ V,
W,X}}&U\otimes (V\otimes (W\otimes X))
\end{diagram}\hspace*{0.5cm}
\begin{diagram}[h=2.5em,w=2.5em] (V\otimes
\mathbf{1})\otimes W&
&\rTo^{a_{V,\mathbf{1},W}}&&V\otimes(\mathbf{1}\otimes W)\\
&\rdTo<{r_V\otimes W}&&\ldTo>{V\otimes l_W}& \\ & & V\otimes W & &
\end{diagram}
\end{equation*}
For details on monoidal categories we refer to \cite[Chapter
XI]{Ka}. A monoidal category is called \emph{strict} if the
associativity constraint and unit constraints are the
corresponding identity morphisms.
\end{claim}

\begin{claim}
\label{cl:CohThm}As it is noticed in \cite[p. 420]{Maj2}, the
Pentagon Axiom solves the consistency problem that appears because
there are two ways to go from $((U\otimes V)\otimes W)\otimes X$
to $U\otimes (V\otimes (W\otimes X)).$ The coherence theorem, due
to S. Mac Lane, solves the similar problem for the tensor product
of an arbitrary number of objects in $\mathcal{M}.$ Accordingly
with this theorem, we can always omit all brackets and simply
write $X_{1}\otimes \cdots \otimes X_{n}$ for any object obtained from $%
X_{1},\ldots ,X_{n}$ by using $\otimes $ and brackets. Also as a
consequence of the coherence theorem, the morphisms $a,$ $l,$ $r$
take care of themselves, so they can be omitted in any computation
involving morphisms in $\mathcal{M.}$
\end{claim}

\begin{claim}
\label{cl:MonFun}A \emph{monoidal functor} between two monoidal categories $(%
\mathcal{M},\otimes ,\mathbf{1},a,l,r\mathbf{)}$ and $\mathbf{(}\mathcal{M}%
'\mathcal{,}\otimes ,\mathbf{1},a,l,r\mathbf{)}$ is a triple $%
(F,\phi _{0},\phi _{2}),$ where $F:\mathcal{M}\rightarrow \mathcal{M}%
'$ is a functor, $\phi _{0}:\mathbf{1\rightarrow F(1)}$ is an
isomorphism such that
\begin{equation*}
\begin{diagram}[h=2em,w=1.5em]  \mathbf{1}\otimes F(U) &
\rTo^{ l_{F(U)} } & F(U) \\ \dTo<{
\phi_0\otimes F(U) } & & \uTo>{ F(l_U) } \\ F(\mathbf{1})\otimes
F(U) &
\rTo_{\phi_2( \mathbf{1} ,U) } & F( \mathbf{1}\otimes U) \\
\end{diagram}%
\hspace*{1cm}
\begin{diagram}[h=2em,w=1.5em] F(U)
\otimes \mathbf{1}
& \rTo^{ r_{F(U)} } & F(U) \\ \dTo<{F(U) \otimes\phi_0 } & &
\uTo>{ F(r_U) } \\ F(U)\otimes F(\mathbf{1})& \rTo_{\phi_2(U,
\mathbf{1}) } &
F(U\otimes \mathbf{1}) \\
\end{diagram}
\end{equation*}
are commutative, and $\phi _{2}(U,V):F(U\otimes V)\rightarrow
F(U)\otimes F(V)$ is a family of functorial isomorphisms such that
the following diagram is commutative.
\begin{equation*}
\begin{diagram}[h=2em,w=1.5em] (F(U)\otimes F(V))\otimes
F(W)&\rTo^{a_{F(U),F(V),F(W)}}&F(U)\otimes (F(V)\otimes F(W))\\
\dTo<{\phi_2(U,V)\otimes F(W)} & & \dTo>{F(U)\otimes
\phi_2(V,W)}\\ F(U\otimes V)\otimes F(W)& &F(U)\otimes F(V\otimes
W) \\ \dTo<{\phi_2(U\otimes V,W)} & & \dTo>{\phi_2(U,V\otimes
W)}\\ F((U\otimes V)\otimes W)&\rTo_{F(a_{ U,V, W})}&F(U\otimes
(V\otimes W))\\ \end{diagram}
\end{equation*}
\end{claim}

\begin{examples}
\label{ex:categorii} a) The category $\mathfrak{M}_K$ of all
modules over a commutative ring $K,$ is a monoidal category with
the tensor product of $K$\lin{}modules, that will be denoted by
$\otimes _{K}.$

b) Suppose that $H$ is a Hopf algebra over a commutative ring $K.$
The category $\mathfrak{M}_{H}$ of right $H$\lin{}modules is a
monoidal category with respect to the tensor product defined as
follows. For any right $H$\lin{}modules $M$ and $N,$ let $M\otimes
N$ be $M\otimes _{K}N,$ regarded as an right $H$\lin{}module with
the structure:
\begin{equation*}
(m\otimes n)h:=\sum mh_{(1)}\otimes nh_{(2)},\text{ }\forall m\in
M,\;\forall n\in N,\;\forall h\in H
\end{equation*}
where $\Delta _{H}(h)=\sum h_{(1)}\otimes h_{(2)}$ is the
$\Sigma$\lin{}notation
that we use for the comultiplication of $H$. The unit object in $\mathfrak{M}%
_{H}$ is $K,$ which is a right $H$\lin{}module via $\varepsilon
_{H},$ the counit of $H.$

c) The category $\mathfrak{M}^{H},$ of right $H$\lin{}comodules,
is a monoidal category. The structures on the tensor product of
two bicomodules are obtained by duality from the previous example.

d) Suppose that $B$ is an arbitrary associative ring with unity.
The category ${}_B\mathfrak{M}_{B}$ of all $(B,B)$\lin{}bimodules
is a monoidal
category with the tensor product $\otimes _{B}$ and unit object $\mathbf{1}%
:=B.$
\end{examples}

\begin{claim}
\label{cl:Algebra}Following \cite[Definition 9.2.11]{Maj2}, let us
recall the
definition of associative algebras in a monoidal category $(\mathcal{M}%
,\otimes ,1,a,l,r).$ Let $A$ be an object in $\mathcal{M.}$ Suppose that $%
m:A\otimes A\rightarrow A$ and $u:\mathbf{1}\rightarrow A$ are morphisms in $%
\mathcal{M}$. If $m$ and $u$ obey the \emph{associativity} and
\emph{unity} axioms:
\begin{equation*}
\begin{diagram}[h=2em,w=1.5em] A\otimes (A\otimes A)&\rTo^{A\otimes m}& A\otimes A\\
\uTo<{a_{A,A,A}}& & \\ (A\otimes A)\otimes A& &\dTo>m \\
\dTo<{m\otimes A}&
& \\ A\otimes A&\rTo_m&A\\ \end{diagram}\hspace*{1.5cm}%
\begin{diagram}[h=2.em,w=1.5em] \mathbf{1}\otimes A&\rTo^{l_A}&A& \lTo^{r_A}&A\otimes
\mathbf{1}\\ \dTo<{u\otimes A}& &\dEqto& &\dTo>{A\otimes u} \\
A\otimes A&\rTo_m&A&\lTo_m&A\otimes A \\ \end{diagram}
\end{equation*}
we say that $(A,m,u)$ is an (associative) algebra with
\emph{multiplication}
$m$ and \emph{unit} $u$ in $\mathcal{M}$. As we explained in (\ref{cl:CohThm}%
), we can omit the maps $a,$ $l$ and $r,$ so we shall draw these
diagrams in a more simple way as follows.
\begin{equation*}
\begin{diagram}[h=2em,w=1.5em]
A\otimes A\otimes A&\rTo^{A\otimes m}& A\otimes A\\
\dTo<{m\otimes A}& & \dTo>m\\ A\otimes A&\rTo_m&A\\ \end{diagram}%
\hspace*{1.5cm}
\begin{diagram}[h=2em,w=1.5em]\mathbf{1}\otimes
A&\rEqto&A& \lEqto&A\otimes \mathbf{1}\\ \dTo<{u\otimes A}&
&\dEqto& &\dTo>{A\otimes u}
\\ A\otimes A&\rTo_m&A&\lTo_m&A\otimes A \\ \end{diagram}
\end{equation*}
\end{claim}

\begin{examples}
\label{ex:AlgTensCat} a) An algebra in $(\mathfrak{M}_K,$ $\otimes
_{K},K)$
is an unitary associative ring $A$ together with a ring morphism $%
i:K\rightarrow A$ such that the image of $i$ is included in the center of $%
A. $ We recognize the usual definition of algebras over a
commutative ring.

b) Let $H$ be a Hopf algebra. An algebra in the category
$\mathfrak{M}_{H}$ is an associative algebra $A,$ in the usual
sense, which is a right $H$\lin{}module such that
\begin{eqnarray*}
(xy)h &=&\sum xh_{(1)}\otimes yh_{(2)},\;\forall x,y\in
A,\;\forall h\in H,
\\
1_Ah &=&\varepsilon _{H}(h)1,\;\forall h\in H.
\end{eqnarray*}
We recognize the definition of $H$\lin{}\emph{module algebras}
\cite[Definition 4.1.1]{Mo}, sometimes called
$H$\lin{}{}\emph{differential algebras} (see for example
\cite{BDK}).

c) An algebra in $\mathfrak{M}^{H}$ is a right
$H$\lin{}{}\emph{comodule algebra}, see \cite[Definition
4.1.2]{Mo}. Recall that this means an algebra $A$ which is a right
$H$\lin{}{}comodule such that
\begin{eqnarray*}
\rho (xy) &=&\sum x_{\left\langle 0\right\rangle }y_{\left\langle
0\right\rangle }\otimes x_{\left\langle 1\right\rangle
}y_{\left\langle 1\right\rangle },\;\forall x,y\in A, \\ \rho (1)
&=&1\otimes 1.
\end{eqnarray*}

d) A bimodule over $B,$ say $A,$ is an algebra in $({}_B\mathfrak{M}%
_{B},\otimes _{B},B)$ iff $A$ is an associative ring with unity
$1_{A}$ such that $1_{A}\in \{a\in A\mid ba=ab,$ $\forall b\in
B\}.$ This set will be denoted $A^{B}$ (more generally, if $M\in
{}_B\mathfrak{M}_{B}$ then $M^{B}$ will denote the set of all
$m\in M$ such that $bm=mb,$ $\forall b\in B$). For example, any
morphism of associative rings $i:B\rightarrow A$ gives an algebra
in ${}_B\mathfrak{M}_{B}$ where $A$ is a $(B,B)$\lin{}{}bimodule
with the restriction of scalars via $i.$
\end{examples}

\begin{claim}
\label{cl:Module}Now we are going to define the representations of
algebras
in monoidal categories. We shall proceed as in the case of algebras in $%
\mathfrak{M}_K.$ Let us assume that $(A,m,u)$ is an algebra in the
monoidal category $(\mathcal{M},\otimes ,\mathbf{1}).$ By a left
$A$\lin{}{}\emph{module} we mean an object $M\in \mathcal{M}$
together with a morphism $\mu :A\otimes M\rightarrow M$ such that
\begin{equation*}
\begin{diagram}[h=2.em,w=1.5em]
 A\otimes A \otimes M&\rTo^{A\otimes \mu}& A\otimes M\\
\dTo<{m\otimes M}& & \dTo>{\mu}\\ A\otimes M&\rTo_{\mu}&M\\
\end{diagram}%
\hspace*{1.5cm}
\begin{diagram}[h=2.em,w=1.5em] \mathbf{1}\otimes M&\rEqto&M\\
\dTo<{u\otimes M}& &\dEqto\\ A\otimes M&\rTo_{\mu}&M \\
\end{diagram}
\end{equation*}
are commutative. If $(M,\mu )$ and $(N,\nu )$ are two left
$A$\lin{}{}modules a
morphism of modules from $M$ to $N$ is a morphism $f:M\rightarrow N$ in $%
\mathcal{M}$ such that $\nu (A\otimes f)=f\mu .$ The category of
left $A$\lin{}{}modules will be denoted by $_{A}\mathcal{M}$ Let
us remark that $_{A} \mathcal{M}$ is an abelian category if
$\mathcal{M}$ is so.

Similarly, we construct the category of right modules
$\mathcal{M}_{A}.$ Combining left and right modules we get
$(A,A)$\lin{}{}bimodules. More precisely, an
$(A,A)$\lin{}{}bimodule is an object in $\mathcal{M}$ together
with two maps, $\mu _{l}:A\otimes M\rightarrow M$ and $\mu
_{r}:M\otimes A\rightarrow M,$ such
that $(M,\mu _{l})\in $ ${}_{A}\mathcal{M}$ and $(M,\mu _{r})\in \mathcal{M}%
_{A}$ and the structures are compatible, that is the following
diagram is commutative.
\begin{equation*}
\begin{diagram}[h=2.em,w=1.5em] A\otimes M \otimes A&\rTo^{A\otimes \mu_r}& A\otimes M\\
\dTo<{\mu_l\otimes A}& & \dTo>{\mu_l}\\ M\otimes
A&\rTo_{\mu_r}&M\\
\end{diagram}
\end{equation*}
A morphism $f:M\rightarrow N$ between two bimodules is a morphism in $%
\mathfrak{M}$ which is both a morphism of left ant right modules.
For the category of $(A,A)$\lin{}{}bimodules we shall use the
notation ${}_{A}\mathcal{M}_{A}.$ Of course, if $\mathcal{M}$ is
abelian then ${}_{A}\mathcal{M}_{A}$ is abelian too.
\end{claim}

\begin{examples}
\label{ex:Modules}a) $A$ always is an $(A,A)$\lin{}{}bimodule,
having both left and right module structures defined by the
multiplication $m.$

b) Suppose that $(A,m,u)$ is an algebra in $(\mathcal{M},\otimes ,\mathbf{1}%
).$ Then $A\otimes X\in {}_{A}\mathcal{M},$ for any $X\in
\mathcal{M}$, where the left structure is given by $\mu :=m\otimes
X.$ Thus we have a functor $_{A}F:\mathcal{M}\rightarrow
{}_{A}\mathcal{M}$ , which is defined by $_{A}F(X)=A\otimes X$ and
$_{A}F(f)=A\otimes f.$

Similarly $X\otimes A$ is a right $A$\lin{}{}module, so we obtain a functor: $%
F_{A}:\mathcal{M}\rightarrow {}\mathcal{M}_{A}$ given by
$F_{A}(X)=X\otimes A $ and $F_{A}(f)=f\otimes A.$

c) Let $A$ be as above, and let $M\in {}_{A}\mathcal{M}.$ Then
$M\otimes A$ is a right $A$\lin{}{}module as in the previous
example, and is a left $A$\lin{}{}module via $\nu =\mu \otimes A.$
These two structures are compatible, defining an
$(A,A)$\lin{}{}bimodule on $M\otimes A.$ Similarly, if $M\in
\mathcal{M}_{A}$ then $ A\otimes M$ is an $(A,A)$\lin{}{}bimodule.

In particular $(A\otimes X)\otimes A$ is an
$(A,A)$\lin{}{}bimodule and $_{A}F_{A}: \mathcal{M}\rightarrow
{}_{A}\mathcal{M}_{A}$, $_{A}F_{A}(X)=(A\otimes X)\otimes A$ and
$_{A}F_{A}(f)=(A\otimes f)\otimes A$ is a functor.

Analogously, $A\otimes (X\otimes A)$ can be regarded as an
$(A,A)$\lin{}{}bimodule, and one can easily prove that
$a_{A,X,A}:(A\otimes X)\otimes A\rightarrow A\otimes (X\otimes A)$
is a functorial isomorphism of bimodules$.$
\end{examples}

\begin{proposition}
\label{pr:Adjunction}a) \ $_{A}F$ is a left adjoint of $_{A}U:{}_{A}%
\mathcal{M\rightarrow M},$ the functor that ``forgets'' the module
structure.

b) \ $F_{A}$ is a left adjoint of $U_{A}:{}\mathcal{M}_{A}%
\mathcal{\rightarrow M},$ the functor that ``forgets'' the module
structure.

c) \ $_{A}F_{A}$ is a left adjoint of $_{A}U_{A}:{}_{A}\mathcal{M}%
_{A}\rightarrow \mathcal{M},$ the functor that ``forgets'' the
bimodule structure.
\end{proposition}

\begin{proof}
a) To prove that $_{A}F$ is a left adjoint of $_{A}U:{}_{A}\mathcal{M}%
\rightarrow {}\mathcal{M}$ we need functorial morphisms:
\begin{equation*}
{}_{A}\mathcal{M}(A\otimes X,M)\pile{\rTo^{\phi_l(X,M)}\\ \lTo_{\psi_l(X,M)}}%
\mathcal{M}(X,M),
\end{equation*}
that are inverses each other. We define $\phi
_{l}(X,M)(f):=f(u\otimes X)l_{X}^{-1}$ and $\psi _{l}(X,M)(g):=\mu
(A\otimes g),$ where $\mu $ is the module structure of $M.$ It is
easy to prove that $\psi _{l}(X,M)(g)$ is a morphism of left
modules, and that $\psi _{l}(X,M)$ is the inverse of $\phi
_{l}(X,M).$

b) The isomorphisms
\begin{equation*}
\mathcal{M}_{A}(X\otimes A,M)\pile{\rTo^{\phi_r(X,M)}\\
\lTo_{\psi_r(X,M)}}\mathcal{M}(X,M)
\end{equation*}
are now given by $\phi _{r}(X,M)(f):=f(X\otimes u)r_{X}$ and $\psi
_{l}(X,M)(g):=\mu (g\otimes A),$ where $\mu $ is the module
structure of $M.$

c) The isomorphisms $_{A}\mathcal{M}_{A}((A\otimes X)\otimes A,M)%
\pile{\rTo^{\phi(X,M)}\\ \lTo_{\psi(X,M)}}\mathcal{M}(X,M)$ are
obtained by combining the isomorphisms constructed above: $\phi
(X,M)=\phi _{l}(X,M)\phi _{r}(A\otimes X,M),$ and similarly for
$\psi (X,M).$ For future references, we explicitly write them
down:
\begin{eqnarray}
\phi (X,M)(f) &=&f[(A\otimes X)\otimes u]r_{A\otimes
X}^{-1}(u\otimes X)l_{X}^{-1},  \label{ec:Fi(X,M)} \\ \psi
(X,M)(g) &=&\mu _{r}(\mu _{l}\otimes A)[(A\otimes g)\otimes A],
\label{ec:Psi(X,M)}
\end{eqnarray}
where $\mu _{r}$ and $\mu _{l}$ give respectively the right and
left $A$\lin{}{}module structures of $M.$
\end{proof}

\begin{corollary}
The functors $_{A}F,$ $F_{A}$ and $_{A}F_{A}$ are right exact.
\end{corollary}

\begin{claim}
\label{cl:SimplObj}Since one of our main goals is to investigate
the relative derived functors of ${}_{A}\mathcal{M}_{A}(A,-),$
with respect to a
certain projective class of epimorphisms in $\mathcal{M},$ we need ${}_{A}%
\mathcal{M}_{A}$ to be an abelian category. One can prove easily that ${}_{A}%
\mathcal{M}_{A}$ is abelian, if we assume that $\mathcal{M}$ is
so. Thus,
from now on, we shall assume that $\mathcal{M}$ \emph{is an abelian category}%
.

In order to produce projective resolutions we will apply the
machinery of
bar resolutions, see \cite[Chapter 8.6]{We}. The pair of adjoint functors $%
(F_{A},U_{A})$ defines a cotriple $(\perp _{A},\varepsilon
_{A},\delta _{A})$ on $\mathcal{M}_{A}\mathfrak{.}$ The functor
$\bot _{A}$ is defined by $\bot _{A}:=F_{A}U_{A}.$ The functorial
morphism $\varepsilon _{A} $ is the counit
of the adjunction, and $\delta _{A}(M):=F_{A}(\eta _{U_{A}(M)}),$ where $%
\eta $ is the unit of the adjunction. A quick computation shows us
that, for
any $(M,\mu )\in \mathcal{M}_{A},$ we have $\varepsilon _{A}(M)=\mu ,$ and $%
\delta _{A}(M)=[(M\otimes u)r_{M}^{-1}]\otimes A.$

By following the construction in \cite[8.6.4]{We}, for any $(M,\mu
)\in
\mathcal{M}_{A}$ we obtain a simplicial object $(\beta_{\bullet}(A,M),%
\partial _{\bullet },\sigma _{\bullet })$, where $\beta_{n}(A,M)=\bot
_{A}^{n+1}M$. Its face and degeneracy operators are:
\begin{equation*}
\partial _{i}=\bot _{A}^{i}(\varepsilon _{A}(\bot _{A}^{n-i}M))\quad\text{and%
}\quad\sigma _{i}=\bot _{A}^{i}(\delta _{A}(\bot _{A}^{n-i}(M)).
\end{equation*}
For $i<n$ the module structure on $\bot _{A}^{n-i}(M)$ is defined
by $(\bot _{A}^{n-i}M\otimes m)a_{\bot _{A}^{n-i}M,A,A},$ so we
have:
\begin{eqnarray}
\partial _{i} &=&\left\{
\begin{array}{ll}
\bot _{A}^{i}((\bot _{A}^{n-i-1}M\otimes m)a_{\bot
_{A}^{n-i-1}M,A,A}), & \text{for }0\leq i<n; \\ \perp _{A}^{n}(\mu
), & \text{for }i=n.
\end{array}
\right.  \label{ec:Face} \\ \sigma _{i} &=&\bot _{A}^{i}((\bot
_{A}^{n-i}M\otimes u)r_{\bot _{A}^{n-i}M}^{-1}\otimes A).
\label{ec:Degeneracy}
\end{eqnarray}
By \cite[Proposition 8.6.10]{We}, the augmented simplicial object $%
U(\beta_{\bullet}(A,M))\overset{U(\mu )}{\rightarrow }U(M)$ is
aspherical, so the associated augmented chain complex is exact.
Since $U$ is faithfully exact, it results that
$\beta_{\bullet}(A,M)\overset{\mu }{\rightarrow }M$ is aspherical
too. Its associated exact sequence $\beta_{\bullet}(A,M)$ will be
called, as in the classical case, \emph{the} \emph{bar resolution
of} $M$ \emph{in} $\mathcal{M}_{A}$. Our next aim is to give a new
interpretation of $\beta (A,A)$ in terms of
$\mathcal{E}$\lin{}{}projective resolutions, where $ \mathcal{E}$
is an appropriate class of projective epimorphisms.
\end{claim}

\begin{lemma}
\label{le:Beta=Bimodul}If $M$ is an $(A,A)$\lin{}{}bimodule in
$\mathcal{M}$ then $ \beta_{\bullet} (A,M)$ is an exact complex in
${_{A}\mathcal{M}_{A}}.$
\end{lemma}

\begin{proof}
Since $M\in {_{A}\mathcal{M}_{A}}$ it follows that $U_{A}(M)$ is a
left module, so $F_{A}U_{A}(M)$ is an $(A,A)$\lin{}{}bimodule,
with the structures as in example \ref{ex:Modules}c). By induction
$\bot ^{n}M$ is an $(A,A)$\lin{}{}bimodule for any $n\geq 0$. If
remains to show that the differential maps are morphisms of
$(A,A)$\lin{}{}bimodules. But $\partial _{i}:\bot
^{n+1}M\rightarrow \bot ^{n}M$ is given by $\partial _{i}=\bot
_{A}^{i}(\varepsilon _{A}(\bot _{A}^{n-i}M))$ and in our case
$\varepsilon _{A}(\bot _{A}^{n-i}M)$ defines the right structure
on $\bot _{A}^{n-i}M$. Obviously for any bimodule $N$ the maps
$\mu _{r}:$ $N\otimes A\rightarrow N$ and $\mu _{l}:A\otimes
N\rightarrow N$ defining the module structures are
morphisms of left, respectively right, $A$\lin{}{}modules. Moreover, if $%
f:N\rightarrow P$ is a morphism of bimodules, then $\perp
_{A}^{n}(f)$ is a morphism of bimodules too. Then $\partial_i $ is
a morphism of bimodules,
which ends the proof of the lemma because the differential maps of $%
\beta_{\bullet} (A,M)$ are defined by $d_{n}=
\sum_{i=0}^{n}(-1)^{i}\partial _{i}$.
\end{proof}

\begin{claim}\label{cl:E-Proj}
Let $\mathcal{M}$ be an abelian category and let $\mathcal{E}$ be
a class of epimorphisms in $\mathcal{M}$. We recall that an object
$P$ in $\mathcal{M}$ is called projective rel $\varepsilon $,
where $\varepsilon :X\rightarrow Y$ is an epimorphism in
$\mathcal{E}$, if $\mathcal{M}(P,\varepsilon ):\mathcal{
M}(P,X)\rightarrow \mathcal{M}(P,X)$ is surjective. $P$ is called
$\mathcal{E }$\lin{}{}projective if it is projective rel
$\varepsilon $ for every $\varepsilon $ in $\mathcal{E}$.
The\emph{\ closure} of $\mathcal{E}$ is the class $
\mathcal{C(E)}$ containing all epimorphism $\mathcal{E}$ in
$\mathcal{M}$ such that every $\mathcal{E}$\lin{}{}projective
object is also projective rel $ \varepsilon $. The class
$\mathcal{E}$ is called \emph{closed} if $\mathcal{E }$ is
$\mathcal{C(E)}$. The class $\mathcal{E}$ is called
\emph{projective} if for any object $M$ in $\mathcal{M}$ there is
an epimorphism $\varepsilon :P\rightarrow M$ in $\mathcal{E}$ such
that $P$ is $\mathcal{E}$\lin{}{}projective. Suppose now that
$\mathcal{E}$ is a closed class of epimorphisms in $\mathcal{M} $.
A morphism $f:X\rightarrow Y$ in $\mathcal{M} $ is called
$\mathcal{E}$\lin{}{}admissible if in the canonical splitting $f=i
p$ $ ,$ $i$ monic and $p$ epic, we have $p\in \mathcal{E}$.
Finally an $\mathcal{E }$\lin{}{}projective resolution of $M$ is
an exact sequence:
\begin{equation*}
\longrightarrow P_{n}\overset{d_{n}}{\longrightarrow
}P_{n-1}\overset{d_{n-1} }{\longrightarrow }\cdots \rightarrow
P_{1}\overset{d_{1}}{\longrightarrow }
P_{0}\overset{d_{0}}{\longrightarrow }M\longrightarrow 0
\end{equation*}
such that all maps are $\mathcal{E}$\lin{}{}admissible and $P_{n}$
is $\mathcal{E}$\lin{}{}projective, for every $n\geq 0$. As in the
usual case ($\mathcal{E}$ is the class of all epimorphisms) one
can show that any object in $\mathcal{M}$ has an
$\mathcal{E}$\lin{}{}projective resolution, which is unique up to
a homotopy. The theory of derived functors can be adapted to the
relative context without difficulties. For details the reader is
referred to \cite[Chapter XI] {HiSt}.

\begin{claim}
We are now going to define the projective class of epimorphisms
that we are interested in. Let $(\mathcal{M},\otimes ,1)$ be a
monoidal category and let $(A,m,u)$ be an algebra in
$\mathcal{M}$. In the abelian category ${_{A} \mathcal{M}_{A}}$
(recall that we always assume that $\mathcal{M}$ is abelian) we
consider the class $\mathcal{E}$ of all epimorphisms that have a
section in $\mathcal{M}$. To prove that $\mathcal{E}$ is
projective we note that
$\mathcal{E}={_{A}}U{_{A}^{-1}}(\mathcal{E}_{0})$, where $
(_{A}F_{A},{}_{A}U_{A})$ is the pair of adjoint functors from
Proposition \ref{pr:Adjunction}, and $\mathcal{E}_{0}$ is the
class of all epimorphism in $\mathcal{M}$ that have a section in
$\mathcal{M}$. Obviously $\mathcal{E} _{0}$ is projective since
any object in $\mathcal{M}$ is
$\mathcal{E}_{0}$\lin{}{}projective. Theorem IX.4.1 in \cite{HiSt}
reads in our situation as follows.
\end{claim}
\end{claim}

\begin{proposition}
\label{cl:Beta=E-projective}The class $\mathcal{E}$ of all epimorphisms in ${%
_{A}\mathcal{M}_{A}}$ that split in $\mathcal{M}$ is projective.
The objects ${_{A}F_{A}}P$, where $P \in \mathcal{M}$, are
$\mathcal{E}$\lin{}{}projective and are sufficient for
$\mathcal{E}$\lin{}{}presenting objects of $_{A}\mathcal{M}_{A}$,
so that the $\mathcal{E}$\lin{}{}projectives are precisely the
direct summands of objects ${_{A}F_{A}}P$.
\end{proposition}

\begin{claim}
Following \cite{HiSt}, in view of foregoing proposition, we can
now consider, for every $M \in {_A\mathcal{M}_A}$, the right
$\mathcal{E}$\lin{}{}derived functors
$\mathbf{R}^\bullet_{\mathcal{E}} F_M$ of
$F_M:={}_A\mathcal{M}_A(-,M)$. Then, for every $M,N \in
{_A\mathcal{M}_A}$, we set:
\begin{equation}
\mathbf{Ext}^\bullet_{\mathcal{E}}
(N,M)=\mathbf{R}^\bullet_{\mathcal{E}} F_M(N).
\end{equation}
\end{claim}
\noindent The following well known result can be proved as in the
non-relative case.
\begin{proposition}\label{te:verysep}
Let $(A,m,u)$ be an algebra in a monoidal category $(\mathcal{M},\otimes ,%
\mathbf{1})$ and let $N\in ${{$_{A}\mathcal{M}_{A}$}}. The
following assertions are equivalent:

(a) $N$ is $\mathcal{E}$\lin{}{}projective.

(b) $\mathbf{Ext}^1_{\mathcal{E}} (N,M)=0,$ for all $M\in
${{$_{A}\mathcal{M}_{A}$}}.

(c) $\mathbf{Ext}^n_{\mathcal{E}} (N,M)=0,$ for all $M\in
${{$_{A}\mathcal{M}_{A}$}}, for all $n>0.$
\end{proposition}

\begin{claim}
We now want to prove that the bar resolution
$\beta_{\bullet}(A,A)$ is made of $\mathcal{E}$\lin{}{}projective
modules. To this aim, let us prove that we have a canonical
isomorphism of bimodules:
\begin{equation*}
\varphi _{n}:\perp ^{n+1}A\rightarrow (A\otimes \perp ^{n-1}A)\otimes A,%
\text{ for all }n\geq 1.
\end{equation*}
Indeed, if $n=1$ then $a_{A,A,A}:(A\otimes A)\otimes A\rightarrow
A\otimes (A\otimes A)$ is an isomorphism of bimodules, as we
noticed in example \ref {ex:Modules}(c). So we can take $\varphi
_{1}:=a_{A,A,A}.$ Let us assume that we have constructed $\varphi
_{1},\cdots ,\varphi _{n}.$ Then:
\begin{equation*}
\perp _{A}^{n+2}A=\perp _{A}(\perp _{A}^{n+1}A) \overset{\perp
_{A}(\varphi
_{n})}{\simeq }\perp _{A}((A\otimes \perp _{A}^{n-1}A)\otimes A) \overset{%
\perp _{A}(a_{A,\perp ^{n-1}A,A})}{\simeq }\perp _{A}(A\otimes
(\perp _{A}^{n-1}A\otimes A)).
\end{equation*}

As $\perp _{A}(A\otimes (\perp _{A}^{n-1}A\otimes A)) =\perp
_{A}(A\otimes \perp _{A}^{n}A)=(A\otimes \perp _{A}^{n}A)\otimes
A$ we can take by definition $\varphi _{n+1}:=\perp
_{A}(a_{A,\perp _{A}^{n-1}A,A})\perp _{A}(\varphi _{n}).$ Since
$\varphi _{n}$ and $a_{A,\perp _{A}^{n-1}A,A}$ are morphisms of
$(A,A)$\lin{}{}bimodules we deduce that $\perp _{A}(\varphi _{n})$
and $\perp _{A}(a_{A,\perp _{A}^{n-1}A,A})$ are so. Hence $\varphi
_{n+1}$ is an isomorphism of bimodules. Now we can prove the
following lemma.
\end{claim}

\begin{lemma}
\label{le:Beta=E-projective}$\beta _n(A,A)$ is
$\mathcal{E}$\lin{}{}projective for every $n\in \mathbb{N}$.
\end{lemma}

\begin{proof}
By (\ref{cl:Beta=E-projective}) we know that $\beta
_{n}(A,A)\simeq
(A\otimes \perp ^{n-1}A)\otimes A={_{A}}F{_{A}}(\perp ^{n-1}A)$. Since $%
\mathcal{E} _{0} $ is the class of all splitting epimorphisms, it
follows that any object in $\mathcal{M}$ is
$\mathcal{E}_{0}$\lin{}{}projective. By
\cite[Theorem IX.4.1]{HiSt} we deduce that $_{A}F_{A}(X)$ is $\mathcal{E}$%
\lin{}{}projective (see also the definition of $\mathcal{E}$) for
every $X\in \mathcal{M}$. In particular $\beta _{n}(A,A)\simeq
{_{A}}F{_{A}}(\perp ^{n-1}A)$ is $\mathcal{\
E}$\lin{}{}projective.
\end{proof}

\begin{theorem}
\label{th:BarResolution}$\beta_{\bullet}(A,A)$ is an
$\mathcal{E}$\lin{}{}projective resolution of $A$ in
$_A\mathcal{M}_A$.
\end{theorem}

\begin{proof}
We already know that $\beta_{\bullet}(A,A)$ is an exact sequence in ${}_{A}%
\mathcal{M}_{A}$, see Lemma \ref{le:Beta=Bimodul}. By Lemma \ref
{le:Beta=E-projective} it results that $\beta _{n}(A,A)$ is $\mathcal{E-}$%
projective for every $n\in \mathbb{N}.$ It remains to show that
the differential maps of $\beta_{\bullet}(A,A)$ are
$\mathcal{E}$\lin{}{}admissible. By \cite[Proposition 8.6.10]{We}
it follows that the augmented simplicial object
$U_{A}(\beta_{\bullet}(A,A))\overset{U_{A}(m)}{\longrightarrow
}U_{A}(A)$ constructed in \ref{cl:SimplObj} is contractible in
$\mathcal{M}$. Here, $m$ denotes the multiplication of $A$. It
follows that $\beta_{\bullet}(A,A)$ is
split exact in $\mathcal{M},$ see \cite[Exercise 8.4.6]{We}$.$ Let $%
s_{n}:\beta_n(A,A)\rightarrow \beta_{n+1}(A,A)$ be a morphism such that $%
d_{n}=d_{n}s_{n-1}d_{n},$ where $d_{n}:\beta_{n}(A,A)\rightarrow
\beta_{n-1}(A,A)$ are the differentials of $\beta_{\bullet}(A,A).$ If $%
d_{n}=i_{n}p_{n}$ is the canonical decomposition in $\mathcal{M},$ with $%
p_{n}$ an epimorphism and $i_{n}$ a monomorphism, then $p_{n}(s_{n-1}i_{n})=%
\mathrm{Id}.$ Thus $d_{n}$ splits in $\mathcal{M},$ so $d_{n}$ is
admissible for any $n\in \mathbb{N}.$
\end{proof}

\begin{definition}
Let $(\mathcal{M},\otimes,\mathbf{1})$ be a monoidal category. Suppose that $%
A$ is an algebra in $\mathcal{M}$ and that $M$ is an
$(A,A)$\lin{}{}bimodule. The \emph{Hochschild cohomology} of $A$
with coefficients in $M$ is:
\begin{equation}
\mathbf{H}^\bullet(A,M)=\mathbf{Ext}^\bullet_{\mathcal{E}} (A,M).
\end{equation}
\end{definition}

\begin{claim}
In order to compute $\mathbf{H}^\bullet(A,M)$ we can apply the functor ${}_A%
\mathcal{M}_A(-,M)$ to $\beta_{\bullet}(A,A)$, the bar resolution of $A$ in $%
{}_A\mathcal{M}_A$. We obtain \emph{the standard complex}:
\begin{equation*}
0 \longrightarrow \mathcal{M}(\mathbf{1},M)\overset{b^0}{%
\longrightarrow }\mathcal{M}({A},M)\overset{b^1}{\longrightarrow
}\mathcal{M}({A\otimes A},M)\overset{b^2}{\longrightarrow
}\mathcal{M}({A\otimes A\otimes A},M)\overset{b^3}{\longrightarrow
} \cdots
\end{equation*}
where for $n\in\{0,1,2\}$ the differentials $b^n$ are given by:
\begin{eqnarray*}
b^0(f)&=&\mu_l(A\otimes f)r_A^{-1}-\mu_r(f\otimes A)l_A^{-1}; \\
b^1(f)&=&\mu_l(A\otimes f)-f m+\mu_r(f\otimes A); \\
b^2(f)&=&\mu_l(A\otimes f) a_{A,A,A}-f (m\otimes A)+f (A\otimes m)
a_{A,A,A}-\mu_r(f\otimes A).
\end{eqnarray*}
We will not write down the formula for $b^n$ in general, because
we shall need an explicit computation only in degree
$n\in\{0,1,2\}$. Note that if we
omit the associativity constraint $a_{A,A,A}$ then the formulas for $b^0$, $%
b^1$ and $b^2$ becomes the usual ones, well\lin{}{}known from the
case $\mathcal{M}=
\mathfrak{M}_K$. This observation holds true in general, for all $b^n$, $n\in%
\mathbb{N}$.
\end{claim}

\section{Hochschild dimension of algebras in monoidal categories}

In this section we define the Hochschild dimension of an algebra
$A$ in a monoidal category $(\mathcal{M},\otimes,\mathbf{1})$.
Then we will give characterizations of algebras of Hochschild
dimension less than or equal to 1. Our next goal will be to study
separable algebras in monoidal categories and to prove a
Wedderburn-Malcev type theorem. This last theorem will be used in
the next sections to study the structure of Hopf algebras having a
``nice'' radical (Section 3) or coradical (Section 5).

\begin{definition}
An algebra $(A,m,u)$ in a monoidal category $(\mathcal{M},\otimes
,\mathbf{1} )$ is called \emph{separable} if the multiplication
$m:A\otimes A\rightarrow A$ has a section in the category of
$(A,A)$\lin{}{}bimodules ${{_{A}\mathcal{M}_{A} }}$.
\end{definition}

\begin{remark}
The multiplication $m$ always has a section in $_{A}\mathcal{M}$
and in $ \mathcal{M}_{A}$, namely $A\otimes u$ and respectively
$u\otimes A$.
\end{remark}

\begin{proposition}\label{prop:sep}
$(A,m,u)$ is separable iff $A$ is $\mathcal{E}$\lin{}{}projective.
\end{proposition}

\begin{proof}
Recall that $\mathcal{E}$ is the projective class of all epimorphism in ${{%
_{A}\mathcal{M}_{A}}}$ that have a section in $\mathcal{M}$. Therefore an $%
(A,A)$\lin{}{}bimodule $P$ is $\mathcal{E}$\lin{}{}projective if
there is an object $X$ in $\mathcal{M}$ and an epimorphism $\pi
:A\otimes (X\otimes A)\rightarrow P$
in ${{_{A}\mathcal{M}_{A}}}$ that splits in ${{_{A}\mathcal{M}_{A}}}$. Thus $%
A$ is $\mathcal{E}$\lin{}{}projective if $m:A\otimes A\rightarrow
A$ has a section
in ${{_{A}\mathcal{M}_{A}}}$, since $A\otimes A\simeq A\otimes (\mathbf{1}%
\otimes A)$. Conversely, if $A$ is
$\mathcal{E}$\lin{}{}projective, then $m$ has a
section in ${{_{A}\mathcal{M}_{A}}}$ since $m$ is an epimorphism in $%
\mathcal{E}$.
\end{proof}

\begin{theorem}\label{te:sep}
Let $(A,m,u)$ be an algebra in a monoidal category $(\mathcal{M},\otimes ,%
\mathbf{1})$. The following assertions are equivalent:

(a) $A$ is separable.

(b) $\mathbf{H}^{1}(A,M)=0,$ for all $M\in
${{$_{A}\mathcal{M}_{A}$}}.

(c) $\mathbf{H}^{n}(A,M)=0,$ for all $M\in
${{$_{A}\mathcal{M}_{A}$}}, for all $n>0.$

\end{theorem}

\begin{proof} Follows by Proposition \ref{te:verysep}, in the case
when $N=A$, and by Proposition \ref{prop:sep}.
\end{proof}

\begin{definition}
The \emph{Hochschild dimension} of an algebra $A$ in the monoidal category $%
\mathcal{M} $ is the smallest $n\in \mathbb{N}$ (if it exists) such that $%
\mathbf{H}^{n+1}(A,M)=0,\forall M\in {{_{A}\mathcal{M}_{A}}}$. If
such an $n$ does not exist, we will say that the Hochschild
dimension of $A$ is
infinite. We shall denote the Hochschild dimension of $A$ by $\mathrm{%
\mathrm{Hdim}}(A)$.
\end{definition}

\begin{corollary}
An algebra $(A,m,u)$ in $(\mathcal{M},\otimes ,\mathbf{1})$ is
separable iff $\mathrm{\mathrm{Hdim}}(A)=0$.
\end{corollary}

\noindent For the characterization of algebras of Hochschild
dimension $1$ we need the interpretation of $\mathbf{H}^{2}(A,-)$
in terms of algebra extensions. First some definitions.

\begin{definition}
Let $A$ and $B$ be two algebras in a monoidal category $(
\mathcal{M},\otimes ,\mathbf{1})$. A morphism $\sigma
:B\rightarrow A$ in $ \mathcal{M}$ is called \emph{unital} if
$\sigma u_{B}=u_{A}$, where $u_{A}$ and $u_{B}$ are the units of
$A$ and $B,$ respectively. Moreover, if $f:A\ra B$ is a morphism
of algebras in $\calM$ we shall say that $\sigma$ is an
\emph{unital section} of $f$ if $f \sigma =\mathrm{ Id}_{B}$ and
$\sigma$ is an unital morphism.
\end{definition}

\begin{definitions}
a) \ An ideal of an algebra $(A,m,u)$ in $(\mathcal{M},\otimes
,\mathbf{1})$ is a pair $(I,i)$ such that $I$ is an
$(A,A)$\lin{}{}bimodule and $i:I\rightarrow A $ is a monomorphism
of $(A,A)$\lin{}{}bimodules.

b) If $(I,i)$ is an ideal in $A$ and $n\geq 2,$ we define $I^{n}$
to be the image of $m_{n} {i^{\otimes n}}$ where $m_{n}:A^{\otimes
n}\rightarrow A$ is the $n^{\text{th}}$ iterated multiplication of
$A$ $(m_{2}:=m)$.

c) An ideal $(I,i)$ is called nilpotent if there is $n\geq 2$ such that $%
I^{n}=0$ (equivalently $m_{n} {i^{\otimes n}}=0$).
\end{definitions}

\begin{remarks}
a) If $(I,i)$ is a ideal then there is a unique algebra structure on $A/I=%
\mathrm{\mathrm{Coker}}$ $i$ such that the canonical morphism
$A\rightarrow A/I$ is an algebra map. \newline
b) There is a canonical monomorphism $i_{n}:I^{n}\rightarrow A$ in ${{_{A}%
\mathcal{M}_{A}}}$. Thus $(I_{n},i_{n})$ is an ideal of $A$. The
map $i_{n}$ factors through a morphism of bimodules
$i_{n}':I^{n}\rightarrow I$ $($i.e. $i i_{n}'=i_{n})$. Moreover,
for every $m\geq n\geq 2$,
there is a morphism of bimodules $i_{n,m}:I^{m}\rightarrow I^{n}$ such that $%
i_{n} i_{n,m}=i_{m}$.
\end{remarks}

\begin{lemma}
\label{X lem 1.5.8}Let $\left( A,m,u\right) $ and $\left(
E,m_{E},u_{E}\right) $ be algebras. Let $\pi :E\rightarrow A$ be a
morphism
of algebras in $(\mathcal{M},\otimes ,\mathbf{1})$ that has a section $%
\sigma $ in $\mathcal{M}$. Assume that $(\mathrm{Ker\,} \pi
)^{2}=0$.

a) If $\sigma :A\rightarrow E$ is a section of $\pi $ then:$\quad$
$m_{E} \left( \sigma u\otimes \sigma u\right)
l_{\mathbf{1}}^{-1}=2\sigma u-u_{E}$.

b) The morphism $\pi $ admits a unital section.

c) $\mathrm{Ker\,} \pi $ has a natural structure of
$(A,A)$\lin{}{}bimodule given
by $\mu _{l}:A\otimes {\mathrm{Ker\,} }\pi $ $\rightarrow {\mathrm{Ker\,} }%
\pi $ and $\mu _{r}:{\mathrm{Ker\,} }\pi $ $\otimes A\rightarrow {\mathrm{%
Ker\,} }\pi $ that are uniquely defined by:
\begin{eqnarray}
i \mu _{l} &=&m_{E} (E\otimes i) (\sigma \otimes {\mathrm{Ker\,}
}\pi ), \label{eq:LeftMod} \\ i \mu _{r} &=&m_{E} (i\otimes E)
({\mathrm{Ker\,} }\pi \otimes \sigma ), \label{eq:RightMod}
\end{eqnarray}
where $i:{\mathrm{Ker\,} }\pi \rightarrow E$ is the canonical
inclusion. The
morphisms $\mu _{l}$ and $\mu _{r}$ do not depend on the choice of the section $%
\sigma $.
\end{lemma}

\begin{proof}
a) \ The relation $\pi \left( \sigma u-u_{E}\right) =0$ tell us
that there
exists a unique morphism $\lambda :\mathbf{1}\rightarrow E$ such that $%
\sigma u-u_{E}=i \lambda $. On the other hand, $(\mathrm{Ker\,}\pi)^2=0$ so $%
m_E(i\otimes i)=0$. Thus:
\begin{equation*}
m_{E} \left[ (\sigma u-u_{E})\otimes (\sigma u-u_{E})\right]  l_{\mathbf{1}%
}^{-1}=0.
\end{equation*}
By expanding this relation and using that $u_E$ is the unit of $E$
we obtain the formula from the first part of the lemma.

b) \ The unital section is given by $\sigma ':=2\sigma-m_{E}
\left( \sigma \otimes \sigma \right) \left( A\otimes u\right)
r_{A}^{-1}$. The details are left to the reader.

c) Left also as an exercise.
\end{proof}

\begin{definitions}
1) Let $(A,m,u)$ be an algebra in $(\mathcal{M},\otimes
,\mathbf{1})$ and let $M $ be an $(A,A)$\lin{}{}bimodule. An
\emph{Hochschild extension }of $A$ with kernel $M$, is an algebra
homomorphism $\pi :E\rightarrow A$ that satisfies the following
conditions:

a) there is a section $\sigma $ of $\pi ;$

b) there is a morphism $i:M\rightarrow E$ such that $\left(
M,i\right) $ is the kernel of $\pi $ in $\mathcal{M};$

c) $m \left( i\otimes i\right) =0$ (i.e. $M^{2}=0$);

d) the $(A,A)$\lin{}{}bimodule structure of $M$ coincides with the
one induced by $ i$ (by the previous lemma ${M}$ is an
$(A,A)$\lin{}{}bimodule with the module structure
(\ref{eq:LeftMod}) and (\ref{eq:RightMod})).

2) Two Hochschild extensions $\pi :E\rightarrow A$ and $\pi
':E'\rightarrow A$ of $A$ with kernel $M$ are equivalent if
there is a morphism of algebras $f:E\rightarrow E'$ such that $%
\pi'f=\pi $ and $f':{\mathrm{Ker\,} }\pi \rightarrow {%
\mathrm{Ker\,} }\pi'$, the restriction of $f$, is an isomorphism in $%
\mathcal{M}$.
\end{definitions}

\begin{remark}
The morphism $f'$ is an isomorphism of $(A,A)$\lin{}{}bimodules. By $5$%
\lin{}{}Lemma $f$ is always an isomorphism of algebras.
\end{remark}

\begin{lemma}
\label{X lem 1.5.7}Let $(A,m,u)$ be an algebra and let $(M,\mu
_{r},\mu _{l}) $ an $(A,A)$\lin{}{}bimodule. Suppose that $\omega
:A\otimes A\rightarrow M$ is a Hochschild 2\lin{}{}cocycle. If
$m_{\omega }:\left( A\oplus M\right) \otimes
\left( A\oplus M\right) \rightarrow A\oplus M$ and $u_{\omega }:\mathbf{1}%
\rightarrow A\oplus M$ are defined by:
\begin{eqnarray*}
m_{\omega }:&=&i_{A} m \left( p_{A}\otimes p_{A}\right)+i_M \left[
\mu _{r} \left( p_{M}\otimes p_{A}\right) +\mu _{l} \left(
p_{A}\otimes p_{M}\right) -\omega \left( p_{A}\otimes p_{A}\right)
\right] , \\
u_{\omega }:&=&i_{A} u+i_{M} \omega \left( u\otimes u\right) l_{\mathbf{1}%
}^{-1},
\end{eqnarray*}
where $i_A,\ i_M$ are the canonical injections in $A\oplus M$ and
$p_A,\ p_M$ are the canonical projections. Then $\left( A\oplus
M,m_{\omega },u_{\omega }\right)$ is an algebra in $\mathcal{M}$.
Moreover $p_{A}:A\oplus M\rightarrow A$ is a Hochschild extension
of $A$ with kernel $\mathrm{Ker\,} p_{A}=(M,i_{M}).$ This
extension will be denoted by $E_{\omega }$.
\end{lemma}

\begin{proof}
Very tedious computations, that will be skipped. We just remark that the $%
m_\om$ defines an associative multiplication because $\omega$ is a
cocycle.
\end{proof}

\begin{definitions}
a) \ The Hochschild extension $p_{A}:E_{\omega }\rightarrow A,$
introduced
in the Lemma above, is called the \emph{Hochschild extension associated to }$%
\omega .$

b) \ If $(A,m_{A},u_{A})$ and $(E,m_{E},u_{E})$ are algebras and
$\sigma
:A\rightarrow E,$ is a unital morphism in $\mathcal{M}$ we define the \emph{%
curvature} of $\sigma$ to be the morphism:
\begin{equation}
\theta _{\sigma }:A\otimes A\rightarrow E,\quad\theta _{\sigma
}:=\sigma m_{A}-m_{E} (\sigma \otimes \sigma )  \label{X formula
(curvature)}
\end{equation}
\end{definitions}

\begin{proposition}
\label{X lem 1.5.9}Let $\pi :E\rightarrow A$ be a Hochschild
extension of $A $ with kernel $(M,i)$, let $\sigma :A\rightarrow
E$ be a section of $\pi $ and let $\theta _{\sigma }$ be the
morphism defined by formula (\ref{X formula (curvature)}). Then
there exists a unique morphism $\omega :A\otimes A\rightarrow M$,
such that $i \omega =\theta _{\sigma }.$ Moreover $\omega $
is a $2$\lin{}{}cocycle whose class $\left[ \omega \right] $ in $\mathbf{H}%
^{2}\left( A,M\right) $ does not depend on the choice of $\sigma .$ If $%
p_{A}:E_{\omega }\rightarrow A$ is the Hochschild extension associated to $%
\omega ,$ the morphism
\begin{equation*}
f_{\omega }:=\sigma p_{A}+i p_{M}:E_{\omega }\rightarrow E
\end{equation*}
defines an equivalence of Hochschild extensions.
\end{proposition}

\begin{proof}
The morphism $\pi $ is an algebra homomorphism, so that $\pi
\theta _{\sigma }=0$. Then there exists a unique morphism $\omega
:A\otimes A\rightarrow M$ such that $i \omega =\theta _{\sigma }.$
Let $\mu _{l}$ and $\mu _{r}$ be the morphisms that define the
module structure of $M$ and let $m$ and $m_{E}$ be the
multiplications of $A$ and $E$. We have:
\begin{equation*}
i b^{2}\left( \omega \right) =m_{E} \left( \sigma \otimes \theta
_{\sigma }\right) -\theta _{\sigma } (m \otimes A)+\theta _{\sigma
}(A\otimes m)-m_E(\theta _\sigma\otimes \sigma) ,
\end{equation*}
and, in view of the definition of $\theta _{\sigma },$ we obtain
$i b^{2}\left( \omega \right) =0$ so that $b^{2}\left( \omega
\right) =0.$

Let $\sigma ':A\rightarrow E$ be another section of $\pi $. Since $%
\pi \left( \sigma-\sigma '\right) =0,$ there exists a unique
morphism $\tau :A\rightarrow M$ such that $i \tau =\sigma-\sigma
'. $ If $\omega '$ is the associated to $\sigma ',$ a
straightforward computation shows us that:
\begin{equation*}
i \left( \omega '-b^{1}\left( \tau \right) \right) =i \omega ,
\end{equation*}
so $\omega '=b^{1}\left( \tau \right) +\omega .$ Thus $\left[
\omega \right] =\left[ \omega '\right] .$

As $m \left( i\otimes i\right) =0$ one can check that $f_{\omega
}$ is an algebra homomorphism. Moreover $\pi f_{\omega }=p_{A}$
and $f_{\omega }
i_{M}=i$ (the restriction of $f_{\omega }$ to $M$ is the isomorphism $%
\mathrm{Id}_{M}$). Thus $f_{\omega }$ is an equivalence of
Hochschild extensions.
\end{proof}

\begin{definition}
With the notations of the previous Proposition, the class $\left[
\omega \right] $ is called the \emph{cohomology class associated
to the Hochschild extension} $\pi :E\rightarrow A.$
\end{definition}

\begin{lemma}
\label{X rem associated class}Let $\omega :A\otimes A\rightarrow M$ be a $2$%
\lin{}{}cocycle and let $p_{A}:E_{\omega }\rightarrow A$ be the
Hochschild extension associated to $\omega$. Then the cohomology
class associated to the Hochschild extension $p_{A}:E_{\omega
}\rightarrow A$ is exactly $\left[ \omega \right] .$
\end{lemma}

\begin{proof}
Since $i_{A}:A\rightarrow E_{\omega }$ is a section of $p_{A}$, we
have:
\begin{equation*}
\theta _{i_{A}} =i_{A} m-m_{\omega } \left( i_{A}\otimes
i_{A}\right) =i_{A} m-i_{A} m+i_{M} \omega =i_{M} \omega .
\end{equation*}
Thus, in view of Proposition \ref{X lem 1.5.9}, the cohomology
class associated to this extension is $[\omega]$.
\end{proof}

\begin{claim}
Let $A$ be an algebra and let $M$ be an $(A,A)$\lin{}{}bimodule.
If $\pi
:E\rightarrow A$ is an Hochschild extension, we will denote by $\left[ E%
\right] $ the class of the Hochschild extensions equivalent to it.
We define:
\begin{equation*}
\mathrm{Ext}\left( A,M\right) :=\left\{ \left[ E\right] \mid E\rightarrow A%
\text{ is a Hochschild extension of }A\text{ with kernel
}M\right\} .
\end{equation*}
\end{claim}

\begin{proposition}
\label{X pro 1.5.12}Let $A$ be an algebra and let $M$ be an $(A,A)$%
\lin{}{}bimodule. If $\omega ,\omega ':A\otimes A\rightarrow M$ are $2$%
\lin{}{}cocycle, then:
\begin{equation*}
\left[ \omega \right] =\left[ \omega '\right] \Longleftrightarrow %
\left[ E_{\omega }\right] =\left[ E_{\omega '}\right] .
\end{equation*}
\end{proposition}

\begin{proof}
Suppose that $\left[ E_{\omega }\right] =\left[ E_{\omega
'}\right]. $ There exists an algebra homomorphism $g:E_{\omega
}\rightarrow E_{\omega '}$ that is an equivalence of Hochschild
extensions. As $g i_{A}$ is a section of $p_{A}':E_{\omega
'}\rightarrow A,$ we have:
\begin{equation*}
\theta _{g i_{A}}' =i_{M}' \omega,
\end{equation*}
so that, by Lemma \ref{X rem associated class}, $\left[ \omega \right] =%
\left[ \omega '\right] .$

If $\left[ \omega \right] =\left[ \omega '\right] ,$ there exists
a morphism $\tau :A\rightarrow M$ such that $\omega =\omega
'+b^{1}\left( \tau \right) .$ The morphism $\sigma :=i_{A}+i_{M}
\tau :A\rightarrow E_{\omega }$ is a section of $p_{A}:E_{\omega
}\rightarrow A.$ Thus:
\begin{equation*}
\theta _{\sigma } =i_{M} \omega '.
\end{equation*}
Applying Proposition \ref{X lem 1.5.9} in the case when
$E:=E_{\omega },$ we get that there is an equivalence between
$p_{A}':E_{\omega '}\rightarrow A$ and $p_{A}:E_{\omega
}\rightarrow A,$ so that $\left[ E_{\omega }\right] =\left[
E_{\omega '}\right] .$
\end{proof}

\begin{theorem}
\label{X teo 1.5.13}Let $A$ be an algebra and let $M$ be an $(A,A)$%
\lin{}{}bimodule. The map:
\begin{equation*}
\Phi :\mathbf{H}^{2}( A,M) \rightarrow \mathrm{Ext}\left(
A,M\right) ,
\end{equation*}
where $\Phi \left( \left[ \omega \right] \right) :=\left[
E_{\omega }\right] ,$ is well defined and is a bijection.
\end{theorem}

\begin{proof}
$\Phi $ is well defined and injective by Proposition \ref{X pro
1.5.12}, and it is surjective by Proposition \ref{X lem 1.5.9}.
\end{proof}

\begin{definition}
An extension $\pi :E\rightarrow A$ is a \emph{trivial extension
}whenever it admits a section that is an algebra homomorphism.
\end{definition}

\begin{corollary}
\label{X coro 1.5.16}Let $A$ be an algebra and let $M$ be an $(A,A)$%
\lin{}{}bimodule. Then a Hochschild extension of $A$ with kernel
$M$ is trivial if and only if the associated cohomology class is
zero.
\end{corollary}

\begin{proof}
Let $\pi :E\rightarrow A$ be a Hochschild extension of $A$ with
kernel $M$, and let $i:M\rightarrow E$ the canonical injection. By
the definition of trivial extensions, there exists a section
$\sigma :A\rightarrow E$ of $\pi $ that is an algebra
homomorphism. Thus $i \omega =\theta _{\sigma }=\sigma m-m_{E}
\left( \sigma \otimes \sigma \right) =0$, so that $[\omega] =0.$

If $\left[ \omega \right] =0,$ where $\omega $ is the
$2$\lin{}{}cocycle
associated to $\pi :E\rightarrow A,$ then $\left[ E\right] =\left[ E_{0}%
\right] $ that is there exists an algebra homomorphism
$f:E_{0}\rightarrow E$ that is an equivalence of Hochschild
extensions. Let $\sigma _{0}:A\rightarrow E_{0}$ be a unital
section of $p_{0}:E_{0}\rightarrow A$. Then $0=i 0=\theta _{\sigma
_{0}}={\sigma _{0}} m-m_{E} \left( {\sigma _{0}} \otimes {\sigma
_{0}} \right) $ so that $\sigma _{0}$ is an algebra homomorphism.
It is easy to see that $f \sigma _{0}$ is a section of $\pi $ that
is an algebra homomorphism.
\end{proof}

\begin{corollary}
\label{X coro 1.5.17}Let $A$ be an algebra and let $M$ be an $(A,A)$%
\lin{}{}bimodule. Then $\mathbf{H}^{2}( A,M)=0$ if and only if any
Hochschild extension of $A$ with kernel $M$ is trivial.
\end{corollary}

\begin{proof}
It follows by Theorem \ref{X teo 1.5.13} and Corollary \ref{X coro
1.5.16}.
\end{proof}

\begin{definition}
Let $\left( A,m,u\right) $ be an algebra and let $f:X\rightarrow
A$ be a morphism in $\mathcal{M}$. If $\Lambda _{f}:=m \left(
m\otimes A\right) \left( A\otimes f\otimes A\right)$ then the
\emph{two\lin{}{}sided ideal of} $A$ \emph{generated by} $f$ is
defined by:
\begin{equation*}
({}_{A}\!\left\langle f\right\rangle _{A},i_{f})
:=\mathrm{Im}\left( \Lambda _{f}\right).
\end{equation*}
\end{definition}

\begin{lemma}
\label{UniE lem 12.2}Let $\left( A,m,u\right) $ be an algebra and let $%
f:X\rightarrow A$ be a morphism in $\mathcal{M}$. Then:

a) $\left( _{A}\left\langle f\right\rangle _{A},i_{f}\right) \ $is
an ideal of $A{.}$

b) If $\xi :A\rightarrow Y$ is a morphism such that $\xi m=m (\xi
\otimes \xi )$ and $\xi f=0,$ then $\xi i_{f}=0.$

c) If $\xi :A\rightarrow Y$ is a morphism such that $\xi i_{f}=0,$
then $\xi f=0.$
\end{lemma}

\begin{claim}
Let $(A,m_{A},u_{A})$ be an algebra in $(\mathcal{M},\otimes
,\mathbf{1}).$
Let us consider the tensor algebra $T\left( A\right) :=\oplus _{n\in \mathbb{%
N}}T_{n}\left( A\right) ,$ where $T_{0}\left( A\right) :=\mathbf{1}$ and $%
T_{n+1}\left( A\right) :=T_{n}\left( A\right) \otimes A,\forall
n>0.$ We set:
\begin{equation*}
\left( I,\varsigma _{I}\right) :={{_{T\left( A\right)
}}}\left\langle u_{T\left( A\right) }-i_{A} u_{A}\right\rangle
_{T\left( A\right) },
\end{equation*}
where $u_{T\left( A\right) }:\mathbf{1}\rightarrow T\left( A\right) $ and $%
i_{A}:A\rightarrow T\left( A\right) $ are the canonical morphisms.
Moreover we set:
\begin{equation*}
\left( E_{A},\rho _{A}'\right) :=\mathrm{Coker}\left( \varsigma
_{I}\right) .
\end{equation*}
Since $I$ an ideal of $T\left( A\right) $, $E_{A}$ is an algebra
and $\rho _{A}'$ an algebra homomorphism: by the previous Lemma,
$\rho _{A}' \left( u_{T\left( A\right) }-i_{A} u_{A}\right) =0.$
Let $\rho _{A}=\rho _{A}' i_{A}:A\rightarrow E_{A}$. Then we have:
\begin{equation*}
\rho _{A} u_{A}=\rho _{A}' i_{A} u_{A}=\rho _{A}'
u_{T(A)}=u_{E_{A}}.
\end{equation*}
So, by construction, $\rho _{A}$ is a unital morphism.
\end{claim}

\begin{definition}
$\left( E_{A},\rho _{A}\right) $ is called \emph{the universal
extension}. The following proposition justifies this name.
\end{definition}

\begin{proposition}
\label{universal} Let $A,B$ be algebras in $(\mathcal{M},\otimes ,\mathbf{1}%
).$ Given a unital morphism $\rho :A\rightarrow B,$ there exists a
unique algebra homomorphism $v:E_{A}\rightarrow B$ such that $v
\rho _{A}=\rho .$
\end{proposition}

\begin{proof}
By the universal property of the tensor algebra $T\left( A\right)
,$ there exists a unique algebra homomorphism $\xi :T\left(
A\right) \rightarrow B$ such that $\xi i_{A}=\rho .$ Then:
\begin{equation*}
\xi \left( u_{T\left( A\right) }-i_{A} u_{A}\right) =\xi
u_{T\left( A\right) }-\rho u_{A}=u_{B}-u_{B}=0,
\end{equation*}
where $u_{B}$ is the unit of $B.$ By Lemma \ref{UniE lem 12.2}(b),
$\xi \varsigma _{I}=0$ so that there exists a unique morphism
$v:E_{A}\rightarrow B$ such that $v \rho _{A}'=\xi $ and hence $v
\rho _{A}=v \rho _{A}' i_{A}=\xi i_{A}=\rho $. Moreover such a
morphism is an algebra map. The uniqueness is due to the universal
property of $T\left( A\right) $.
\end{proof}

\begin{corollary}
There exists a unique algebra map $\pi _{A}:E_{A}\rightarrow A$ such that $%
\pi _{A} \rho _{A}=\mathrm{Id}_{A}.$
\end{corollary}

\begin{claim}
Let $\left( e_{A},i_A\right) =\mathrm{Ker\,} \pi _{A}.$ We have
the exact sequence:
\begin{equation*}
0\rightarrow e_{A}\overset{i_A}{\rightarrow }E_{A}\overset{\pi _{A}}{%
\rightarrow }A\rightarrow 0.
\end{equation*}
From this sequence, we obtain an Hochschild extension of $A$, namely$:$%
\begin{equation}
0\rightarrow \frac{e_{A}}{e_{A}^{2}}\rightarrow \frac{E_{A}}{e_{A}^{2}}%
\rightarrow A\rightarrow 0,  \label{X uni Hoch ext}
\end{equation}
where the section of $E_{A}/e_{A}^{2}\rightarrow A$ is given by
the composition of $E_{A}\rightarrow E_{A}/e_{A}^{2}$ and $\rho
_{A}:A\rightarrow E_{A}.$ The extension (\ref{X uni Hoch ext}) is
called the \emph{universal Hochschild extension of} $A$.
\end{claim}

\begin{proposition}
Let $A,B$ be algebras in $(\mathcal{M},\otimes ,\mathbf{1}),$ let $%
0\rightarrow M\overset{i}{\rightarrow }E\overset{\pi }{\rightarrow }%
B\rightarrow 0$ be an Hochschild extension of $B$ with kernel $M$ and let $%
f:A\rightarrow B$ be an algebra homomorphism. Then, there exists
an algebra
homomorphism $\pi _{f}:E_{A}/e_{A}^{2}$ $\rightarrow E$ and an $(A,A)$%
\lin{}{}bimodule homomorphism $g:e_{A}/e_{A}^{2}\rightarrow M$
such the following diagram commutes:
\begin{equation*}
\begin{diagram}[h=2em,w=2.5em] 0&\rTo&e_A/e_A^2& \rTo& E_A/e_A^2
&\rTo&A&\rTo & 0 \\ & &\dDotsto<{g}& &\dDotsto^{\pi_f} & &
\dTo>{f} & & \\ 0&\rTo& M& \rTo^{i}& E &\rTo^{\pi}&B&\rTo & 0 \\
\end{diagram}
\end{equation*}
\end{proposition}

\begin{proof}
Let $\rho :B\rightarrow E$ be a unital section of $\pi .$ By
Proposition \ref {universal} there exists a unique algebra
homomorphism $\pi _{f}':E_{A}\rightarrow E$ such that $\rho f=\pi
_{f}' \rho _{A}$. Therefore we get:
\begin{equation*}
f \pi _{A} \rho _{A}=f=\pi \rho f=\pi \pi _{f}' \rho _{A}.
\end{equation*}
Thus, by Proposition \ref{universal}, we get $f \pi _{A}=\pi \pi
_{f}'.$ Since $M=\mathrm{Ker\,} \pi $ the relation $\pi \pi _{f}'
i_A=f \pi _{A} i_A=0 $ implies the existence of a unique morphism
$\gamma :e_{A}\rightarrow M$ such that $i \gamma =\pi _{f}' i_A$.
Then, from $M^{2}=0$ we deduce:
\begin{equation*}
\pi _{f}' m_{E_{A}} \left( i_A\otimes i_A\right) =m_{E} \left( \pi
_{f}'\otimes \pi _{f}'\right) \left( i_A\otimes i_A\right) =m_{E}
\left( i\otimes i\right) \left( \gamma \otimes \gamma \right) =0,
\end{equation*}
so that there exists a unique morphism $\pi _{f}:E_{A}/e_{A}^{2}$ $%
\rightarrow E$ which, composed with the canonical projection $%
E_{A}\rightarrow E_{A}/e_{A}^{2},$ gives $\pi _{f}'.$ Since $%
e_A/e_A^2$ is the kernel of $\pi_A$ there is a unique $g$ such
that the left square of the above diagram in commutative.
\end{proof}

\begin{theorem}
\label{X coro UniE}Let $A$ be an algebra in $(\mathcal{M},\otimes ,\mathbf{1}%
)$. Then the following conditions are equivalent:

(a) The universal Hochschild extension of $A$ is trivial.

(b) If $\pi :E\rightarrow B$ is an algebra homomorphism that splits in $%
\mathcal{M}$ and $\left( \mathrm{Ker\,} \pi \right) ^{2}=0,$ then
any algebra homomorphism $f:A\rightarrow B$ can be lifted to an
algebra homomorphism $A\rightarrow E.$

(c) Given an algebra homomorphism $f:A\rightarrow B$ which admits
a section that is an algebra homomorphism, then any Hochschild
extension of $B$ is trivial.

(d) Any Hochschild extension of $A$ is trivial.

(e) $\mathbf{H}^{2}\left( A,M\right) =0,\forall M\in
{{_{A}\mathcal{M}_{A}.}}$

(f) Given an algebra homomorphism $f:A\rightarrow B$ which admits
a section
that is an algebra homomorphism, then the universal Hochschild extension of $%
B$ is trivial.
\end{theorem}

\begin{proof}
(a) $\Rightarrow$ (b) \ If $\pi :E\rightarrow B$ is an algebra
homomorphism that splits in $\mathcal{M}$ and $\left(
\mathrm{Ker\,} \pi \right) ^{2}=0,$ then $\pi $ is an Hochschild
extension of $B$. If $\sigma :A\rightarrow
E_{A}/e_{A}^{2}$ is an algebra morphism that is a section of the morphism $%
E_{A}/e_{A}^{2}\rightarrow A,$ then $\pi _{f} \sigma :A\rightarrow
E$ is the algebra morphism that lifts $f$.

(b)$\Rightarrow$ (c) \ If $\pi :E\rightarrow B$ is an Hochschild
extension,
then $f:A\rightarrow B$ can be lifted to an algebra homomorphism $%
g:A\rightarrow E.$ If $\sigma :B\rightarrow E$ is a section of $f$
that is an algebra homomorphism, then $g \sigma $ is a section of
$\pi $ that is an algebra homomorphism, so that $\pi $ is trivial.

(c) $\Rightarrow$ (d) \ The identity of $A$ is an algebra
homomorphism and is its own section, so that any Hochschild
extension of $A$ is trivial.

(d) $\Rightarrow$ (a) \ It is obvious.

(d)$\Leftrightarrow$ (e) \ It follows by Corollary \ref{X coro
1.5.17}, and obviously (f) is equivalent to the others.
\end{proof}

\begin{lemma}
Let $A$ be an algebra and let $\theta _{A}:A\otimes A\rightarrow
E_{A}$ be the curvature associated to the canonical morphism $\rho
_{A}:A\rightarrow E_{A}.$ Then $e_{A}\simeq {_{E_{A}}\left\langle
\theta _{A}\right\rangle _{E_{A}}.}$
\end{lemma}

\begin{proof}
Let us denote by $(X,\phi)$ the cokernel of $\Lambda_{\theta_ A}$.
By definition, $_{E_{A}}\left\langle \theta _{A}\right\rangle
{}_{E_{A}}= \mathrm{Im}\,\Lambda _{\theta _{A}}$. As $\pi
_{A}:E_{A}\rightarrow A$ is an algebra homomorphism, then $\pi
_{A} \theta _{A}=0$. So, by Lemma \ref{UniE lem 12.2}, $\pi _{A}
i_{\theta _{A}}=0$.

Let $\beta : X\rightarrow A$ be such that $\pi _{A}=\beta \phi$. Since $%
X\simeq \mathrm{Coker}\,i_{\theta _{A}}$ and
$({_{E_{A}}\left\langle \theta _{A}\right\rangle
_{E_{A}},}i_{\theta _{A}})$ is an ideal of $E_{A}$ it follows that
$X$ has an algebra structure so that $\phi $ is an algebra
homomorphism. As, by definition, $\phi i_{\theta _{A}}=0$ we have $%
\phi \theta _{A}=0.$ This relation, the fact that $\phi $ is an
morphism and $\rho _{A}$ a unital morphism imply that $\phi \rho
_{A}:A\rightarrow X $ is an algebra homomorphism. We have:
\begin{equation*}
\phi \rho _{A} \pi _{A} \rho _{A}=\phi \rho _{A}.
\end{equation*}
By Proposition \ref{universal} we deduce that
\begin{equation*}
\phi \rho _{A} \pi _{A}=\phi.
\end{equation*}
In particular, $\phi \rho _{A}$ is an epimorphism. As $\beta \phi
\rho _{A}=\mathrm{Id}_{A}$ then $\phi \rho _{A}$ is a monomorphism
too.
Therefore we get $\left( A,\pi _{A}\right) \simeq X$, so that $%
_{E_{A}}\left\langle \theta _{A}\right\rangle _{E_{A}}\simeq
\mathrm{Ker\,} \pi _{A}=e_{A}.$
\end{proof}

\begin{theorem}
\label{X formally smooth teo}Let $A$ be an algebra in
$\mathcal{M}$. Then the following assertions are equivalent:

(a) $\mathrm{Ker\,} m$ is $\mathcal{E}$\lin{}{}projective, where
$m$ is the multiplication of $A$.

(b) $\mathbf{H}^{2}(A,M) =0,\forall M\in {{_{A}\mathcal{M}_{A}.}}$

(c) Let $\pi :E\rightarrow B$ be an algebra epimorphism and let
$I$ denote
the kernel of $\pi$. Assume that there is $n\in \mathbb{N}$ so that $I^{n}=0$%
. If for any $r=1,\cdots ,n-1$ the canonical projection $%
p_{r}:E/I^{r+1}\rightarrow E/I^{r}$ splits in $\mathcal{M}$ then
any algebra
homomorphism $f:A\rightarrow B$ can be lifted to an algebra homomorphism $%
g:A\rightarrow E.$

(d) Let $\pi :E\rightarrow A$ be an algebra epimorphism and let
$I$ denote
the kernel of $\pi$. Assume that there is $n\in \mathbb{N}$ so that $I^{n}=0$%
. If for any $r=1,\cdots ,n-1$ the canonical projection $%
p_{r}:E/I^{r+1}\rightarrow E/I^{r}$ splits in $\mathcal{M}$ then
$\pi $ has a section which is an algebra homomorphism.
\end{theorem}

\begin{proof}
(a) $\Leftrightarrow$ (b) Let $\left( L,j\right) :=\mathrm{Ker\,}
m$ and let us consider the exact sequence:
\begin{equation*}
0\rightarrow L\overset{j}{\rightarrow }A\otimes A\overset{m}{\rightarrow }%
A\rightarrow 0.
\end{equation*}
We know that $m$ has a section in $\mathcal{M}$ so that $m\in
\mathcal{E}$
and the sequence above is $\mathcal{E}$\lin{}{}exact. Given any $M\in {{_{A}%
\mathcal{M}_{A},}}$ we apply the functor
$F:={}_A\mathcal{M}_A(-,M)$ to the sequence above and find:
\begin{equation*}
\mathbf{Ext}^1_{\mathcal{E}}\left( A\otimes A,M\right) \rightarrow
\mathbf{Ext}^1_{\mathcal{E}}\left( L,M\right) \rightarrow
\mathbf{Ext}^2_{\mathcal{E}}\left( A,M\right) \rightarrow
\mathbf{Ext}^2_{\mathcal{E}}\left( A\otimes A,M\right) .
\end{equation*}
Since $A\otimes A$ is $\mathcal{E}$\lin{}{}projective, we get that
$\mathbf{Ext}^1_{\mathcal{E}}\left( L,M\right) \simeq
\mathbf{Ext}^2_{\mathcal{E}}\left( A,M\right) = \mathbf{H} ^{2}(
A,M)$. Then (a) and (b) are equivalent in view of Proposition
\ref{te:verysep}.

(b) $\Leftrightarrow$ (c) If we assume that (c) holds then, in
particular, we have the lifting property from Theorem \ref{X coro
UniE}(b). Hence, by the same theorem, $\mathbf{H}^{2}(A,M)=0$, for
every $(A,A)$\lin{}{}bimodule $M$.

Now let us assume that the second Hochschild cohomology group of
$A$ with coefficients in $M$ is trivial, for any
$M\in{}_A\mathcal{M}_A$. By Theorem \ref{X coro UniE} we know that
we have the required lifting property for all epimorphisms $\pi$
splitting in $\mathcal{M}$ and satisfying $I^2=0$. Let now $\pi$
be an arbitrary epimorphisms as in (c). Since $p_{r}:E/I^{r+1}
\rightarrow E/I^{r}$ splits in $\mathcal{M}$ and the square of its
kernel is trivial, inductively we can construct a sequence of
algebra morphisms $ f_1:=f, f_2,\cdots,f_n$ such that
$f_r:A\rightarrow E/I^r$ and $p_rf_{r+1} =f_r$. We conclude this
implication by remarking that $E=E/I^n$, so the lifting of $f$ can
be chosen to be $f_n$.

(b) $\Leftrightarrow$ (d) Similarly to the proof of (b)
$\Leftrightarrow$ (c), by using the fact that in Theorem \ref{X
coro UniE} the second and the fourth assertions are equivalent.
\end{proof}

\begin{definition}
Any algebra $\left( A,m,u\right) $ in $(\mathcal{M},\otimes
,\mathbf{1}),$ satisfying one of the conditions of Theorem \ref{X
coro UniE} or of Theorem \ref{X formally smooth teo}, is called
\emph{formally smooth.}
\end{definition}

\begin{corollary}
\label{sep-smooth}Any separable algebra is formally smooth.
\end{corollary}

\begin{proposition}
If $A$ is a formally smooth algebra and $M$ is an
$\mathcal{E}$\lin{}{}projective bimodule in $_A\mathcal{M}_A$,
then the tensor algebra $T_{A}\left( M\right) $ is also formally
smooth.
\end{proposition}

\begin{proof}
Let $\pi :E\rightarrow B$ be an algebra homomorphism that splits in $%
\mathcal{M} $ and such that $\left( \mathrm{Ker\,} \pi \right) ^{2}=0.$ Let $%
f:T_{A}\left( M\right) \rightarrow B$ be an algebra homomorphism.
Since $A$ formally smooth, by the second condition from Theorem
\ref{X coro UniE}, there exists an algebra homomorphism
$g_{0}:A\rightarrow E$ such that $\pi g_{0}=f i_{A},$ where
$i_{A}:A\rightarrow T_{A}\left( M\right) $ is the
canonical inclusion. The objects $E$ and $B$ have a natural $(A,A)$%
\lin{}{}bimodule induced by $g_{0}$ and $f i_{A}$, respectively.
Thus $\pi $ and $f $ become homomorphisms of
$(A,A)$\lin{}{}bimodules. Let $i_{M}:M\rightarrow
T_{A}\left( M\right) $ be the canonical inclusion. Since $M$ is $\mathcal{E}$%
\lin{}{}projective there exists a morphism of $(A,A)$\lin{}{}bimodules $%
g_{1}:M\rightarrow E$ such that $\pi g_{1}=f i_{M}.$ By the
universal property of $T_{A}\left( M\right) ,$ there exists a
unique algebra homomorphism $g:T_{A}\left( M\right) \rightarrow E$
such that $g i_{A}=g_{0}$ and $g i_{M}=g_{1}:$ then $\pi g
i_{A}=\pi g_{0}=f i_{A}$ and $\pi g i_{M}=\pi g_{1}=f i_{M}.$
Finally $\pi g=f.$
\end{proof}

\begin{corollary}
If $(A,m,u)$ is a formally smooth algebra, the tensor algebra
$T_{A}\left( \mathrm{Ker\,} m\right) $ is also formally smooth.
\end{corollary}

\section{Separable and formally smooth algebras in ${}^H \mathfrak{M}^H$}

In this section we shall apply the results of the previous section
to study separability and formally smoothness of algebras in the
monoidal category of all $(H,H)$\lin{}{}bicomodules, where $H$ is
a given Hopf algebra.

\begin{claim}
Let $H$ be a Hopf algebra. Let us consider the category $\mathcal{M} :={}^{H}%
\mathfrak{M}^{H} $ of $(H,H)$\lin{}{}bicomodules with the tensor
product $\left( -\right) \otimes _{K}\left(-\right) $ as in
Example \ref{ex:categorii}(c).
Hence an algebra in $\mathcal{M}$ is an algebra $A$ which is an $(H,H)$%
\lin{}{}bicomodule such that $A$ is a left and a right
$H$\lin{}{}comodule algebra. We shall say that $A$ is an
$H$\lin{}{}\emph{bicomodule algebra}.

Let $A$ be an $H$\lin{}{}bicomodule algebra. The category of all $(A,A)$%
\lin{}{}bimodules in $\mathcal{M}$ will be denoted by
$_{A}^{H}\mathfrak{M}_{A}^{H}
$. An $(H,H)$%
\lin{}{}bicomodule $M$ is an object in
$_{A}^{H}\mathfrak{M}_{A}^{H}$ if it is an $
(A,A)$\lin{}{}bimodule too such that $\mu _{l}:A\otimes
M\rightarrow M$ and $\mu _{r}:M\otimes A\rightarrow M$ are
morphisms of $(H,H)$\lin{}{}bicomodules. Here $ \mu_{l}$ and
$\mu_{r}$ define the module structures on $M$ and $A\otimes M$ is
an $(H,H)$\lin{}{}bicomodule with the diagonal coactions. For
$A=K$ with trivial $H$\lin{}{}comodule structures we get the
category of $(H,H)$\lin{}{}bicomodules.
Also for the trivial Hopf algebra $H=K$ we get that $A$ is just a $K$%
\lin{}{}algebra, and
$_{A}^{H}\mathfrak{M}_{A}^{H}={}_{A}\mathfrak{M}_{A}$ .

${_{A}^{H}}\mathfrak{M}{_{A}^{H}}$ is a monoidal category with the
usual tensor product of two $(A,A)$\lin{}{}bimodules $(-)\otimes
_{A}(-)$. If $V,W\in _{A}^{H} \mathfrak{M}_{A}^{H}$ then the left
structures on $V\otimes _{A}W$ are given by:
\begin{eqnarray*}
r\left( v\otimes _{A}w\right) &=&rv\otimes _{A}w \\ \rho
_{V\otimes _{H}W}^{l}\left( v\otimes _{A}w\right) &=&\sum
v{}_{\langle -1\rangle}w{}_{\langle-1\rangle}\otimes (v{}_{\langle
0\rangle}\otimes _{A}w{}_{\langle 0\rangle}).
\end{eqnarray*}
The right structures are defined similarly. The unit in $_{A}^{H}\mathfrak{M}%
_{A}^{H}$ is $A$. 

By definition, an algebra $A$ in ${}^H \mathfrak{M}^H$ is
separable if and only if the multiplication $m:A\otimes
A\rightarrow A$ has a section $\sigma :A\rightarrow A\otimes A$
which is a morphism of $(A,A)$\lin{}{}bimodules and $
(H,H)$\lin{}{}comodules. Obviously, then $A$ is separable as an algebra in $%
\mathfrak{M}_K$, but the converse does not hold in general.
Nevertheless, if
the forgetful functor $F:{_{A}^{H}\mathfrak{M} _{A}^{H}\rightarrow }${{$_{A}%
\mathfrak{M}_{A}$ is separable, then }}$A$ is separable as an algebra in $%
\mathcal{M}$. Before to prove this result, let us recall the
definition and basic properties of separable functors.
\end{claim}

\begin{claim}
A functor $F:\mathfrak{C}\rightarrow \mathfrak{D}$ is called
separable if,
for all objects $C_{1},C_{2}\in \mathfrak{C}$, there is a map $\varphi :Hom_{%
\mathfrak{D}}\left( FC_{1},FC_{2}\right) \rightarrow Hom_{\mathfrak{C}%
}\left( C_{1},C_{2}\right) $ such that:

1) For all $f\in Hom_{\mathfrak{C}}\left( C_{1},C_{2}\right) ,$
$\varphi \left( F\left( f\right) \right) =f$

2) We have $\varphi(k) f=g\varphi(h)$ for every commutative diagram in $%
\mathfrak{D}$ of type:
\begin{equation*}
\begin{diagram}[h=2.em,w=1.5em]
 F(C_1) & \rTo^{h } & F(C_2) \\ \dTo<{F(f) } & & \dTo>{ F(g)
} \\ F(C_3) & \rTo_{k} & F(C_4) \\ \end{diagram}
\end{equation*}
If $F:\mathfrak{C}\rightarrow \mathfrak{D}$ has a right adjoint $G:%
\mathfrak{D}\rightarrow \mathfrak{C}$, then $F$ is separable (see \cite{Raf}%
) iff the unit $\sigma :1_{\mathfrak{C}}\rightarrow GF$ splits,
i.e. there
is a functorial morphism $\mu :GF\rightarrow 1_{\mathfrak{C}}$ such that $%
\mu \sigma =\mathrm{Id}_{1_{\mathfrak{C}}}.$ If $F$ is separable and $%
F\left( f\right) $ has a section in $\mathfrak{D}$, then $f$ has a
section in $\mathfrak{C}.$
\end{claim}

\begin{claim}
The forgetful functor $F:{_{A}^{H}\mathfrak{M}_{A}^{H}\rightarrow }${{$_{A}\mathfrak{M}%
_{A}$ has a right adjoint }}$G:{{_{A}\mathfrak{M}_{A}}\rightarrow {_{A}^{H}%
\mathfrak{M}_{A}^{H},G}}\left( M\right) =H\otimes M\otimes H$, where $%
G\left( M\right) $ is a bicomodule via $\Delta _{H}\otimes M\otimes H$ and $%
H\otimes M\otimes \Delta _{H}$, and $G\left( M\right) $ is a
bimodule with diagonal actions:
\begin{eqnarray*}
a\left( h\otimes m\otimes k\right) &=&\sum a_{<-1>}h\otimes
a_{<0>}m\otimes a_{<1>}k \\ \left( h\otimes m\otimes k\right) a
&=&\sum ha_{<-1>}\otimes ma_{<0>}\otimes ka_{<1>}.
\end{eqnarray*}
Here we used the $\Sigma$\lin{}{}notation: $(\rho _{A}^{l}\otimes
A)\rho _{A}^{r}\left( a\right) =\sum a_{<-1>}\otimes
a_{<0>}\otimes a_{<1>}.$ For any $M\in $
${{_{A}^{H}\mathfrak{M}_{A}^{H}}}$ the unit of the adjunction is
the map $\sigma _{M}:M\rightarrow H\otimes M\otimes H,\sigma
_{M}=(\rho _{M}^{l}\otimes H) \rho _{M}^{r}.$
\end{claim}

\begin{lemma}
\label{le:semicose} Let $H$ be a semisimple and cosemisimple Hopf
algebra over a field $K$. Then there is a left and a right
integral $\lambda $ such that $\lambda (1)=1$ and
\begin{equation}
\lambda \left( \sum h_{\left( 1\right) }xSh_{\left( 2\right)
}\right) =\lambda \left( \sum Sh_{\left( 1\right) }xh_{\left(
2\right) }\right) =\varepsilon \left( h\right) \lambda \left(
x\right) ,\forall h,x\in H. \label{ec:AdInv}
\end{equation}
\end{lemma}

\begin{proof}
First let us note that any semisimple Hopf algebra is finite
dimensional.
Hence, in view of \cite[Corollary 3.2]{EtGe}, we get that $S^{2}=\mathrm{Id}%
_{H}$. On the other hand, by \cite[Exercises 5.5.9 and
5.5.10]{DNC}, $H$ is unimodular and there is a (unique) right and
left integral $\lambda \in H^{\ast }$ such that $\lambda \left(
1\right) =1.$ Hence equation 1(a) in \cite[Theorem 3]{Radf}
becomes in this particular case $\lambda (hk)=\lambda (kh),$
$\forall h,k\in H,$ as $S^{2}=\mathrm{Id}_{H}$ and $H$ is
unimodular. Therefore :
\begin{equation*}
\lambda (\sum Sh_{\left( 1\right) }xh_{\left( 2\right) })=\sum
\lambda (xh_{\left( 2\right) }Sh_{\left( 1\right) })=\varepsilon
(h)\lambda (x),
\end{equation*}
where for the last equality we used $\sum h_{\left( 2\right)
}Sh_{\left(
1\right) }=\varepsilon (h),$ relation that holds since $S^{2}=\mathrm{Id}%
_{H}.$ The second equation of (\ref{ec:AdInv}) can be proved
similarly.
\end{proof}

\begin{remark}
Suppose that $H$ is a Hopf algebra. Then $H$ is cosemisimple and
has a non\lin zero left and right integral $\lambda$ verifying
(\ref{ec:AdInv}) if and only if there is a (necessarily unique)
left and right integral $\lambda$ such that (\ref{ec:AdInv}) holds
true and $\lambda(1)=1$. Indeed, two non\lin zero integrals are
proportional, hence any non\lin zero integral verifies
(\ref{ec:AdInv}). On the other hand $H$ is cosemisimple if and
only if there is a (unique) integral $\lambda$ such that
$\lambda(1)=1$.
\end{remark}

\begin{definition}\label{ad-invariant}
A left and right integral $\lambda$ verifying (\ref{ec:AdInv}) and
$\lambda(1)=1$ will be called an $ad$\lin invariant integral.
\end{definition}

\begin{theorem}
\label{pr:Fsep}Let $H$ be a Hopf algebra with an $ad$%
\lin{}{}invariant integral $\lambda$. Then $F:{_{A}^{H}\mathfrak{M}%
_{A}^{H}\rightarrow } ${{$_{A}\mathfrak{M}_{A}$ is a separable
functor}}.
\end{theorem}

\begin{proof}
We have to construct a functorial section of $\left( \sigma
_{M}\right) _{M\in {_{A}^{H}\mathfrak{M}_{A}^{H}}}.$ Let $\lambda$
be an $ad$\lin invariant integral. Let
\[\mu _{M}:H\otimes
M\otimes H\rightarrow M,\quad\mu _{M}\left( h\otimes m\otimes
k\right) =\sum \lambda \left( Sh\,m_{<-1>}\right) m_{<0>}\lambda
\left( m_{<1>}Sk\right) .
\]
Obviously $\left( \mu _{M}\right) _{M\in {_{A}^{H}%
\mathfrak{M}_{A}^{H}}}$ is a functorial morphism. Let us check
that $\mu _{M} $ is a morphism in
${_{A}^{H}\mathfrak{M}_{A}^{H}}$, i.e. $\mu _{M}$ is a morphism of
$(A,A)$\lin{}{}bimodules and a morphism of
$(H,H)$\lin{}{}bicomodules. Let $x=\mu _{M}(a\left( h\otimes
m\otimes k\right) )$. Then we have:
\begin{equation*}
x =\sum \lambda \left( ShSa_{<-2>}\,a_{<-1>}m_{<-1>}\right)
a_{<0>}m_{<0>}\lambda \left( a_{<1>}m_{<1>}SkSa_{<2>} \right).
\end{equation*}
Hence $\mu _{M}(a\left( h\otimes m\otimes k\right) )=a \mu
_{M}(h\otimes m\otimes k)$. This relation proves that $\mu _{M}$
is left $A$-linear.
Similarly, using the second equality of (\ref{ec:AdInv}), one can show that $%
\mu _{M}$ is right $A$\lin{}{}linear. We have:
\begin{equation*}
\sum h_{\left( 1\right) }\otimes \mu _{M}\left( h_{\left( 2\right)
}\otimes m\otimes k\right) =\sum h_{\left( 1\right) }\lambda
\left( Sh_{\left( 2\right) }\,m_{<-1>}\right) \otimes
m_{<0>}\lambda \left( m_{<1>}Sk\right) .
\end{equation*}
Let $y:=\sum h_{\left( 1\right) }\lambda \left( Sh_{\left(
2\right) }m_{<-1>}\right)$. Then, since $\lambda$ is a right
integral, we have:
\begin{equation*}
Sy =\sum \lambda \left( Sh\,m{\un{-1}}_{(1)}\right)
Sm{{\un{-1}}_{(2)}}.
\end{equation*}
Thus $y =\sum \lambda \left( Sh\,m{\un{-1}}_{(1)}\right)
m{{\un{-1}}_{(2)}}$, so:
\begin{eqnarray*}
\sum h_{\left( 1\right) }\otimes \mu _{M}\left( h_{\left( 2\right)
}\otimes m\otimes k\right) &=&\sum \lambda \left( Sh\,
m_{<-2>}\right) m_{<-1>}\otimes m_{<0>}\lambda \left(
m_{<1>}Sk\right).
\end{eqnarray*}
As $\rho _{M}^{l}\left( \mu _{M}\left( h\otimes m\otimes k\right)
\right)$ equals the right hand side of the above equation, we have
shown that $\mu _{M}$ is left\lin{}{}colinear. Analogously it can
be proved that $\mu _{M}$ is right $H$\lin{}{}colinear. It remains
to show that $\mu _{M}$ is a retraction of $\sigma _{M}.$ But:
\begin{eqnarray*}
\left( \mu _{M}\sigma _{M}\right) \left( m\right) =\sum \lambda
\left( Sm_{<-2>}\,m_{<-1>}\right) m_{<0>}\lambda \left(
m_{<1>}Sm_{<2>}\right)=m,
\end{eqnarray*}
so the theorem is proved.
\end{proof}

\begin{theorem}
\label{HH separable} Let $H$ be a Hopf algebra over a field $K$
and assume that $H$ has a $ad$\lin{}{}invariant integral. An
algebra $A$ in the category
${}^H \mathfrak{M}^H$ is separable iff $A$ is separable as an algebra in $(%
\mathfrak{M}_K,\otimes_{K},K)$, i.e. as an usual algebra.
\end{theorem}

\begin{proof}
It is enough to prove that if $A$ is separable as an algebra in $\mathfrak{M}%
_K$ then it is separable as an algebra in ${}_A^H\mathfrak{M}^H_A$. If $%
m:A\otimes A\rightarrow A$ is the multiplication of the algebra
$A$ in the monoidal category ${}_A^H\mathfrak{M}^H_A$, then $m$
also defines the multiplication of $A$ as an algebra in
$\mathfrak{M}_K$. Thus $F(m)=m$ has a section in
$_{A}\mathfrak{M}_{A}.$ Since $F$ is separable, then $m$ has a
section in ${}_A^H\mathfrak{M}^H_A$. Thus $A$ is separable in ${}^H \mathfrak{M}%
^H$.
\end{proof}

\begin{corollary}
\label{co:separable} Let $H$ be a semisimple and cosemisimple Hopf
algebra over a field $K$. If $A$ is an algebra in the category $^H
\mathfrak{M}^H$ then $A$ is separable as an algebra in
$^H\mathfrak{M}^H$ iff $A$ is separable as an algebra in
$(\mathfrak{M}_K,\otimes_{K},K)$.
\end{corollary}

\begin{proof}
By Lemma \ref{le:semicose}, $H$ has a non\lin{}{}zero
$ad$\lin{}{}invariant integral.
\end{proof}

\begin{proposition}
\label{pr:sectiona} Let $H$ be a semisimple and cosemisimple Hopf
algebra. If $\pi :A\rightarrow B$ is a surjective morphism of
algebras in ${}^H \mathfrak{M}^H$ such that $B$ is separable (as
an algebra in $\mathfrak{M}_K$) and the kernel of $\pi$ is
nilpotent then there is a section $\sigma:B\rightarrow A$ of $\pi$
which is a morphism of algebras in ${}^H \mathfrak{M}^H$.
\end{proposition}

\begin{proof}
By assumption $H$ is a semisimple and a cosemisimple Hopf algebra
and hence $H$ is separable and coseparable by \cite[Exercise
5.2.12 ]{DNC}.  Moreover by Corollary \ref{co:separable}, $B$ is
separable as an algebra in the category ${}^H \mathfrak{M}^H$.
Let $n$ be a natural number such that $I^{n}=0$, where $I=\mathrm{Ker\,}\pi$%
. Since $H$ is coseparable, any epimorphism in the category ${^{H}} %
\mathfrak{M}{^{H}}$ splits in ${^{H}}\mathfrak{M}{^{H}}$, see
\cite{Doi}. In particular, for every $r=1,\cdots ,n-1$ the
canonical morphism $\pi
_{r}:A/I^{r}\rightarrow A/I^{r+1}$ has a section in the category ${^{H}}%
\mathfrak{M}{^{H}}$. We can now conclude by applying Theorem
\ref{X formally smooth teo} to the algebra homomorphism $\pi
:A\rightarrow B.$
\end{proof}

Let $H$ be a semisimple and cosemisimple Hopf
algebra. In particular $H$ is separable and coseparable by \cite[Exercise 5.2.12 ]{DNC} (note that $%
H$ is necessarily finite dimensional). Hence, in the previous
proposition, we can choose $B=H$. \newline Actually, in this case,
we can relax the assumptions made on $H$. In order to do this we
first prove the following lemma.

\begin{lemma}\label{le:SepInM^H}
Let $H$ be a Hopf algebra.

a) $H$ is separable as an algebra in $\mm^H$ if and only if $H$ is
semisimple.

b) $H$ is separable as an algebra in $\hm^H$ if and only if there
is an integral $t\in H$ such that $\ep{}(t)=1$ and $\sum
t\ro{1}St\ro{3}\ot t\ro{2} =1\ot t$.
\end{lemma}

\begin{proof}
We prove first (b). The category of $(H,H)$\lin  bimodules in
$\hm^H$ is $\hmh$. Hence we need a section of the multiplication
of $H$ in $\hmh$. By the equivalence $\hmh\simeq\yd$ (see \ref{YD
category}), it results that: $$ \hmh(H,H\ot H)\simeq
\yd\left(H^{co(H)},(H\ot H)^{co(H)}\right), $$ the isomorphism
being given by the restriction to $H^{co(H)}$. Since the right
comodule structure of $H\ot H$ is defined by $H\ot\Delta$, we have
$(H\ot H)^{co(H)}=H$, where $H$ is regarded as a left module via
the multiplication of $H$, and as a left comodule via the adjoint
coaction. Hence: $$ \yd(H^{co(H)},(H\ot H)^{co(H)})\simeq\yd(K,H),
$$ In conclusion, there is an one$\,$\lin to$\,$\lin one
correspondence: $$ \hmh(H,H\ot H)\simeq\{t\in H\mid \sum
t\ro{1}St\ro{3}\ot t\ro{2} =1\ot t\quad \text{and}\quad
ht=\ep{}(h)t, \forall h\in H \}. $$ Through this bijection a
section of the multiplication corresponds to an  element $t$ such
that $\ep{}(t)=1$.

To prove (a) we first remark that the category  of $(H,H)$\lin
bimodules in $\mm^H$ is $_H\mm^H_H$. Proceeding as in the proof of
(b), but neglecting the left comodule structure, one can show that
there is a bijection between $_H\mm^H_H(H,H\ot H)$ and $\{t\in
H\mid ht=\ep{}(h)t, \forall h\in H\}$. Moreover, the set of
sections of the multiplication is bijectively equivalent with the
set of all $t$ as above such that $\ep{}(t)=1$.
\end{proof}

\begin{definition}
An integral $t$ in a Hopf algebra $H$ will be called $ad$\lin
coinvariant if $\ep{}(t)=1$ and $\sum t\ro{1}St\ro{3}\ot t\ro{2}
=1\ot t$.
\end{definition}

\begin{remark}
Note that if $H$ has an $ad$\lin coinvariant integral then $H$ is
semisimple, by Maschke Theorem.
\end{remark}

\begin{theorem}
\label{te:section} \ Let $H$ be a semisimple  Hopf algebra. Let
$\calM$ be either the monoidal category $\mm^H$ or $^H\mm^H$.
Suppose that $\pi :A\ra H$ is a morphisms of algebras in $\calM$
such that $\mathrm{Ker}\,\pi$ is the Jacobson radical $J$ of $A$
and $J$ is nilpotent.

a) Let $\calM=\mm^H$. Then $\pi:A\rightarrow A/J\simeq H$ has  a
section $\sigma$ in $\mm^H$ which is an algebra map.

b) Let $\calM={}^H\mm^H$. Assume that $H$ has an $ad$\lin
coinvariant integral and that every canonical map $A/J^{n+1}\ra
A/J^n$ splits in $\hm^H$. Then $\pi:A\rightarrow H$ has a section
$\sigma$ in $^H\mm^H$ which is an algebra map.
\end{theorem}

\begin{proof}
a) \ The Jacobson radical $J$ of $A$ is an $H$\lin%
subcomodule of $A$ since $\pi$ is a morphism of $H$\lin comodules.
Hence, for every $n>0$, $J^n$ is a subcomodule of $A$ too such
that the canonical map $A/J^{n+1}\ra A/J^n$ is $H$\lin colinear.
Furthermore, $J^n/J^{n+1}$ has a natural module structure over
$A/J\simeq H$, and with respect to this structure $J^n/J^{n+1}$ is
an object in $\mm^H_H$. Hence $J^n/J^{n+1}$ is a cofree right
comodule (i.e. $J^n/J^{n+1}\simeq V\ot H$). In particular
$J^n/J^{n+1}$ is an injective comodule, so the canonical map
$A/J^{n+1}\ra A/J^n$ has  a section in $\mm^H$. By the previous
lemma we know that $H$ is separable as an algebra in $\mm^H$,
therefore we can apply Theorem \ref{X formally smooth teo}.

b) We first remark that $J^n$ is an $(H,H)$\lin subbicomodule of
$A$ and that the canonical maps $A/J^{n+1}\ra A/J^n$ are morphisms
of bicomodules.  By the preceding lemma it results that $H$ is
separable in $^H\mm^H$, so we conclude by applying Theorem \ref{X
formally smooth teo}.
\end{proof}

\begin{corollary}\label{co:SectonBicolin}
\ Let $A$ be a Hopf algebra such that $J$, the Jacobson radical of
$A$ is a nilpotent coideal in $A$. Let $H:=A/J$, and let $\pi:A\ra
H$ be the canonical projection.

a) If $H$ is semisimple then there is an algebra morphism in
$\mm^H$ that is a section of $\pi$.

b) If $H$ has an $ad$\lin coinvariant integral and every canonical
map $A/J^{n+1}\ra A/J^n$ splits in $\hm^H$ then there is an
algebra morphism in $^H\mm^H$ that is a section of $\pi$.

c) If $H$ has an $ad$\lin coinvariant integral  and any object in
$\hmh$ is injective as an $(H,H)$\lin bicomodule (note that this
holds whenever $H$ is both semisimple and cosemisimple; e.g. when
$H$ is semisimple over a field of characteristic $0$) then there
is a section of $\pi$ as in (b).

\end{corollary}
\begin{proof}
The first two assertions follows directly from the previous
theorem, since  we can regard $A$ as an algebra in $\mm^H$ and as
an algebra in $^H\mm^H$, as $\pi$ is a morphism of bialgebras.

Let us prove (c). In view of (b) it is enough to show  that the
canonical epimorphisms $A/J^{n+1}\ra A/J^n$ splits in $^H\mm^H$.
Since $A/J^n$ is an object in $\hmh$ and the canonical epimorphism
$A/J^{n+1}\ra A/J^n$ is a morphism in $\hmh$ it follows that
$J^n/J^{n+1}\in \hmh$, so it is an injective $(H,H)$\lin
bicomodule. Therefore $A/J^{n+1}\ra A/J^n$ has a section in
$\hmh$.

If $H$ is semisimple and cosemisimple then $H$ is finite
dimensional.  By applying Lemma \ref{le:semicose} to the dual Hopf
algebra $H^\ast$ it results that $H$ has an $ad$\lin coinvariant
integral. Furthermore, $H$ is coseparable, so every bicomodule is
injective. In particular $J^n/J^{n+1}$ is injective. Thus
$A/J^{n+1}\ra A/J^n$ splits in $^H\mm^H$.
\end{proof}

\section{Splitting morphisms of bialgebras}
Let $H$ be a Hopf algebra and let $( A,m,u,\Delta ,\varepsilon)$
be a bialgebra. Motivated by the result that we obtained in
(\ref{co:SectonBicolin}) we are going to investigate those
bialgebras $A$ with the property  that there is a pair of
$K$\lin{}{}linear maps:
\begin{equation*}
\pi :A\rightarrow H\text{\qquad and\qquad }\sigma :H\rightarrow A
\end{equation*}
such that $\pi$ is a morphism of bialgebras and $\sigma$ is an
$(H,H)$\lin{}{}bicolinear algebra section of $\pi$ such that $\pi
\sigma = \mathrm{Id} _{H}$.

\begin{claim}
\label{pr:BialgebraInBimod}Our approach is based on the
observation that
such a bialgebra can be viewed in a natural way as an object $A\in{}%
_{H}^{H}\mathfrak{M}_{H}^{H}$ such that $A$ is an algebra in $({}_{H}^{H}%
\mathfrak{M}_{H}^{H},\otimes_H,H)$ and a coalgebra in $({}^H \mathfrak{M}%
^H,\square_H,H)$. Let us explain the exact meaning of this
sentence.

Since $\pi $ is a morphism of coalgebras, $A$ is an
$(H,H)$\lin{}{}bicomodule with the structures induced by $\pi $.
Similarly $\sigma $ defines an $(H,H)$\lin{}{}bimodule structure
on $A$.

Let us prove that these structures make $A$ a Hopf bimodule. We
will check that $\rho ^{l}:A\rightarrow H\otimes A$ and $\rho
^{r}:$ $A\rightarrow A\otimes H$ are bimodule morphisms. By
definition $\rho ^{l}(a)=\sum \pi (a_{(1)})\otimes a_{(2)}.$
Hence:
\begin{equation*}
\rho ^{l}(ha) =\rho ^{l}(\sigma (h)a)=\sum \pi (\sigma
(h)_{(1)}a_{(1)})\otimes \sigma (h)_{(2)}a_{(2)} =\sum \pi (\sigma
(h)_{(1)})\pi (a_{(1)})\otimes \sigma (h)_{(2)}a_{(2)},
\end{equation*}
where the last equality has been deduced by using the fact that
$\pi$ is a morphism of algebras. Thus, by the definition of the
left $H$\lin{}{}coaction on $A $ and the fact that $\sigma$ is
left $H$\lin{}{}colinear, we get:
\begin{equation*}
\rho ^{l}(ha){=}\sum \sigma (h){}_{\langle-1\rangle}a{}_{\langle
-1\rangle}\otimes \sigma (h){}_{\langle 0\rangle}a{}_{\langle 0\rangle} {=}%
\sum h_{1}a{}_{\langle-1\rangle}\otimes \sigma (h_{2})a{}_{\langle
0\rangle} {=}\sum h_{1}a{}_{\langle-1\rangle}\otimes
h_{2}a{}_{\langle 0\rangle}.
\end{equation*}
In a similar way one can prove that $\rho ^{l}$ is right
$H$\lin{}{}linear and that $\rho ^{r}$ is a morphism of bimodules.

By assumption, $\sigma $ is a morphism of
$(H,H)$\lin{}{}bicomodules. Moreover, since $\sigma $ is also an
algebra morphism, we get that $\sigma $ is morphism of
$(H,H)$\lin{}{}bimodules and $m:A\otimes A\rightarrow A$
factorizes to a map $\overline{m}:A\otimes _{H}A\rightarrow A.$
Furthermore, $m$ is left $H $\lin{}{}colinear. Indeed, we have:
\begin{equation*}
\left( H\otimes m\right) \rho _{A\otimes _{H}A}^{l}(a\otimes b)
=\sum a{}_{\langle-1\rangle}b{}_{\langle-1\rangle}\otimes
a{}_{\langle 0\rangle}b{}_{\langle 0\rangle}=\sum \pi
((ab)_{(1)})\otimes (ab)_{(2)}=\rho ^{l}(m(ab)).
\end{equation*}
In a similar way one proves that $m$ is also right
$H$\lin{}{}colinear. Clearly $m$ is also $(H,H)$\lin{}{}bilinear.
We have proved that $(A,\overline{m},\sigma )$ is an algebra in
$({}_{H}^{H}\mathfrak{M}_{H}^{H},\otimes_H,H)$.

One can easily check that the image of $\Delta $ is contained in $A{\square }%
_{H}A$. Let $\overline{\Delta }$ be the corestriction of $\Delta $ to $A{%
\square } _{H}A$. We have to prove that $\overline{\Delta }$ is a
morphism in ${^{H}} \mathfrak{M}{^{H}}$. Since the left
$H$\lin{}{}coaction on $A\square_H A $ is given by the restriction
of $\rho^l_A\otimes A$ to $A\square_H A$ it results:
\begin{equation*}
\rho _{A\square_H A}^l(\overline{\Delta}(a)) = \sum
\pi(a{}_{\left( 1\right)})\otimes a{}_{\left( 2\right)}\otimes
a{}_{\left( 3\right)}= ( H\otimes \overline{\Delta}) \rho
_{A}^{l}(a),
\end{equation*}
and similarly for the right coaction. Since $\pi $ is a coalgebra
map, we
have that $\varepsilon _{H} \pi =\varepsilon .$ Moreover $\pi $ is $(H,H)$%
\lin{}{}bicolinear. Hence $\left( A,\overline{\Delta} ,\pi
\right)$ is a coalgebra in $({^{H}}\mathfrak{M}{^{H}},{\square
}_{H},H)$.

These considerations lead us to the following definition (see also
Definition \ref{de:BecomesCoalg}).
\end{claim}

\begin{definition}
Let $R$ be an $H$\lin{}{}bicomodule algebra. Let $A$ be an object
in $_{R}^{H} \mathfrak{M}_{R}^{H}$ which is an algebra in the
category of vector spaces with multiplication $m:A\otimes
A\rightarrow A$ and unit $u:K\rightarrow A.$ We say that $(A,m,u)$
\emph{becomes an} \emph{algebra in}
$({_{R}^{H}}\mathfrak{M}{_{R}^{H},\otimes }_{R},R)$ if $(A,m,u)$
is an $H$\lin{}{}bicomodule algebra and $m$ factorizes to a
morphism $\overline{m}:A\otimes _{R}A\rightarrow A$ in
${}_{R}^{H}\mathfrak{M}_{R}^{H}.$
\end{definition}

\begin{claim}
Note that $A$ becomes an algebra in $_{R}^{H}\mathfrak{M}_{R}^{H}$
iff $A$ is an $H$\lin{}{}bicomodule algebra and $m$ is a morphism
of $(R,R)$\lin{}{}bimodules which is $R$\lin{}{}balanced. Let us
denote $m(a\otimes b)=\overline{m}(a\otimes _{R}b)$ by $ab$. Then,
for $a,b\in A$ and $r\in R$, we have
\begin{equation}
(ar) b=a(rb)\qquad ra=\left( r1_{A}\right) a\qquad ar=a\left(
1_{A}r\right), \label{acca}
\end{equation}
since by definition $m$ is an $R$\lin{}{}balanced morphism of
$(R,R)$\lin{}{}bimodules. In particular the first relation gives
us $r1_A$=$1_Ar$, for all $ r\in R$, so the unique left
$R$\lin{}{}linear map $\overline{u}:R\rightarrow A$,
$\overline{u}(r)=r1_A$ is a morphism of $(R,R)$\lin{}{}bimodules.
Since $A$ is an object in ${}_{R}^{H}\mathfrak{M}_{R}^{H}$ and $u$
is $(H,H)$\lin{}{}bicolinear one can check easily that
$\overline{u}$ is a morphism of $(H,H)$\lin{}{}bicomodules too, so
$(A,\overline{m},\overline{u})$ is an algebra in the monoidal
category $({}_{R}^{H}\mathfrak{M}_{R}^{H},\otimes_R,R)$. This
remark explains the terminology we introduced in the previous
definition.
\end{claim}

\begin{proposition}
\label{prob}Let $R$ be an $H$\lin{}{}bicomodule algebra. Let $\phi
:(A,m,u)\rightarrow \left( B,n,v\right) $ be an isomorphism of
algebras in the category of vector spaces. If $A\in
{}_{R}^{H}\mathfrak{M}_{R}^{H}$ then $B$ can be endowed, via $\phi
,$ with obvious Hopf bimodule structures and $\phi :A\rightarrow
B$ is an isomorphism in ${_{R}^{H}}\mathfrak{M}{_{R}^{H}}$ .
Moreover if $A$ becomes an algebra in
$({_{R}^{H}}\mathfrak{M}{_{R}^{H},\otimes }_{R},R)$, then $\left(
B,n,v\right) $ also becomes an algebra in
$({_{R}^{H}}\mathfrak{M}{_{R}^{H},\otimes }_{R},R)$ such that
$\phi :(A,\overline{m},\overline{u})\rightarrow
(B,\overline{n},\overline{v})$ is an algebra isomorphism in the
category $({ _{R}^{H}}\mathfrak{M}{_{R}^{H},\otimes }_{R},R)$.
\end{proposition}

\begin{proof}
Obvious.
\end{proof}

\begin{claim}\label{YD category}
It is well known that ${}_{H}^{H}\mathfrak{M}_{H}^{H}$ is a
monoidal category equivalent to the category of
Yetter\lin{}{}Drinfeld modules ${_{H}^{H} \mathcal{YD}}$. This
category consists of left $H$\lin{}{}modules $V$ which are left
$H$\lin{}{}comodules such that
\begin{equation*}
\rho \left( {}^{h}v\right) =\sum h_{\left( 1\right) }v{}_{\langle
-1\rangle}Sh_{\left( 3\right) }\otimes {}^{h_{\left( 2\right)
}}v{}_{\langle 0\rangle},\forall h\in H,v\in V,
\end{equation*}
where $^{h}v$ is the notation that we shall use for the
multiplication of $ v\in V$ by $h\in H$ in an
Yetter\lin{}{}Drinfeld module. The tensor product is $ (-)\otimes
_{K}(-)$, endowed with diagonal action and coaction. The
equivalence between ${}_{H}^{H}\mathfrak{M}_{H}^{H}$ and
${}_H^H\mathcal{YD}$ is given by $V\rightarrow V^{co\left(
H\right) },$ where
\begin{equation*}
V^{co\left( H\right) }=\left\{ v\in V\mid \rho _{V}^{r}\left(
v\right) =v\otimes 1\right\} .
\end{equation*}
The structures making $V^{co\left( H\right) }$ a left
Yetter\lin{}{}Drinfeld module are the left adjoint action and the
restriction of the left comodule structure of $V$:
\begin{eqnarray}
^{h}v &=&\sum h_{\left( 1\right) }vSh_{\left( 2\right) }
\label{eq: left YD-1} \\ \rho &=&\rho _{V}^{l}\rest{{V^{co\left(
H\right)}}}.  \label{eq: left YD-2}
\end{eqnarray}
Conversely if $W\in {_{H}^{H}\mathcal{YD}}$, then $W\otimes H$
becomes an object in ${_{H}^{H}}\mathfrak{M}{_{H}^{H}}$ with the
canonical right structures (coming from $H$) and with diagonal
left action and coaction:
\begin{eqnarray*}
{}^h\left( w\otimes k\right) &=&\sum {}^{h_{\left( 1\right)
}}w\otimes h_{\left( 2\right) }k,\ \forall w\in W,\ \forall h,k\in
H \\ \rho ^{l}\left( w\otimes h\right) &=&\sum
w{}_{\langle-1\rangle}h_{\left( 1\right) }\otimes w{}_{\langle
0\rangle}\otimes h_{\left( 2\right) },\ \forall w\in W,\ \forall
h\in H
\end{eqnarray*}
The functor that associates to $W\in {{}_{H}^{H}\mathcal{YD}}$ the
Hopf bimodule $W\otimes H$ is an inverse of the monoidal functor
defined above. Let us remark that, for every
$V\in{}_{H}^{H}\mathfrak{M}_{H}^{H}$, the isomorphism in
${}_{H}^{H}\mathfrak{M}_{H}^{H}$ between $V^{Co(H)}\otimes H$ and
$V$ is given by:
\begin{equation}  \label{eq:fiV}
\phi_V:V^{Co(H)}\otimes H\rightarrow V,\quad \phi_V(v\otimes
h)=vh.
\end{equation}
\end{claim}

\begin{example}
\label{Gold claim}Let $R$ be a left $H$\lin{}{}module algebra.
Recall that the \emph{smash product} $R\#H$ of $R$ and $H$ is the
associative algebra defined on $R\otimes H$ by setting:
\begin{equation}
\left( r\#h\right) \left( s\#k\right) =\sum r\left( {}^{h_{\left(
1\right) }}s\right) \#h_{\left( 2\right) }k.  \label{eq:Smash}
\end{equation}
This algebra is unitary, with unit $1_{R}\#1_{H}.$ Here
$r\#h:=r\otimes h.$
Moreover, if we assume in addition that $R$ is an algebra in ${_{H}^{H}%
\mathcal{YD}}$ then $R\#H$ becomes an algebra in $({_{H}^{H}}%
\mathfrak{M}{_{H}^{H},\otimes }_{H},H){,}$ with respect to the
structures:
\begin{equation*}
\begin{tabular}{ll}
$\rho _{R\#H}^{l}\left( r\#h\right) =\sum
r_{\left\langle-1\right\rangle }h_{(1)}\otimes \left(
r_{\left\langle 0\right\rangle }\#h_{(2)}\right) $ & $ \rho
_{R\#H}^{r}\left( r\#h\right) =\sum \left( r\#h_{(1)}\right)
\otimes h_{(2)}$ \\ $h\left( r\#k\right) =\sum {}^{h_{\left(
1\right) }}r\#h_{\left( 2\right) }k$
& $\left( s\#k\right) h=s\#kh.$%
\end{tabular}
\end{equation*}
Our next aim is to prove that \emph{any} algebra $A$ that becomes
an algebra in $({_{H}^{H}}\mathfrak{M}{_{H}^{H},\otimes }_{H},H)$
is of this type, i.e. there is an algebra $R$ in
$_{H}^{H}\mathcal{YD}$ such that $A\simeq R\#H.$
\end{example}

\begin{definition}
Let $V$ be a Hopf bimodule. The space of right coinvariant
elements of $V$ will be called the \emph{diagram} of $V$ and it
will be denoted by $R_V$, or shortly by $R$ if there is no danger
of confusion.
\end{definition}

\begin{proposition}
\label{pr:AlgInBi} \label{smash}Let $(A,m,u)$ be an algebra.
Suppose that $A$ is an object in ${}_{H}^{H}\mathfrak{M}_{H}^{H}$
such that $A$ becomes an algebra in
$({}_{H}^{H}\mathfrak{M}_{H}^{H}, \otimes_H,H)$. If $R=A^{Co\left(
H\right) }$ is the diagram of $A$ then $R$ is an algebra in $
{}_H^H\mathcal{YD}$ and the canonical isomorphism $\phi_A :R\#
H\rightarrow A $ is a morphism of algebras in
$({}_{H}^{H}\mathfrak{M}_{H}^{H}, \otimes_H,H) $.
\end{proposition}

\begin{proof}
Since $A$ becomes an algebra in
$({}_{H}^{H}\mathfrak{M}_{H}^{H},\otimes_H,H) $ the multiplication
$m$ factors to a map $\overline{m}:A\otimes_H A\rightarrow A$
which is a morphism of Hopf bimodules. As $\left(-\right)
^{Co\left( H\right) }:{}_{H}^{H}\mathfrak{M}_{H}^{H}\rightarrow{}_H^H%
\mathcal{YD}$ is a monoidal equivalence, we have
\begin{equation*}
(A{\otimes }_{H}A)^{Co\left( H\right) }\simeq A^{Co\left( H\right) } {\otimes%
}A^{Co\left( H\right). }
\end{equation*}
The morphisms of Yetter\lin{}{}Drinfeld modules that corresponds
$\overline{m}$ through this equivalence is the restriction of $m$
to $R\otimes R$. Therefore we shall denoted by $m:R\otimes
R\rightarrow R$ too. Obviously $m$ defines an associative
multiplication on $R$.

Let $\overline{u}:H\rightarrow A$ be the unique left
$H$\lin{}{}linear map such that $\overline{u}(1_H)=u(1_K)=1_A$. By
assumption $\overline{u}$ is a
morphism in ${}_{H}^{H}\mathfrak{M}_{H}^{H}$ and $(A,\overline{m},\overline{u%
})$ is an algebra in
$({}_{H}^{H}\mathfrak{M}_{H}^{H},\otimes_H,H)$. In particular
$1_A$ is right coinvariant, so we can regard $u$ as a
$K$\lin{}{}linear map from $K$ to $R$. As $u:K\rightarrow R$ is
the morphism of Yetter\lin{}{}Drinfeld modules that corresponds to
$\overline{u}$ through the
equivalence $(-)^{Co(H)}$ we deduce that $(R,m,u)$ is an algebra in ${}_H^H%
\mathcal{YD}$.

It remains to prove that $\phi_A$ is an isomorphism of algebras in $%
({}_{H}^{H}\mathfrak{M}_{H}^{H},\otimes_H,H)$, see (\ref{eq:fiV})
for the definition of $\phi_A$. By \cite{Sch1} we know that
$\phi_A$ is a bijective morphism in
${}_{H}^{H}\mathfrak{M}_{H}^{H}$. Therefore it is enough to see
that $\phi_A$ is a morphism of algebras. We have:
\begin{equation*}
\phi_A ( ( r\#h)( s\#k) ) {=}\sum [r ( ^{h_{( 1) }}s)]( h_{( 2)
}k)= \sum r[ ( h_{( 1) }s Sh_{( 2) })( h_{( 3) }k)] =r[ ( hs)k]=(r
h)( sk).
\end{equation*}
To deduce the first equality we used the definitions of the
multiplication in $R\# H$ and of $\phi_A$, and the fact that
$\overline{m}:A\otimes _{H}A\rightarrow A$ is right
$H$\lin{}{}linear. The second equality comes from the definition
of the left module structure on $R$, see (\ref{eq: left YD-1}). To
obtain the last equalities we applied the associativity relations
(\ref {acca}) and the fact $m$ is $H$\lin{}{}balanced. We conclude
the proof of the proposition by remarking that $\phi (
r\#h)\phi(s\#k)=(r h)( sk)$.
\end{proof}

\begin{claim}
\label{cl:D^H_M_H^D}We can dualize all construction above. In
particular, the category $_H\mathfrak{M}_H$ of
$H$\lin{}{}bimodules is a monoidal category with the tensor
product $\left(-\right) \otimes _{K}\left(-\right) $. A coalgebra
$D $ in this category is a coalgebra which is an $(H,H)$-bimodule
such that $D$ is a left and a right $H$\lin{}{}module coalgebra,
i.e. an $H$\lin{}{}bimodule
coalgebra. For such a coalgebra we can consider the category ${_{H}^{D}}%
\mathfrak{M}{_{H}^{D}}$ of all $(D,D)$\lin{}{}bicomodules which
are also $(H,H)$\lin{}{}bimodules such that $\mu ^{l}:H\otimes
M\rightarrow M$ and $\mu ^{r}:M\otimes H\rightarrow M$ are
morphism of bicomodules ($\mu ^{l}$ and $\mu ^{r}$ define the
module structures on $M$ and $H\otimes M$ is an
$(D,D)$\lin{}{}bicomodule
with the diagonal coactions). For $D=K$ we get ${}_{H}^{D}\mathfrak{M}%
_{H}^{D}={}_H\mathfrak{M}_H$, and for $H=K$ we have ${}_{H}^{D}\mathfrak{M}%
_{H}^{D}={}^D\mathfrak{M}^D$.

${_{H}^{D}}\mathfrak{M}{_{H}^{D}}$ is a monoidal category with
respect to the tensor product given by $\left(-\right) {\square
}_{D}\left(-\right) ,$
the cotensor product of two $(D,D)$\lin{}{}bicomodules. If $V,W\in {_{H}^{D}}%
\mathfrak{M}{_{H}^{D}}$, then $V{\square }_{D}W$ is an
$(H,H)$\lin{}{}bimodule with diagonal actions, and its comodule
structures are defined by $\rho _{V}^{l}\otimes W$ and $V\otimes
\rho _{W}^{r}$

The monoidal category $({_{H}^{H}}\mathfrak{M}{_{H}^{H},\square
}_{H},H) $ is also monoidal equivalent to
${_{H}^{H}\mathcal{YD}}$. We are not going into details, the
reader can find them in \cite{Sch1}.
\end{claim}

\begin{definition}
\label{de:BecomesCoalg} \label{cl:epsR} Let $H$ be a Hopf algebra
and let $D$ be an $(H,H)$\lin{}{}bimodule coalgebra. Let $\left(
C,\Delta_C ,\varepsilon_C
\right) $ be a coalgebra and assume that $C\in {_{H}^{D}}%
\mathfrak{M}{_{H}^{D}}$. We say that $(C,\Delta_C ,\varepsilon_C )$ \emph{%
becomes a coalgebra} in
$({_{H}^{D}}\mathfrak{M}{_{H}^{D}},{\square }_{D},D)$ iff $\left(
C,\Delta_C ,\varepsilon_C \right) $ is an $(H,H)$\lin{}{}bimodule
coalgebra, the image of $\Delta_C $ is included in $C{\square }_{D}C$ and $%
\overline{\Delta}_C:C\rightarrow C{\square }_{D}C$ (the corestriction of $%
\Delta_C $) is a morphism in ${_{H}^{D}}\mathfrak{M}{_{H}^{D}}$.
\end{definition}

\begin{claim}
\label{cl:BecomesCoalg} Note that $(C,\Delta_C ,\varepsilon_C )$
becomes a coalgebra in ${}_{H}^{D}\mathfrak{M}_{H}^{D}$ iff
$\left( C,\Delta_C ,\varepsilon_C \right) $ is an
$(H,H)$\lin{}{}bimodule coalgebra and $\Delta_C $ is a morphism of
$(D,D)$\lin{}{}bicomodules such that $(\mathrm{Im})\Delta_C
\subseteq C{\square }_{D}C$.

Since $\mathrm{Im}(\Delta_C)\subseteq C\square_D C$, for every
$c\in C$, we have :
\begin{equation*}
\sum\left(c{}_{\left( 1\right)}\right){}_{\langle
0\rangle}\otimes\left(c{}_{\left( 1\right)}\right){}_{\langle
1\rangle}\otimes c{}_{\left( 2\right)}=\sum c{}_{\left(
1\right)}\otimes\left(c{}_{\left( 2\right)}\right){}_{\langle
-1\rangle}\otimes\left(c{}_{\left( 2\right)}\right){}_{\langle
0\rangle},
\end{equation*}
hence by applying $\varepsilon_C$ to the first and third factor we
get:
\begin{equation}  \label{eq:EpsBar}
\sum \varepsilon_{C}(c{}_{\langle 0\rangle})c{}_{\langle
1\rangle}=\sum \varepsilon_{C}(c{}_{\langle
0\rangle})c{}_{\langle-1\rangle}.
\end{equation}
In particular this relation means that the $K$\lin{}{}linear map $\overline{%
\varepsilon}_C:C\rightarrow D$, $\overline{\varepsilon}_C(c)=\sum {%
\varepsilon}_{C}(c{}_{\langle 0\rangle})c{}_{\langle-1\rangle}$ is
a morphism of $(D,D)$\lin{}{}comodules. In fact one can check
easily that it is a morphism of $(H,H)$\lin{}{}bimodules too,
since $C$ is a bimodule coalgebra over $H$ . Thus
$\overline{\varepsilon}_{C}$ is a morphism in
${}_{H}^{D}\mathfrak{M} _{H}^{D}$ such that
\begin{equation*}
\varepsilon_C =\varepsilon _{D} \overline{\varepsilon}_C
\end{equation*}
In this case $(C,\overline{\Delta }_C,\overline{\varepsilon }_C)$
is a coalgebra in{\ }$({_{H}^{D}}\mathfrak{M}{_{H}^{D}},{\square
}_{D},D)$. Indeed, $\overline{\Delta }_C$ is obviously
coassociative, and one can check easily that that the squares in
the following diagram are commutative
\begin{equation*}
\begin{diagram}[small] C{\square }_{D}C&\lTo^{\overline{\Delta }_C}&C&
\rTo^{\overline{\Delta }_C}&C{\square }_{D}C\\
\dTo<{\overline{\varepsilon }_C{\square }_{D}C}& &\dEqto&
&\dTo>{C{\square }_{D}\overline{\varepsilon }_C} \\ D{\square
}_{D}C&\lTo<{\lambda^l}&C&\rTo<{\lambda^r}&C{\square }_{D}D \\
\end{diagram}
\end{equation*}
where $\lambda ^{l}:D\square _{D}C\rightarrow C$ $\lambda
^{r}:C\square _{D}D\rightarrow C$ are the canonical isomorphisms.
\end{claim}

\begin{example}
\label{Gold claim1}Let $C$ be a left $H$\lin{}{}comodule
coalgebra. Recall that the \emph{smash coproduct} $C\#H$ of $C$
and $H$ is the coassociative and counitary coalgebra defined on
$C\otimes H$ by setting:
\begin{eqnarray}
\Delta (c\#h) &=&\sum c_{( 1) }\# ( c_{( 2) })
{}_{\langle-1\rangle}h_{( 1) }\otimes (c_{(2)}) {}_{\langle
0\rangle}\# h_{( 2)}  \label{eq:Cosmash1} \\ \varepsilon (c\#h)
&=&\varepsilon_C (c) \varepsilon_H (h) \label{eq:Cosmash2}
\end{eqnarray}
Moreover, if we assume in addition that $C$ is a coalgebra in ${}_H^H%
\mathcal{YD}$ then $C\#H$ becomes a coalgebra in $({}_{H}^{H}\mathfrak{M}%
_{H}^{H},\square_H,H)$ with respect to the left diagonal action
and diagonal coaction, and right canonical structures (coming from
the corresponding structures of $H$).

The following Proposition is the dual statement of
(\ref{pr:AlgInBi}).
\end{example}

\begin{proposition}
\label{cosmash}Let $(C,\Delta_C,\varepsilon_C)$ be a coalgebra.
Suppose that $C$ is an object in ${}_{H}^{H}\mathfrak{M}_{H}^{H}$
such that $C$ becomes a coalgebra in
$({}_{H}^{H}\mathfrak{M}_{H}^{H},\square_H,H)$. If $R=C^{Co(H)}$
is the diagram of $C$ then $R$ is a coalgebra in
${}_H^H\mathcal{YD}$ and the canonical isomorphism $\phi_C :R\#
H\rightarrow C$ is a morphism of coalgebras in
$({}_{H}^{H}\mathfrak{M}_{H}^{H},\square_H,H)$.
\end{proposition}

\begin{remark}
We keep the notation from the previous proposition. Let us denote
the
comultiplication of $R$ by $\delta :R\rightarrow R\otimes R$. Then $%
\delta=\psi_C \Delta\rest{R},$ where $\psi_C $ is the isomorphism
\begin{equation*}
\left( C\square _{H}C\right) ^{Co\left( H\right) }\simeq
C^{Co\left( H\right) }\otimes C^{Co\left( H\right) }
\end{equation*}
that comes from the equivalence $({}_{H}^{H}\mathfrak{M}_{H}^{H},%
\square_{H},H)\simeq \left( {}_H^H\mathcal{YD},\otimes, K \right)
. $ More generally, the isomorphism $\left( V\square _{H}W\right)
^{Co\left( H\right) }\simeq V^{Co\left( H\right) }\otimes
W^{Co\left( H\right) }$ maps an element
$\sum_{i=1}^{n}v_{i}\otimes w_{i}$ to $\sum_{i=1}^{n}v_{i}S\left(
w_{i{}_{\langle-1\rangle}}\right) \otimes w_{i{}_{\langle
0\rangle}}.$
\end{remark}

\begin{claim}
\label{VW} Suppose that $H$ is a Hopf algebra. Let $V\in\mathfrak{M}_K$ and $%
W\in{}^H \mathfrak{M}$. It is well\lin{}{}known that we have a
functorial isomorphism:
\begin{equation}  \label{eq:IzoCotensor}
(V\otimes H)\square_H W\simeq V\otimes W
\end{equation}
which is given by $\sum_{i=1}^{n}v_{i}\square_H h_i \otimes
w_{i}\mapsto
\sum_{i=1}^{n}\varepsilon(h_i)v_i\otimes w_1$. The inverse of this map is $%
V\otimes\rho_W$.

Furthermore, the functor $F:\mathfrak{M}_K\rightarrow\mathfrak{M}^H$, $%
F(V)=V\otimes H$, has as a left adjoint the functor $G:\mathfrak{M}%
^H\rightarrow\mathfrak{M}_K$ that ``forgets'' the comodule
structure. The maps that define this adjunction are:
\begin{eqnarray*}
\alpha_{V,W}:\mathfrak{M}^H(V,W\otimes H)\rightarrow\mathfrak{M}%
_K(V,W),&&\alpha_{V,W}(f)=(V\otimes\varepsilon) f \\
\beta_{V,W}:\mathfrak{M}_K(V,W)\rightarrow\mathfrak{M}^H,&&\beta_{V,W}(g)=(g%
\otimes H) \rho_V
\end{eqnarray*}
where $\rho_V$ defines the comodule structure on $V$.
\end{claim}

\begin{lemma}
\label{gamma}%
Let  $V$, $W$ be two vector spaces and let $Z$ be a left $H$\lin
comodule. If we regard $Z\ot H$ as a left $H$\lin comodule with
diagonal coaction and a right comodule via $Z\ot\Delta_H$ then
there is an one$\,$\lin to$\,$\lin one
correspondence between $\mathfrak{M}_K(V\otimes H,W\otimes Z)$ and $%
\mathfrak{M}^H(V\otimes H,(W\otimes H)\square_H(Z\otimes H))$. If $\gamma\in%
\mathfrak{M}_K(V\otimes H,W\otimes Z)$ and
$\Gamma\in\mathfrak{M}^H(V\otimes H,(W\otimes H)\square_H(Z\otimes
H))$ correspond to each other through this bijective map then they
are related by the following relations:
\begin{eqnarray}
&\widetilde{\gamma} =(W\otimes \varepsilon_H\otimes
Z\otimes\varepsilon_H)\Gamma,  \label{eq:gama} \\ &\Gamma(v\otimes
h)=\sum\widetilde{\gamma}^1(v\otimes h{}_{\left(
1\right)})\otimes\widetilde{\gamma}^2(v\otimes h{}_{\left(
1\right)}){}_{\langle-1\rangle}h{}_{\left( 2\right)}\otimes\widetilde{\gamma%
}^2(v\otimes h{}_{\left( 1\right)}){}_{\langle 0\rangle}\otimes
h{}_{\left( 3\right)},  \label{G}
\end{eqnarray}
where $\widetilde{\gamma}(v\otimes
h)=\sum\widetilde{\gamma}^1(v\otimes
h)\otimes\widetilde{\gamma}^2(v\otimes h)\in W\otimes Z$ is a
Sweedler\lin{}{}like notation for $\widetilde{\gamma}(v\otimes
h)$.
\end{lemma}

\begin{proof}
By (\ref{eq:IzoCotensor}) we have:
\begin{equation*}
(W\otimes H)\square_H (Z\otimes H)\simeq W\otimes Z\otimes H.
\end{equation*}
Hence
\begin{equation*}
\mathfrak{M}^H(V\otimes H, (W\otimes H)\square_H (Z\otimes H))\simeq %
\mathfrak{M}^H(V\otimes H, W\otimes Z\otimes H).
\end{equation*}
By composing this isomorphism with $\alpha_{V\otimes H,W\otimes
Z}$ we obtain a bijective map:
\begin{equation*}
\mathfrak{M}^H(V\otimes H, (W\otimes H)\square_H (Z\otimes H))\rightarrow%
\mathfrak{M}_K(V\otimes H,W\otimes Z).
\end{equation*}
Suppose now that $\widetilde{\gamma}$ and $\Gamma$ correspond each
other through the above $K$\lin{}{}linear isomorphism. A
straightforward but tedious computation shows us that
$\widetilde{\gamma}$ and $\Gamma$ verifies (\ref {eq:gama}) and
(\ref{G}).
\end{proof}

\begin{claim}
\label{gammac}Let $R\in \hm$ and let $\Delta_{R\#H}
:(R\#H)\rightarrow (R\#H){\square }_{H}(R\#H)$ be a right
$H$\lin{}{}colinear map. By the previous lemma if
\begin{equation}
\tilde{\delta}=(R\otimes \varepsilon _{H}\otimes R\otimes
\varepsilon _{H})\Delta_{R\#H},  \label{eq:delta}
\end{equation}
and for $r\in R$, $h\in H$ we write $\tilde{\delta}(r\# h)=\sum \tilde{\delta%
}^{1}(r\# h)\otimes \tilde{\delta}^{2}(r\# h)\in R\otimes R$,
then:
\begin{equation}
\Delta_{R\#H} (r\# h)=\sum \tilde{\delta}^{1}(r\# h_{(1)})\# \tilde{\delta}%
^{2}(r\# h_{\left( 1\right) }){}_{\langle-1\rangle}h_{(2)}\otimes \tilde{%
\delta}^{2}(r\# h_{(1)}){}_{\langle 0\rangle}\# h_{3}.
\label{eqdd}
\end{equation}
Conversely if $\tilde{\delta}:R\otimes H\rightarrow R\otimes R$ is
a linear
map and $\Delta_{R\#H}$ is defined by $\left( \ref{eqdd}\right) $, then $%
\Delta_{R\#H}$ is a right $H$\lin{}{}colinear map and $\mathrm{Im}%
(\Delta_{R\#H})\subseteq(R\#H)\square_H(R\#H)$.
\end{claim}

\begin{claim}
\label{cl:Settings} Let $A\in{}_{H}^{H}\mathfrak{M}_{H}^{H}$ be a
Hopf
bimodule. We assume that $A$ is a bialgebra with multiplication $m$, unit $u$%
, comultiplication $\Delta$ and counit $\varepsilon$. We are
looking for
necessary and sufficient conditions such that $A$ becomes an algebra in $%
({}_{H}^{H}\mathfrak{M}_{H}^{H},\otimes_H,H)$ and a coalgebra in $(^H%
\mathfrak{M}^H,\square_H,H)$. Note that the latter monoidal
category is equal to $(^H_K \mathfrak{M}^H_K,\square_H,H)$, see
(\ref{cl:D^H_M_H^D}). Thus it makes sense to talk about a
bialgebra that becomes a coalgebra in the category of
$(H,H)$\lin{}{}bicomodules, see (\ref{de:BecomesCoalg}).

By Proposition \ref{smash} the diagram $(R,m,u)$ of $A$ is an algebra in $%
{}_H^H\mathcal{YD}$, $R\#H$ is an algebra and the map $\phi_A
:R\#H\rightarrow A$, $\phi_A(r\otimes h)=rh$, is an isomorphism of
algebras in $({}_{H}^{H}\mathfrak{M}_{H}^{H},\otimes_H,H)$.
Obviously, $R\#H$ is a
bialgebra with comultiplication $\Delta_{R\#H}$ and counit $%
\varepsilon_{R\#H}$ given by:
\begin{equation*}
\Delta_{R\#H}:=(\phi_A^{-1}\otimes\phi_A^{-1})\Delta\phi_A \qquad\text{and}%
\qquad\varepsilon_{R\#H}:=\varepsilon_{}\phi_A.
\end{equation*}
Of course, with respect to this bialgebra structure, $\phi_A$
becomes an isomorphism of bialgebras.

Furthermore, since $A$ becomes an algebra in $({}_{H}^{H}\mathfrak{M}%
_{H}^{H},\otimes_H,H)$ and a coalgebra in ${}^H \mathfrak{M}^H$, the smash $%
R\#H$ has the same properties. In particular the image of
$\Delta_{R\#H}$ is included in $\Delta_{R\#H}\square_H
\Delta_{R\#H}$ and $\Delta_{R\#H}$ can be regarded as a morphism
of right $H$\lin{}{}comodules $\Delta_{R\#H} :R\#H\rightarrow
(R\#H){\ \square }_{H}(R\#H)$.

Hence, by Corollary \ref{gammac}, $\Delta_{R\#H}$ is uniquely
determined by the $K$\lin{}{}linear map
$\widetilde{\delta}:R\#H\rightarrow R\otimes R$. In order to
determine the counit $\varepsilon_{R\#H}$ we consider the
restriction of $\varepsilon_{}$ to $R$. For simplifying the
notation we shall denote it by $\varepsilon_{}$ too.
\end{claim}

\begin{lemma}
We have
$\varepsilon_{R\#H}(r\#h)=\varepsilon_{}(r)\varepsilon_{H}(h)$,
for
all $r\in R$ and $h\in H$, if and only if $\varepsilon_{}(1_Ah)=%
\varepsilon_{H}(h)$, for all $h\in H$ (equivalently,
$\varepsilon_{}$ is right $H$\lin{}{}linear).
\end{lemma}

\begin{proof}
Let us assume that $\varepsilon_{}(1_Ah)=\varepsilon_{H}(h)$, for
all $h\in H
$. By definition and the middle relation in (\ref{acca}) we have: $%
\varepsilon_{R\#H}(r\#h)=\varepsilon_{}(rh)=\varepsilon_{}(r(1_Ah))=%
\varepsilon_{}(r)\varepsilon_{}(1_Ah)=\varepsilon_{}(r)\varepsilon_{H}(h).$

The other implication is trivial since $\varepsilon_{H}(h)=\varepsilon_{R%
\#H}(1_A\#h)$.
\end{proof}

\noindent All considerations above still hold if we work with an
arbitrary algebra $R$ in ${}_H^H\mathcal{YD}$. To be more precise
we reformulate our problem of characterizing algebras $A$ as above
with the additional property that $\varepsilon$ is right
$H$\lin{}{}linear) in the following way.

\begin{problem}
Let $R$ be an algebra in ${}_H^H\mathcal{YD}$. Suppose that $\widetilde{%
\delta}:R\#H\rightarrow R\otimes R$, $\varepsilon_{}:R\rightarrow K$ are $K$%
\lin{}{}linear maps. Let $\Delta_{R\#H}$ be defined by (\ref{eqdd}) and let $%
\varepsilon_{R\#H}:=\varepsilon_{}\otimes\varepsilon_{H}$. Find
necessary and sufficient condition such that
$(R\#H,\Delta_{R\#H},\varepsilon_{R\#H})$ is a bialgebra that
becomes a coalgebra in $({}^H \mathfrak{M}^H,\square_H,H) $.
\end{problem}

\noindent Note that $R\#H$ always becomes an algebra in $({}_{H}^{H}%
\mathfrak{M}_{H}^{H},\otimes_H,H)$. Of course, by solving the
above problem we also get an answer to our initial question (of
finding all bialgebras $A$ that become an algebra in
${}_{H}^{H}\mathfrak{M}_{H}^{H}$ and a coalgebra in ${}^H
\mathfrak{M}^H$). It is enough to take $R$ to be the diagram of
$A$ and $\widetilde{\delta}$, $\varepsilon_{}$ as in
(\ref{cl:Settings}). Therefore, throughout remaining part of this
section we will keep (if not otherwise stated) the following
\textbf{notation}:

$\bullet$ \ $R$ is an algebra in ${}_H^H\mathcal{YD}$;

$\bullet$ \ $\widetilde{ \delta}:R\#H\rightarrow R\otimes R$ and
$\varepsilon_{}:R\rightarrow K$ are $K$\lin{}{}linear maps;

$\bullet$ \  $\Delta_{R\#H}$ is defined by (\ref{eqdd});

$\bullet$ \
$\varepsilon_{R\#H}:=\varepsilon_{}\otimes\varepsilon_{H}$.

$\bullet$ \ $\rho_{\sm}:\sm\ra H\otimes \sm$ denotes the map that
defines the left coaction on $\sm$, coming from the monoidal
structure of $\yd$.

\begin{claim}
To simplify the computation sometimes we shall use the method of
representing morphisms in a braided category by diagrams. The
reader is referred to for details to \cite[Chapter XIV.1]{Ka}.
Here we shall only mention that the morphisms are represented by
arrows oriented downwards.

We shall apply this method in the category ${}_H^H\mathcal{YD}$ of
Yetter\lin{}{}Drinfeld modules. Recall that, for every
$V,W\in{}_H^H\mathcal{YD}$ the braiding is given by:
\begin{equation}  \label{eq:Braiding}
c_{V,W}:V\otimes W\rightarrow W\otimes V\text{\qquad
}c_{V,W}(v\otimes w)=\sum v{}_{\langle-1\rangle}w\otimes
v{}_{\langle 0\rangle}.
\end{equation}
Two examples of diagrams in this category can be found in
Figure~\ref {fig:RhoDelta}. Note that in both pictures the
crossings represent $c_{R,H}$.
\begin{figure}[bt]
\begin{center}
\includegraphics{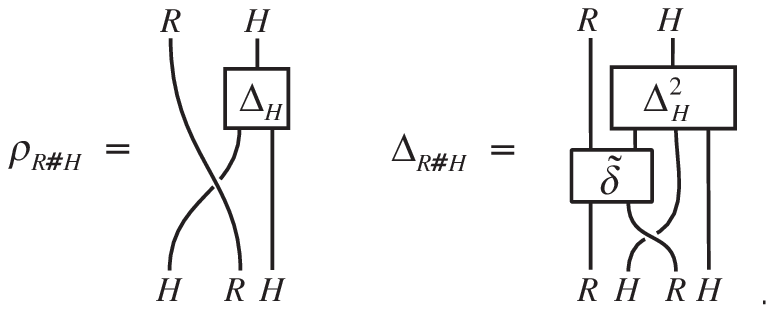}
\end{center}
\par
\vspace*{-0.5cm} \caption{Definitions of $\protect\rho _{R\#H}$
and $\Delta _{R\#H}.$} \label{fig:RhoDelta}
\end{figure}
Throughout the remaining part of this section we shall keep the
notation and the assumptions of this paragraph.
\end{claim}

\begin{lemma}
Let $H$ be a Hopf algebra. Then:

a) $(\varepsilon _{H}\otimes R) c_{R,H}=c_{R,K}(R\otimes
\varepsilon _{H}).$

b) $(\Delta _{H}\otimes R) c_{R,H}=(H\otimes c_{R,H})
(c_{R,H}\otimes H) (R\otimes \Delta _{H}).$
\end{lemma}

\begin{proof}
Trivial.
\end{proof}

\begin{remark}
\label{re:PropEpsDelta} The equations from the previous lemma
admits the representations from Figure~\ref{fig:PropEpsDelta}.
\begin{figure}[tb]
\begin{center}
\includegraphics{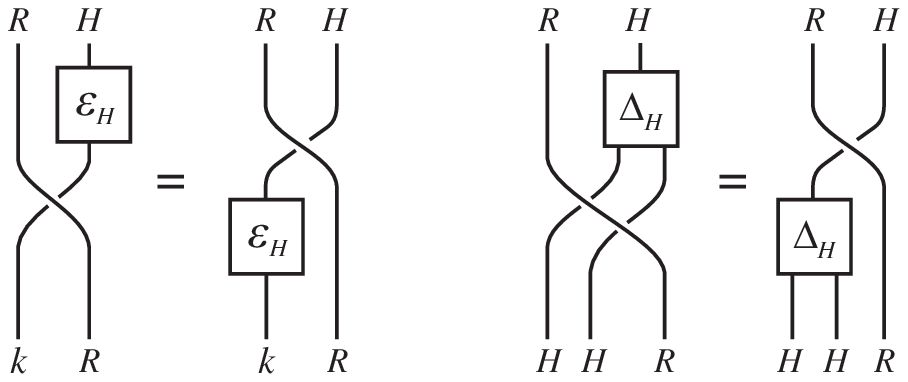}
\caption{Properties of $\protect\varepsilon _{H}$ and $\Delta
_{H}.$} \label{fig:PropEpsDelta}
\end{center}
\end{figure}
\end{remark}

\begin{lemma}
\label{Lemma 03}The following two relations are equivalent.
\begin{eqnarray}
&\lbrack \rho _{R\otimes H}\otimes (R\#H)]\Delta _{R\#H}
=(H\otimes \Delta _{R\#H})\rho _{R\#H}  \label{ec:DeltaColin} \\
&(H\otimes \tilde{\delta})\rho_{R\#H} =\left(c_{R,H}\otimes
R\right) \left( R\otimes c_{R,H}\right) ( \tilde{\delta}\otimes H)
(R\otimes \Delta _{H}) \label{eqdel1}
\end{eqnarray}
\end{lemma}

\begin{proof}
Note that (\ref{ec:DeltaColin}) means that $\Delta _{R\#H}$ is left $H$%
\lin{}{}colinear, and that the equivalence that we have to prove
can be represented as in Figure~\ref{fig:DeltaColin1}.
\begin{figure}[htb]
\begin{center}
\includegraphics{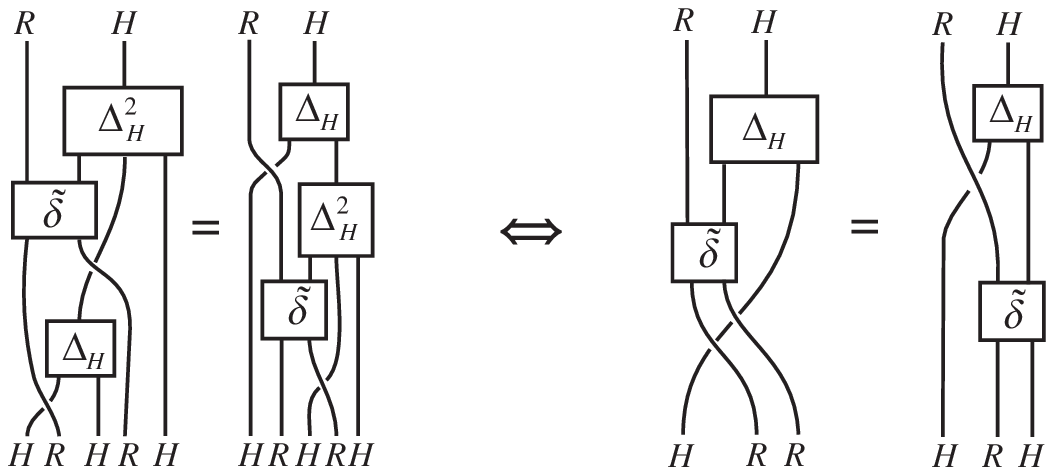}
\end{center}
\par
\vspace*{-0.5cm} \caption{Representation of
$($\ref{ec:DeltaColin}$)\Longleftrightarrow ($\ref {eqdel1}$)$}
\label{fig:DeltaColin1}
\end{figure}
We prove that (\ref{ec:DeltaColin})$\Rightarrow $(\ref{eqdel1}) in Figure~%
\ref{fig:DeltaLColin2}. The first equality there was obtained by
composing with $H\otimes R\otimes \varepsilon _{H}\otimes R\otimes
\varepsilon _{H}$
both sides of (\ref{ec:DeltaColin}). The second equation holds because $%
\varepsilon _{H}$ and $\Delta _{H}$ can be pulled under the string
in a crossing, see Remark \ref{re:PropEpsDelta}. We conclude the
proof of this implication by using that $\varepsilon _{H}$ is the
counit of $H.$

\begin{figure}[tb]
\begin{center}
\includegraphics{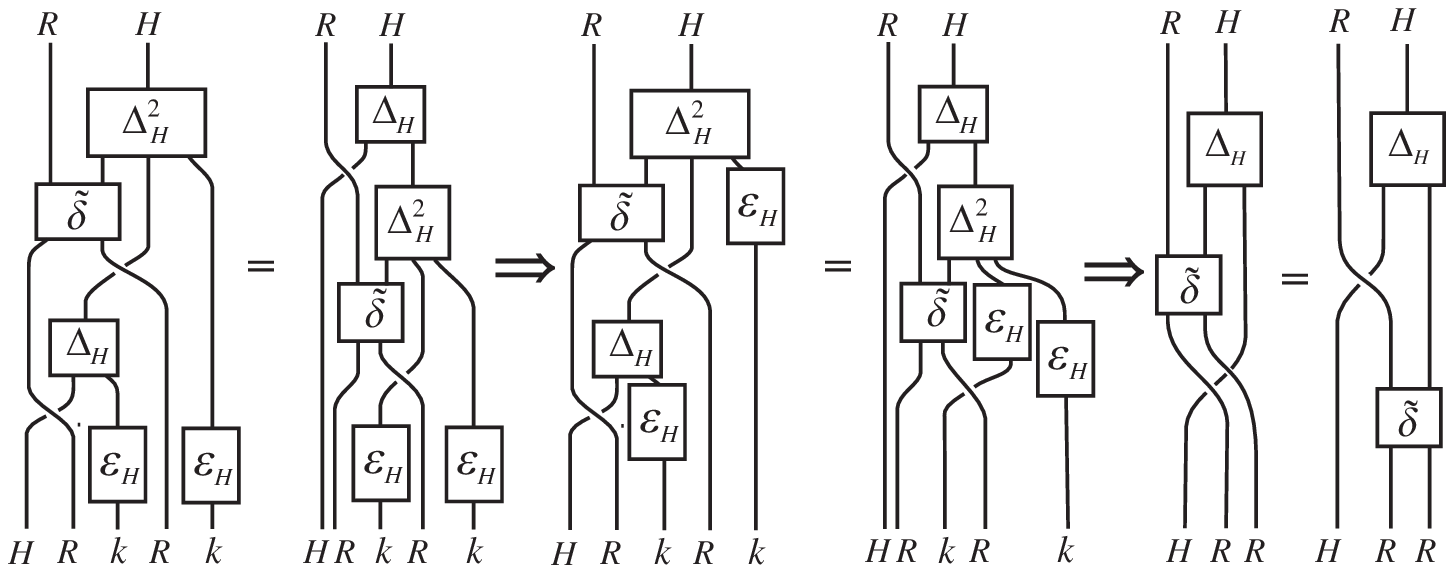}
\end{center}
\par
\vspace*{-0.5cm}
\caption{The proof of $($\ref{ec:DeltaColin}$)\Longrightarrow ($\ref{eqdel1}$%
)$} \label{fig:DeltaLColin2}
\end{figure}
The other implication is proved in Figure~\ref{fig:DeltaLColin3}.
By Remark \ref{re:PropEpsDelta} we can drag $\Delta _{H}$ under
the braiding, so we get the first equality. Since the
comultiplication in $H$ is coassociative we have the second and
last relations. The third one follows since, by assumption,
(\ref{eqdel1}) holds.
\begin{figure}[tb]
\begin{center}
\includegraphics{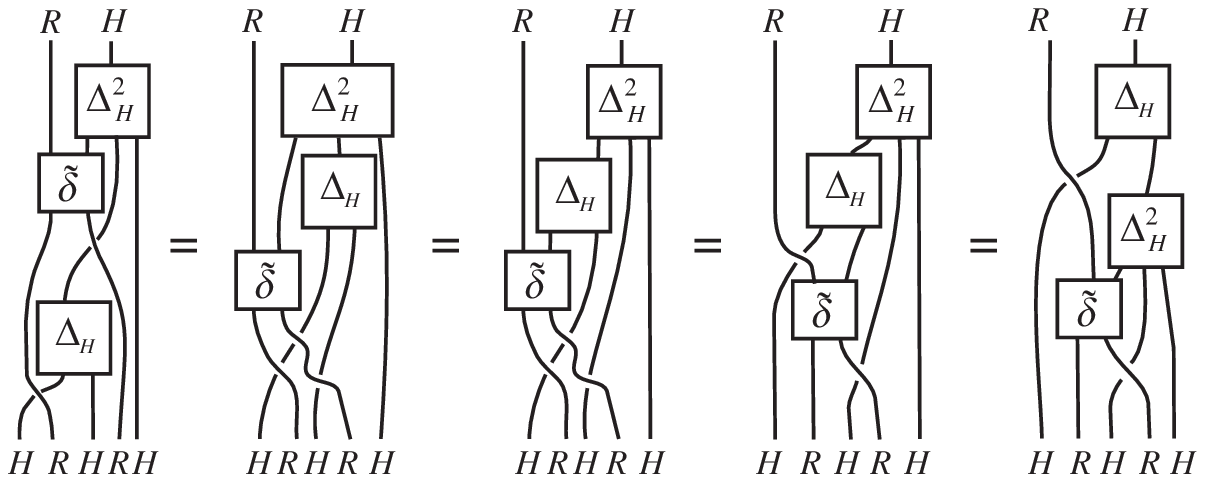}
\end{center}
\par
\vspace*{-0.5cm}
\caption{The proof of (\ref{eqdel1}) $\Longrightarrow $ (\ref{ec:DeltaColin}%
) } \label{fig:DeltaLColin3}
\end{figure}
\end{proof}

\begin{lemma}
\label{Lemma 0.4}Assume that $\Delta _{R\#H}$ is left
$H$\lin{}{}colinear (i.e. satisfies (\ref{ec:DeltaColin})). Then
the following two relations are equivalent:
\begin{eqnarray}
&\lbrack \Delta _{R\#H}\otimes (R\#H)]\Delta _{R\#H}
=[(R\#H)\otimes \Delta _{R\#H})\Delta _{R\#H}
\label{ec:DeltaCoasoc} \\ &(\tilde{\delta}\otimes R)(R\otimes
c_{R,H})(\tilde{\delta}\otimes R)(R\otimes \Delta _{H}) =(R\otimes
\tilde{\delta})(\tilde{\delta}\otimes H)(R\otimes \Delta _{H})
\label{eqdel}
\end{eqnarray}
\end{lemma}

\begin{proof}
The diagrammatic representation of the equivalence is given in
Figure~\ref {fig:DeltaCoass1}.
\begin{figure}[t]
\begin{center}
\includegraphics{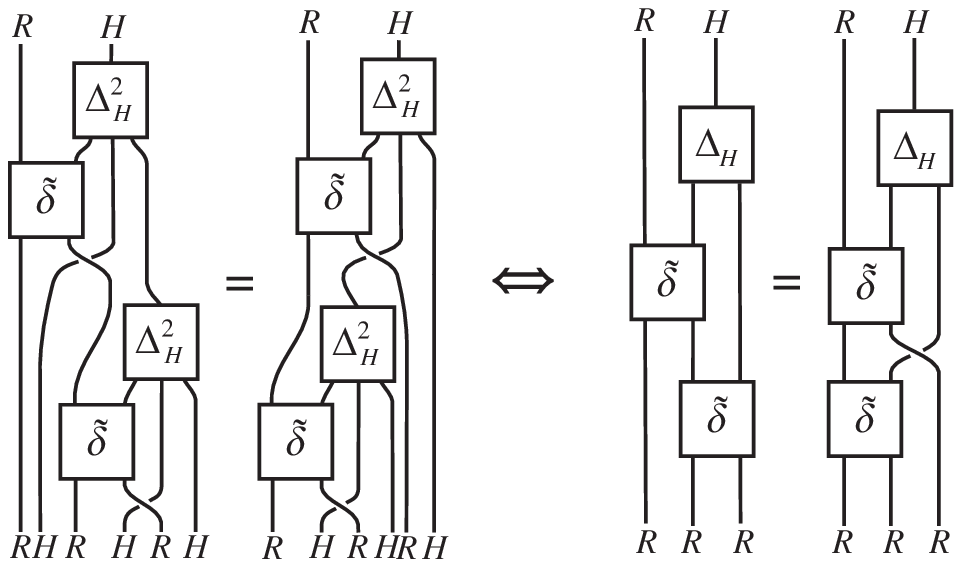}
\end{center}
\par
\vspace*{-0.5cm}
\caption{Representation of $($\ref{ec:DeltaCoasoc}$)\Longleftrightarrow ($%
\ref{eqdel}$)$} \label{fig:DeltaCoass1}
\end{figure}
\noindent It is easy to see that (\ref{ec:DeltaCoasoc}) implies (\ref{eqdel}%
). Indeed it is enough to add $(R\otimes \varepsilon _{H}\otimes
R\otimes \varepsilon _{H}\otimes R\otimes \varepsilon _{H})$ on
the bottom of the diagram representing (\ref{ec:DeltaCoasoc}),
then to drag $\varepsilon _{H}$ under the crossings and to use
that $\varepsilon _{H}$ is a counit. The other implication in
proved in Figure~\ref{fig:DeltaCoass2}.
\begin{figure}[b]
\begin{center}
\includegraphics{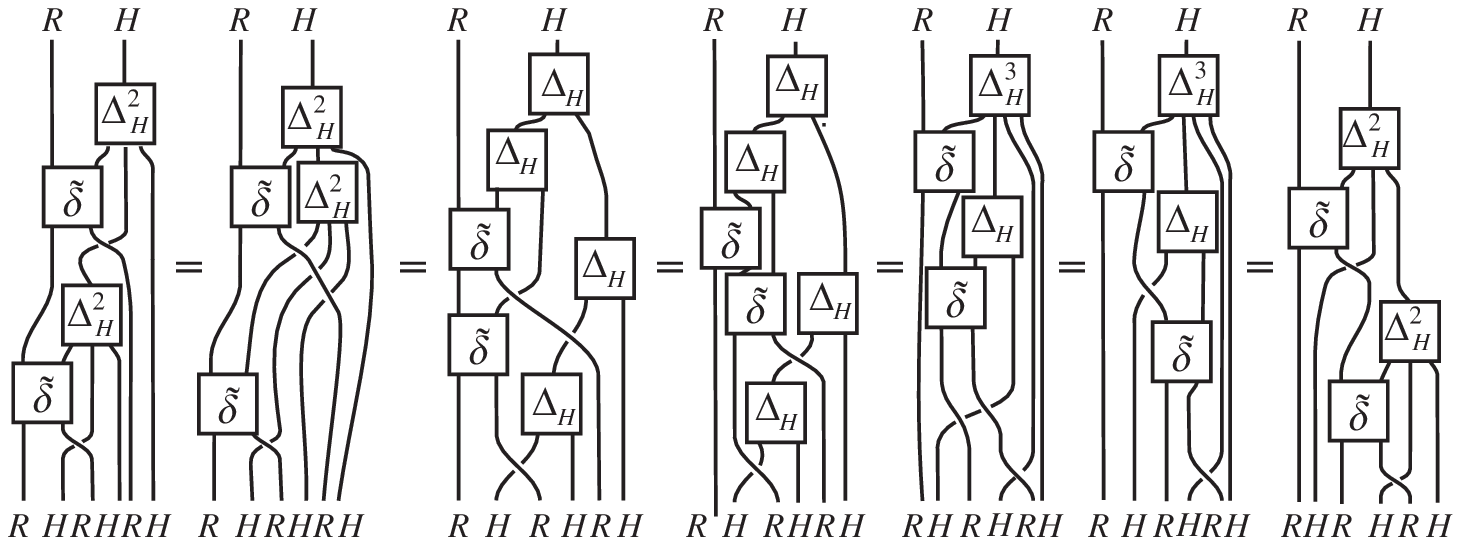}
\end{center}
\par
\vspace*{-0.5cm}
\caption{The proof of $($\ref{eqdel}$)\Longrightarrow ($\ref{ec:DeltaCoasoc}$%
)$} \label{fig:DeltaCoass2}
\end{figure}
\end{proof}

\begin{claim}
\label{cl:BraidTensor} Let $R$ and $S$ be two algebras in the
braided category $\yd$. We can define a new algebra structure on
$R\otimes S$, by using the braiding (\ref{eq:Braiding}), and not
the usual flip morphism. The multiplication in this case is
defined by the formula:
\begin{equation}  \label{eq:BraidTensor}
\left( r\otimes s\right) \left( t\otimes v\right) =\sum
r({}^{s{}_{\langle -1\rangle}}t)\otimes s{}_{\langle 0\rangle}v.
\end{equation}
Let us remark that for any algebra $R$ in ${}_H^H\mathcal{YD}$ the
smash product $R\#H$ is a particular case of this construction.
Just take $S=H$ with the left adjoint coaction and usual left
$H$\lin{}{}module structure. Another example that we are
interested in is $R\otimes R$, where $R$ is the diagram
of a bialgebra $A$ as in (\ref{cl:Settings}). For such an algebra $R$ in $%
{}_H^H\mathcal{YD}$ we shall always use this algebra structure on
$R\otimes R $.
\end{claim}

\begin{lemma}
\label{prodel}Let $\tilde{\delta}:R\otimes H\rightarrow R\otimes R$ be a $K$%
\lin{}{}linear map. Then the following two relations are
equivalent:
\begin{eqnarray}
&\Delta _{R\#H}\left( (r\# h)(s\# k)\right) =\Delta _{R\#H}(r\#
h)\Delta _{R\#H}(s\# k),  \label{eq:DeltaMult} \\
&\tilde{\delta}\left( (r\# h)(s\# k)\right) =\sum
\tilde{\delta}(r\# h_{(1)})\,{}^{h_{(2)}}\tilde{\delta}(s\# k).
\label{eqpdel}
\end{eqnarray}
where, for every $h\in H$ and $r,t\in R$ we have $^h(r\otimes
t)=\sum {}^{h_{\left( 1\right)} }r\otimes {}^{h_{\left( 2\right)}
}t. $
\end{lemma}

\begin{proof}
Let $r\# h$ and $s\# k\in R\#H.$ Thus we have:
\begin{align*}
\Delta (r\# h)=\sum& \tilde{\delta}^{1}(r\# h_{(1)})\# \tilde{\delta}%
^{2}(r\# h_{(1)}){}_{\langle-1\rangle}h_{(2)}\otimes
\tilde{\delta}^{2}(r\# h_{(1)}){}_{\langle 0\rangle}\# h_{(3)} \\
\Delta (s\# k)=\sum &\tilde{\delta}^{1}(s\# k_{(1)})\# \tilde{\delta}%
^{2}(s\# k_{(1)}){}_{\langle-1\rangle}k_{(2)}\otimes
\tilde{\delta}^{2}(s\# k_{(1)}){}_{\langle 0\rangle}\# k_{(3)}, \\
\Delta _{R\#H}\left( (r\# h)(s\# k)\right)=\sum &\tilde{\delta}%
^{1}(r{}^{h{}_{\left( 1\right)}}s\# h{}_{\left( 2\right)}k_{(1)})\# \tilde{%
\delta}^{2}(r{}^{h{}_{\left( 1\right)}}s\# h{}_{\left(
2\right)}k_{(1)}){}_{\langle-1\rangle}h{}_{\left(
3\right)}k_{(2)}\otimes \\ &\otimes
\tilde{\delta}^{2}(r{}^{h{}_{\left( 1\right)}}s\# h{}_{\left(
2\right)}k_{(1)}){}_{\langle 0\rangle}\# h{}_{\left(
4\right)}k_{(3)}.
\end{align*}
By substituting in (\ref{eq:DeltaMult}) the elements involving
$\Delta_{R\#H} $ with the right hand sides of the above three
relations, and then by applying $R\otimes \varepsilon _{H}\otimes
R\otimes \varepsilon _{H}$ it results:
\begin{equation}  \label{eq:xx}
\tilde{\delta}((r\# h)(s\# k))=\sum \tilde{\delta}^{1}(r\# h_{(1)}) \;{}^{%
\tilde{\delta}^{2}(r\# h_{(1)}){}_{\langle-1\rangle}h_{(2)}}\tilde{\delta}%
^{1}(s\# k) \otimes \tilde{\delta}^{2}(r\# h_{(1)}){}_{\langle
0\rangle} \;{}^{h_{(3)}}\tilde{\delta}^{2}(s\# k)
\end{equation}
Since in $R\otimes R$ the multiplication is defined by (\ref{eq:BraidTensor}%
) it follows that the right hand sides of (\ref{eqpdel}) and
(\ref{eq:xx}) are equal, so the equality (\ref{eqpdel}) holds.

Conversely, if (\ref{eqpdel}) holds true then we have
(\ref{eq:xx}). We can
replace the left hand side of this relation by $\sum \tilde{\delta}%
^{1}(r\;{}^{ h_{\left( 1\right) }}s \# h_{(2)}k)\otimes \tilde{\delta}%
^{2}(r\;{}^{ h_{\left( 1\right) }}s \# h_{(2)}k)$. A very long
computation, using this equivalent form of (\ref{eq:xx}), ends the
proof of the proposition.
\end{proof}

\begin{claim}
Let $\tilde{\delta}:R\otimes H\rightarrow R\otimes R$ be a
$K$\lin{}{}linear map. For every $r\in R$ and $h\in H$ we
introduce the notation:
\begin{equation}  \label{eq:DeltOmega}
\delta (r) =\tilde{\delta}(r\#1) \qquad\omega
(h)=\tilde{\delta}(1\#h).
\end{equation}
Then $\delta :R\rightarrow R\otimes R$ and $\omega :H\rightarrow
R\otimes R$
are $K$\lin{}{}linear map. Recall that $R\otimes R$ is an algebra in ${}_H^H%
\mathcal{YD}$ with the multiplication defined in
(\ref{cl:BraidTensor}). For example we can compute the product
$\delta(r)\omega(h)$ in $R\otimes R$. Now, using the notation
above, we can give a new interpretation of (\ref {eq:DeltaMult}).
\end{claim}

\begin{lemma}
\label{delta multi}Let $\tilde{\delta}:R\otimes H\rightarrow
R\otimes R$ be a $K$\lin{}{}linear map. Then $\Delta _{R\#H}$ is a
morphism of algebras iff $\delta(1_R)=1_R\otimes 1_R$,
$\omega(1_H)=1_R\otimes 1_R$ and $\delta , \tilde{\delta}$ and
$\omega$ satisfy the following four relations:
\begin{eqnarray}
&\tilde{\delta}(r\#h)=\delta (r)\omega (h)  \label{del3} \\
&\delta (rs)=\delta (r)\delta (s)  \label{del1} \\ &\omega
(hk)=\sum \omega (h_{(1)})\;{}^{h_{(2)}} \omega (k),  \label{del4}
\\ &\sum \delta \left( {}^{h_{\left( 1\right)} }r\right) \omega
(h_{(2)})=\sum \omega (h_{\left( 1\right) })\;{}^{h_{\left(
2\right) }}\delta (r) \label{del2}
\end{eqnarray}
\end{lemma}

\begin{proof}
By Lemma \ref{prodel}, the map $\Delta_{R\#H}$ is multiplicative
if and only if (\ref{eqpdel}) hold i.e.:
\begin{equation*}
\tilde{\delta}\left( (r\# h)(s\# k)\right) =\sum
\tilde{\delta}(r\# h_{(1)})\,{}^{h_{(2)}}\tilde{\delta}(s\# k).
\end{equation*}
Now assume that (\ref{eqpdel}) holds. Then setting $h=1_{H}=k$ we obtain (%
\ref{del1}), while for $r=1_{R}=s$ we obtain (\ref{del4}). Also
for $h=1_{H}$
and $s=1_{R}$ we get (\ref{del3}) and for $r=1_{R}$ and $k=1_{H}$ we get (%
\ref{del2}), by means of (\ref{del3}). Conversely assume that (\ref{del1}), (%
\ref{del2}), (\ref{del3}) and (\ref{del4}) hold true. Then by (\ref{del3}), (%
\ref{del1}) and (\ref{del4}) we have:
\begin{equation*}
\tilde{\delta}((r\# h)(s\# k)) =\sum \delta (r\,{}^{h_{(1)}}s)
\omega(h_{(2)}k) =\sum \delta (r) \delta ({}^{h_{(1)}}s)
\omega(h_{(2)})\,{}^{h_{(3)}} \omega (k).
\end{equation*}
So, by (\ref{del2}) and by the fact that $R\otimes R$ is an algebra in $%
{}_H^H\mathcal{YD}$ (hence an $H$\lin{}{}module algebra), we get:
\begin{equation*}
\tilde{\delta}((r\# h)(s\# k))=\sum \delta (r)\omega(h_{(1)})
\,{}^{h_{(2)}} \delta (s) \,{}^{h_{(3)}}\omega (k)=\sum
\delta(r)\omega(h_{(1)})\,{}^{h_{(2)}}[\delta(s)\omega(k)].
\end{equation*}
Now we can prove (\ref{eqpdel}) by using (\ref{del3}) once again. Obviously $%
\Delta_{R\#H}$ is a morphism of unitary rings if and only if $%
\delta(1_R)=1_R\otimes 1_R$ and $\omega(1_H)=1_R\otimes 1_R$.
\end{proof}

\begin{remark}
By (\ref{del3}) we can recover $\widetilde{\delta}$ from $\delta
:R\rightarrow R\otimes R$ and $\omega :H\rightarrow R\otimes R$. Equation (%
\ref{del1}) says that $\delta $ is multiplicative with the algebra
structure on $R\otimes R$ introduced in (\ref{cl:BraidTensor}). We
have already noticed that $R\otimes R$ is a left $H$\lin{}{}module
algebra. Now if $A$ is an arbitrary left $H$\lin{}{}module
algebra, then Sweedler in \cite{Sw} defined a noncommutative
$1$\lin{}{}cocycle with coefficient in $A$ to be a
$K$\lin{}{}linear map $\theta :H\rightarrow A$ such that
\begin{equation*}
\theta (hk)=\sum \theta (h_{(1)}){}^{h_{(2)}}\theta (k).
\end{equation*}
Hence (\ref{del4}) means that $\omega $ is a $1$\lin{}{}cocycle
with coefficients in $A.$
\end{remark}

\begin{lemma}
\label{le:colin}Assume that $\Delta _{R\#H}$ is multiplicative.
Then (\ref {eqdel1}) holds iff $\delta $ and $\omega $ are left
$H$\lin{}{}colinear (where $H$ is a left $H$\lin{}{}comodule with
the left adjoint coaction).
\end{lemma}

\begin{proof}
Assume that (\ref{eqdel1}) holds and let $r\in R$, $h\in H$. By evaluating (%
\ref{eqdel1}) at $r\# 1$ we get:
\begin{equation*}
\rho _{R\otimes R}^{l}(\delta (r))= \sum
r{}_{\langle-1\rangle}\otimes \delta (r{}_{\langle 0\rangle}),
\end{equation*}
so $\delta$ is $H$\lin{}{}colinear. Similarly, for $1\# h$ we
have:
\begin{equation*}
\sum h_{(1)}\otimes \omega ^{1}(h_{(2)})\otimes \omega
^{2}(h_{(2)}) =\sum \omega
^{1}(h_{(1)}){}_{\langle-1\rangle}\omega ^{2}(h_{(1)}){}_{\langle
-1\rangle}h_{(2)}\otimes \omega ^{1}(h_{(1)}){}_{\langle
0\rangle}\otimes \omega ^{2}(h_{(1)}){}_{\langle 0\rangle}
\end{equation*}
i.e. we get
\begin{equation}
\sum h_{(1)}\otimes \omega (h_{(2)})=\sum \omega
(h_{(1)}){}_{\langle -1\rangle}h_{(2)}\otimes \omega
(h_{(1)}){}_{\langle 0\rangle}.  \label{eqcl}
\end{equation}
On the other hand:
\begin{equation*}
\sum \omega (h){}_{\langle-1\rangle}\otimes \omega (h){}_{\langle
0\rangle}=\sum \omega (h_{\left( 1\right) }){}_{\langle
-1\rangle}h_{(2)}S(h_{(3)})\otimes \omega (h_{\left( 1\right)
}){}_{\langle 0\rangle}=\sum h_{(1)}S(h_{(3)})\otimes \omega
(h_{(2)}),
\end{equation*}
where the last equality holds in view of (\ref{eqcl}). Hence
$\omega $ is left $H$\lin{}{}colinear.

Conversely assume that $\delta $ and $\omega $ are left
$H$\lin{}{}colinear. The
relation (\ref{eqdel1}) that we have to prove is equivalent to $%
A_l(r,s)=A_r(r,s)$, where:
\begin{eqnarray}
A_l(r,s)&=&\sum\widetilde{\delta}^1(r\#h{}_{\left(
1\right)}){}_{\langle
-1\rangle}\widetilde{\delta}^2(r\#h{}_{\left(
1\right)}){}_{\langle -1\rangle}h{}_{\left(
2\right)}\otimes\widetilde{\delta}^1(r\#h{}_{\left(
1\right)}){}_{\langle 0\rangle}\widetilde{\delta}^2(r\#h{}_{\left(
1\right)}){}_{\langle 0\rangle} \\
A_r(r,s)&=&\sum r{}_{\langle-1\rangle}h{}_{\left( 1\right)}\otimes%
\widetilde{\delta}(r{}_{\left( 0\right)}\#h{}_{\left( 2\right)})
\end{eqnarray}
Then, since $\Delta _{R\#H}$ is multiplicative, by (\ref{del3}) we
have:
\begin{equation*}
A_l(r,s)=\sum \widetilde{\delta}(r\otimes h_{(1)}){}_{\langle
-1\rangle}h_{\left( 2\right) }\otimes \widetilde{\delta}(r\otimes
h_{(1)}){}_{\langle 0\rangle} =\sum \delta
(r){}_{\langle-1\rangle}\omega (h_{(1)}){}_{\langle
-1\rangle}h_{(2)}\otimes \delta (r){}_{\langle 0\rangle}\omega
(h_{(1)}){}_{\langle 0\rangle}.
\end{equation*}
Since $\delta$ and $\omega$ are left colinear it results:
\begin{equation*}
\sum \delta (r){}_{\langle-1\rangle}\omega (h_{(1)}){}_{\langle
-1\rangle}h_{(2)}\otimes \delta (r){}_{\langle 0\rangle}\omega
(h_{(1)}){}_{\langle 0\rangle} =\sum r{}_{\langle
-1\rangle}h_{(1)}\otimes \delta (r{}_{\langle 0\rangle})\omega
(h_{(2)}) =A_r(r,s),
\end{equation*}
so $A_l(r,s)=A_r(r,s)$, thus the lemma has been proved.
\end{proof}

\begin{claim}
To simplify the notation, for every $r\in R$, let $\delta
(r):=\sum r^{(1)}\otimes r^{(2)}$. This is a kind of
$\Sigma$\lin{}{}notation that we shall use for $\delta$.
\end{claim}

\begin{lemma}
\label{le:coass}Assume that $\Delta _{R\#H}$ is a morphism of
algebras such that $\delta $ is left $H$\lin{}{}colinear. Then
(\ref{eqdel}) holds iff the following two relations hold true for
any $r\in R$ and $h\in H$:
\begin{eqnarray}
&\sum r^{(1)}\otimes \delta (r^{(2)})=\sum \delta (r^{(1)})\omega
(r{}_{\langle-1\rangle}^{(2)})\otimes r{}_{\langle 0\rangle}^{(2)}
\label{eqcli} \\ &\sum \omega ^{1}(h_{(1)})\otimes \delta
\left(\omega ^{2}(h_{(1)})\right)\omega (h_{(2)})=\sum \delta
\left(\omega ^{1}(h_{(1)})\right)\omega \left( \omega
^{2}(h_{(1)}){}_{\langle -1\rangle}h_{(2)}\right) \otimes \omega
^{2}(h_{(1)}){}_{\langle 0\rangle} \label{eqclio}
\end{eqnarray}
\end{lemma}

\begin{proof}
Since $\Delta_{R\#H}$ is multiplicative it is straightforward to prove that (%
\ref{eqdel}) holds iff, for every $r\in R$ and $h\in H$, we have $%
B_l(r,h)=B_r(r,h)$, where:
\begin{eqnarray}
&B_l(r,h)=\sum r^{(1)}\left({}^{r^{(2)}_{{}_{\langle
-1\rangle}}}\omega^{1}(h{}_{\left(
1\right)})\right)\otimes\delta\left(r^{(2)}_{{}_{\langle
0\rangle}}\omega^{2}(h{}_{\left(
1\right)}\right)\omega(h{}_{\left( 2\right)}) .
\label{coasocnou1} \\
&B_r(r,h)=\sum\delta\left(r^{(1)}\,{}^{r^{(2)}_{{}_{\langle
-2\rangle}}}\omega^{1}(h{}_{\left(
1\right)})\right)\omega\left(r^{(2)}_{{}_{\langle
-1\rangle}}\omega^{2}(h{}_{\left( 1\right)})\un{-1}h{}_{\left(
2\right)}\right) \otimes r^{(2)}_{{}_{\langle
0\rangle}}\omega^{2}(h{}_{\left( 1\right)}){}_{\langle 0\rangle} .
\label{coasocnou}
\end{eqnarray}
Since $\Delta_{R\#H}$ is a morphism of algebras we have $\delta(1_R)=1_R%
\otimes 1_R$ and $\omega(1_H)=1_R\otimes 1_R$. Hence one can see
easily that
(\ref{eqcli}) and (\ref{eqclio}) are equivalent to $B_l(r,1)=B_r(r,1)$ and $%
B_l(1,h)=B_r(1,h)$, respectively. In particular, (\ref{eqdel})
implies (\ref {eqcli}) and (\ref{eqclio}). In order to probe the
converse, let us denote by $C_l(h)$ and $C_r(h)$ the left and
right hand sides of (\ref{eqclio}). Since $\delta$ is left
$H$\lin{}{}colinear, and by using (\ref{del1}), it results:
\begin{equation*}
B_l(r,h)=\sum \left(r^{(1)}\otimes\delta(r^{(2)})\right)C_l(h),
\end{equation*}
where the product is performed in $R\ot R\ot R$, which is an
algebra with the multiplication:
\[
(r\ot s\ot t)(r'\ot s'\ot t')=\sum r\;{}^{s\un{-1}t\un{-2}}r'\ot
s\un{0}{}^{t\un{-1}}s'\ot t\un{0}t'.
\]
Similarly, by (\ref{del3}) and (\ref{eqpdel}), it follows:
\begin{equation*}
B_r(r,h)=\sum \left(\delta(r^{(1)})\omega(r^{(2)}_{{}_{\langle
-1\rangle}})\otimes r^{(2)}_{{}_{\langle 0\rangle}}\right)C_r(h).
\end{equation*}
We deduce that $B_l(r,h)=B_r(r,h)$ by multiplying (\ref{eqcli}%
) and (\ref{eqclio}) side by side in $R\ot R\ot R$.
\end{proof}

\begin{claim}
\label{cl:EpsBar} \label{claim1}Let $W\in {^{H}}\mathfrak{M}%
{^{H}}$ and let $\varepsilon :W\rightarrow K$ be a
$K$\lin{}{}linear map. Let us define $\overline{\varepsilon
}:W\rightarrow H$ by:
\begin{equation}
\overline{\varepsilon }(w)=\sum \varepsilon (w{}_{\langle
0\rangle})w{}_{\langle 1\rangle}.  \label{def epsbar}
\end{equation}
$\overline{\varepsilon }$ is always right $H$\lin{}{}colinear and $\overline{%
\varepsilon }$ is left $H$\lin{}{}colinear (and hence a morphism
of bicomodules) iff:
\begin{equation*}
\overline{\varepsilon }(w)=\sum w{}_{\langle-1\rangle}\varepsilon
(w{}_{\langle 0\rangle}).
\end{equation*}
Let $V:=W^{coH}.$ Then $V$ is a left $H$\lin{}{}subcomodule of $W.$ Let $%
\varepsilon _{V}$ be the restriction of $\varepsilon $ to $V.$
Then:
\begin{equation*}
\overline{\varepsilon }(v)=\sum \varepsilon (v{}_{\langle
0\rangle})v{}_{\langle 1\rangle}=\varepsilon _{V}(v)1_{H},\forall v\in V%
\text{. }
\end{equation*}
Thus, if $\overline{\varepsilon }$ is left $H$\lin{}{}colinear,
then $\varepsilon _{V}:V\rightarrow K$ is a morphism of left
$H$\lin{}{}comodules.

Now let us consider the following particular case. Take $M$ to be a left $H$%
\lin{}{}comodule and let $W:=M\otimes H$. Then $W$ is an
$(H,H)$\lin{}{}bimodule with respect to the diagonal left coaction
and canonical right comodule
structure. Let $\varepsilon_M :=\varepsilon _{M}:M\rightarrow K$ be a $K$%
\lin{}{}linear map and define $\varepsilon_{}:W\rightarrow K$ by $%
\varepsilon_{}(m\otimes h)=\varepsilon_{M}(m)\varepsilon_{H}(h)$.
Since the
right comodule structure on $W$ is induced by $\Delta_H$ we have $%
W{}^{Co(H)}=\{m\otimes 1_H\mid m\in M\}.$

By the foregoing, if $\overline{\varepsilon}_{M\otimes H}$ is a
morphism of $(H,H)$\lin{}{}bicomodules then:
\begin{equation*}
\sum \varepsilon_M (m{}_{\langle 0\rangle})m{}_{\langle
-1\rangle}=\varepsilon_M (m)1_{H},\forall m\in M,
\end{equation*}
that is $\varepsilon_{M}$ is left $H$\lin{}{}colinear. Actually,
in this particular case, we can prove that the converse also holds
true. In
conclusion, $\overline{\varepsilon}_{M\otimes H}$ is a morphism of $(H,H)$%
\lin{}{}bicomodules iff $\varepsilon_{M}$ is left
$H$\lin{}{}colinear.
\end{claim}

\begin{lemma}
\label{Lemma 0.7}Let $R$ be an algebra in ${_{H}^{H}\mathcal{YD}}$ and let $%
\varepsilon:R\rightarrow K$ be a $K$\lin{}{}linear map. The map
$\varepsilon _{R\#H}:R\#H\rightarrow K,\varepsilon
_{R\#H}(r\otimes
h):=\varepsilon(r)\varepsilon _{H}(h),$ is an algebra map and $\overline{%
\varepsilon}_{R\#H}:R\#H\rightarrow H$ is a left
$H$\lin{}{}colinear map if and only if $\varepsilon$ is an algebra
map in ${_{H}^{H}\mathcal{YD}}$.
\end{lemma}

\begin{proof}
``$\Rightarrow$'' By (\ref{cl:EpsBar}) it follows that $\varepsilon:R%
\rightarrow K$ is left $H$\lin{}{}colinear. By the definition of
the multiplication in $R\#H$ and the definition of
$\varepsilon_{R\#H}$ we get:
\begin{equation*}
\sum \varepsilon(r\,{}^{h_{(1)}}s)\varepsilon _{H}\left(
h_{(2)}v\right) =\varepsilon_{R\#H} \left( (r\#h)(s\#v)\right)
=\varepsilon_{R\#H} (r\#h)\varepsilon_{R\#H}
(s\#v)=\varepsilon(r)\varepsilon _{H}\left( h\right)
\varepsilon(s)\varepsilon _{H}(v).
\end{equation*}
Thus $\varepsilon({}^{h}s)=\varepsilon _{H}(h) \varepsilon(s)$ and
$\varepsilon (rs)= \varepsilon(r)\varepsilon(s)$, i.e.
$\varepsilon $ is an $H$\lin{}{}colinear algebra map.

``$\Leftarrow$'' We know already by (\ref{cl:EpsBar}) that $\overline{%
\varepsilon}_{R\#H}:R\#H\rightarrow H$ is a left
$H$\lin{}{}colinear map. Furthermore:
\begin{equation*}
\overline{\varepsilon}_{R\#H} \left(^{} (r\#h)(s\#k)\right) =\sum \overline{%
\varepsilon}_{R\#H} \left( r\,{}^{h_{(1)}}s s \#h_{(2)}k\right)
=\sum \varepsilon (r)\varepsilon
_{H}(h_{(1)})\varepsilon(s)\varepsilon _{H}(h_{(2)})\varepsilon
_{H}(k)
\end{equation*}
Hence $\overline{\varepsilon}_{R\#H} \left( (r\#h)(s\#k)\right) =\overline{%
\varepsilon}_{R\#H}(r\#h)\varepsilon (s\#k)$, so
$\varepsilon_{R\#H}$ is an algebra map.
\end{proof}

\begin{lemma}
\label{Lemma 0.8}Assume that $\varepsilon$ is an algebra map in ${_{H}^{H}%
\mathcal{YD}}$. Then $\varepsilon _{R\#H}:R\#H\rightarrow K$ is a
counit for $\Delta_{R\#H} $ if and only if, for every $r\in R$ and
$h\in H$, we have:
\begin{equation}
\sum \varepsilon(\tilde{\delta}^{1}(r\otimes
h))\tilde{\delta}^{2}(r\otimes h)=\varepsilon _{H}(h)r=\sum
\tilde{\delta}^{1}(r\otimes h)\varepsilon\left(
\tilde{\delta}^{2}(r\otimes h)\right).  \label{eqdel4}
\end{equation}
\end{lemma}

\begin{proof}
Assume that $\varepsilon_{R\#H}$ is a counit for $\Delta_{R\#H}$.
Then, by the definition of $\Delta_{R\#H}$, see (\ref{eqdd}), it
results:
\begin{equation*}
r\otimes h=\sum\tilde{\delta}^{1}(r\otimes h_{(1)})\otimes \tilde{\delta}%
^{2}(r\otimes h_{(1)}){}_{\langle-1\rangle}h_{(2)}\varepsilon\left( \tilde{%
\delta}^{2}(r\otimes h_{(1)}){}_{\langle 0\rangle}\right)
\varepsilon _{H}\left( h_{(3)}\right).
\end{equation*}
By applying $R\otimes\varepsilon_{H}$ to this relation we get the
second equality of (\ref{eqdel4}). The other one can be proved
similarly.

Conversely assume that the equality (\ref{eqdel4}) holds. Since
$\varepsilon $ is left $H$\lin{}{}colinear, we have
\begin{equation*}
(R\#H\otimes\varepsilon_{R\#H})\Delta_{R\#H}=\sum\tilde{\delta}^{1}(r\otimes
h_{(1)}) \varepsilon\left( \tilde{\delta}^{2}(r\otimes
h_{(1)})\right) \otimes h_{(2)}=\sum r\varepsilon _{H}\left(
h_{(1)}\right) \otimes h_{(2)}=r\otimes h.
\end{equation*}
We can prove the second relation analogously.
\end{proof}

\begin{lemma}
\label{le:counit}Assume that $\Delta _{R\#H}$ is multiplicative and that $%
\varepsilon:R\rightarrow K$ is an algebra map in
${_{H}^{H}\mathcal{YD}} $. Then (\ref{eqdel4}) holds if and only
if:
\begin{eqnarray}
&(\varepsilon\otimes R)\delta =(R\otimes \varepsilon )\delta
=\mathrm{Id}_{R} \label{eqe1} \\ &(\varepsilon \otimes R)\omega
=(R\otimes \varepsilon )\omega =\varepsilon _{H}1_{R}
\label{eqe2}
\end{eqnarray}
\end{lemma}

\begin{proof}
First let us observe that $\varepsilon_{}\otimes R:R\otimes
R\rightarrow R$ and $R\otimes \varepsilon_{}:R\otimes R\rightarrow
R$ are algebra morphisms (recall that $R\otimes R$ is an algebra
with the multiplication $(m_R\otimes
m_R\otimes R)(R\otimes c_{R,R})$, where $c$ is the brainding in ${}_H^H%
\mathcal{YD}~$). Clearly (\ref{eqdel4}) holds if and only if:
\begin{equation*}
(\varepsilon \otimes R)\widetilde{\delta }\left( r\#h\right)
=\varepsilon _{H}\left( h\right) r=(R\otimes \varepsilon
)\widetilde{\delta }\left( r\#h\right) ,\forall r\in R,\ \forall
h\in H.
\end{equation*}
Assume now that (\ref{eqe1}) and (\ref{eqe2}) holds. Then:
\begin{equation*}
(\varepsilon \otimes R)\widetilde{\delta }\left( r\#h\right)
=(\varepsilon \otimes R)\delta \left( r\right) \cdot (\varepsilon
\otimes R)\omega \left( h\right) =\varepsilon _{H}\left( h\right)
r.
\end{equation*}
Analogously we can deduce the second equality of (\ref{eqdel4}).
The other implication is trivial.
\end{proof}

To state easier the main results of this part we collect together
in the
next definition all required properties of $\delta$, $\omega$ and $%
\varepsilon_{}$.

\begin{definition}\label{de:YDq}
Let $H$ be a Hopf algebra and let $R$ be an algebra in
${}_H^H\mathcal{YD}$. Assume that $\varepsilon :R\rightarrow K$ ,
$\delta :R\rightarrow R\otimes R$ and $\omega
:H\rightarrow R\otimes R$ are $K$\lin{}{}linear maps. The quadruple $%
(R,\varepsilon,\delta ,\omega )$ will be called a
\emph{Yetter--Drinfeld quadruple}\ if and only if, for all $r,
s\in R$ and $h,k\in H$, the following relations are satisfied:
\begin{eqnarray}
& \varepsilon_{}({}^hr)=\varepsilon_{}(r)\varepsilon_{H}(h) \qquad \text{and}%
\qquad \sum r{}_{\langle-1\rangle}\varepsilon_{}(r{}_{\langle
0\rangle})=\varepsilon_{}(r)1_H;  \label{eq:YD0} \\
& \varepsilon_{}(rs)=\varepsilon_{}(r)\varepsilon_{}(s) \qquad \text{and}%
\qquad \varepsilon_{}(1_R)=1;  \label{eq:YD1} \\ &\rho_{R\otimes
R}\left(\delta(r)\right)=\sum r{}_{\langle
-1\rangle}\otimes\delta(r{}_{\langle 0\rangle});  \label{eq:YD2}
\\ &\rho_{R\otimes R}\left(\omega(h)\right)=\sum h{}_{\left(
1\right)}S(h{}_{\left( 3\right)})\otimes\omega(h{}_{\left(
2\right)}); \label{eq:YD3} \\ &\delta (rs)=\delta (r)\delta (s)
\qquad \text{and}\qquad \delta \left( 1_{R}\right)= 1_{R}\otimes
1_{R};  \label{eq:YD4} \\ & \omega (hk)=\sum \omega
(h_{(1)})\left( {}^{h_{(2)}}\omega (k)\right)\qquad
\text{and}\qquad \omega\,( 1_{H}); =1_{R}\otimes 1_{R};
\label{eq:YD5} \\ & \sum \delta \left({}^{h_{(1)}}r\right) \omega
(h_{(2)})=\sum \omega (h_{\left( 1\right) })\, {}^{h_{\left(
2\right) }}\delta (r);  \label{eq:YD6}
\\
& \sum r^{(1)}\otimes \delta (r^{(2)})=\sum \delta (r^{(1)})\omega
(r{}_{\langle-1\rangle}^{(2)})\otimes r{}_{\langle
0\rangle}^{(2)}; \label{eq:YD7} \\ &\sum \omega
^{1}(h_{(1)})\otimes \delta (\omega ^{2}(h_{(1)}))\omega
(h_{(2)})=\sum \delta (\omega ^{1}(h_{(1)}))\omega \left( \omega
^{2}(h_{(1)}){}_{\langle-1\rangle}h_{(2)}\right) \otimes \omega
^{2}(h_{(1)}){}_{\langle 0\rangle};  \label{eq:YD8} \\
& (\varepsilon \otimes R)\delta =(R\otimes \varepsilon )\delta =\mathrm{Id}%
_{R};  \label{eq:YD9} \\ & (\varepsilon \otimes R)\omega
=(R\otimes \varepsilon )\omega =\varepsilon _{H}1_{R}.
\label{eq:YD10}
\end{eqnarray}
\end{definition}

\begin{remark}
Note that these relations can be interpreted as follows:

(\ref{eq:YD0}) \lin{}{} $\varepsilon_{}$ is a morphism in
${}_H^H\mathcal{YD}$;

(\ref{eq:YD1}) \lin{}{} $\varepsilon_{}$ is a morphism of
algebras;

(\ref{eq:YD2}) \lin{}{} $\delta$ is left $H$\lin{}{}colinear;

(\ref{eq:YD3}) \lin{}{} $\omega$ is left $H$\lin{}{}colinear,
where $H$ is a comodule with the adjoint coaction;

(\ref{eq:YD4}) \lin{}{} $\delta{}$ is a morphism of algebras
 where on $R\otimes R$ we consider the
algebra structure that uses the braiding $c$;

(\ref{eq:YD5}) \lin{}{} $\omega$ is a normalized cocycle;

(\ref{eq:YD6}) \lin{}{} $\omega$ measures how far $\delta$ is to
be a morphism of left $H$\lin{}{}modules (if $\omega$ is trivial,
i.e. for every $h\in H$ we have $ \omega
(h)=\varepsilon_{}(h)1_R\otimes 1_R$, then $\delta$ is left
$H$\lin{}{}linear); we shall say that $\delta$ is a \emph{twisted
morphism} of left $H$\lin{}{}modules;

(\ref{eq:YD7}) \lin{}{} it was derived from the fact that
$\Delta_{R\#H}$ is coassociative, so we shall say that $\delta$ is
$\omega$\lin{}{}\emph{coassociative} (when $\omega$ is trivial
then (\ref{eq:YD7}) is equivalent to the fact that $\delta$ is
coassociative);

(\ref{eq:YD8}) \lin{}{} is the only property that has not an
equivalent in the
theory of bialgebras; we shall just say that $\delta$ and $\omega$ are \emph{%
compatible};

(\ref{eq:YD9}) \lin{}{} $\delta$ is a \emph{counitary map} with respect to $%
\varepsilon_{}$;

(\ref{eq:YD10}) \lin{}{} $\omega$ is a \emph{counitary map} with respect to $%
\varepsilon_{}$;

Since $\varepsilon_{}$ satisfies the last two relation we shall
call it the \emph{counit} of the Yetter--Drinfeld quadruple\ $R$.
By analogy $\delta $ will be called the \emph{comultiplication} of
$R$. Finally, we shall say that $\omega$ is the \emph{cocycle} of
$R$.
\end{remark}

\begin{claim}
To every Yetter--Drinfeld quadruple\
$(R,\varepsilon_{},\delta,\omega)$ we associate the
$K$\lin{}{}linear maps: $\Delta_{R\#H}:R\#H\rightarrow
(R\#H)\otimes (R\#H)$ and $\varepsilon_{R\#H}:R\#H\rightarrow K$,
which are defined by:
\begin{eqnarray}
&\Delta _{R\#H}(r\otimes h)=\sum \tilde{\delta}^{1}(r\otimes
h_{(1)})\otimes \tilde{\delta}^{2}(r\otimes h_{\left( 1\right)
}){}_{\langle -1\rangle}h_{(2)}\otimes \tilde{\delta}^{2}(r\otimes
h_{(1)}){}_{\langle 0\rangle}\otimes h_{\left( 3\right) }
\label{eq:dsm} \\ & \varepsilon _{R\#H}(r\#h)=\varepsilon
(r)\varepsilon _{H}(h) \label{eq:esm}
\end{eqnarray}
where $\widetilde{\delta }\left( r\#h\right) :=\delta \left(
r\right) \omega
\left( h\right)$ and $\widetilde{\delta}(r\#h)= \sum\tilde{\delta}%
^{1}(r\otimes h)\otimes \tilde{\delta}^{2}(r\otimes h)$.
\end{claim}

\begin{theorem}
\label{pro Boss} Let $R$ be an algebra in ${_{H}^{H}\mathcal{YD}}$. If $%
\varepsilon :R\rightarrow K$, $\delta :R\rightarrow R\otimes R$
and $\omega :H\rightarrow R\otimes R$ are linear maps, then the
following assertions are equivalent:

(a) $(R,\varepsilon_{},\delta,\omega)$ is a Yetter--Drinfeld
quadruple.

(b) The smash product algebra $R\#H$ is a bialgebra with the
comultiplication
$\Delta _{R\#H}$ and the counit $\varepsilon_{R\#H}$ defined by (\ref{eq:dsm}%
) and (\ref{eq:esm}) such that $R\#H$ becomes an algebra in $({_{H}^{H}}%
\mathfrak{M}{_{H}^{H},\otimes }_{H},H)$ and a coalgebra in $({^{H}}%
\mathfrak{M}{^{H}},{\square }_{H},H).$
\end{theorem}

\begin{proof}
\textit{(a) $\Rightarrow$ (b)} \ By Lemma \ref{delta multi} it
results that $\Delta _{R\#H}$ is multiplicative, in view of
(\ref{eq:YD4}), (\ref {eq:YD5}), (\ref{eq:YD6}). Note that we have
(\ref{del3}) by the definition of $\Delta_{R\#H}$. Since
$\delta(1_R)=\omega(1_H)=1_{R}\otimes 1_{R}$ we have $\Delta
_{R\#H}\left( 1_{R}\#1_{H}\right) =\left( 1_{R}\#1_{H}\right)
\otimes \left( 1_{R}\#1_{H}\right)$. In conclusion,
$\Delta_{R\#H}$ is a morphism of unitary algebras.

Since $\Delta _{R\#H}$ is multiplicative we can apply Lemma
\ref{Lemma 03}
and Lemma \ref{le:colin} to deduce that $\Delta _{R\#H}$ is left $H$%
\lin{}{}colinear by using relations (\ref{eq:YD2}) and
(\ref{eq:YD3}), i.e. that $ \delta $ and $\omega $ are left
$H$\lin{}{}colinear. On the other hand, by
(\ref{gammac}), we get that $\Delta _{R\#H}$ is right colinear, so $%
\Delta_{R\#H}$ is a morphism of $(H,H)$\lin{}{}bicomodules. Also
by (\ref{gammac})
it follows that the image of $\Delta_{R\#H}$ is included into $%
(R\#H)\square_H(R\#H)$.

Since $\Delta _{R\#H}$ is multiplicative and left
$H$\lin{}{}colinear and since $\delta$ is also left
$H$\lin{}{}colinear, by (\ref{eq:YD7}) and (\ref {eq:YD8}), it
results that $\Delta _{R\#H}$ is coassociative (use Lemma
\ref{le:coass} and Lemma \ref{Lemma 0.4}).

To prove that $\varepsilon_{R\#H}$ is a morphism of algebras we
use Lemma \ref{Lemma 0.7}, (\ref{eq:YD0}) and (\ref{eq:YD1}).
Finally, in view of (\ref{eq:YD9}) and (\ref{eq:YD10}), Lemma
\ref{le:counit} implies that $\varepsilon _{R\#H}$ is a counit for
$\Delta _{R\#H}.$ All these properties together
mean that $R\#H$ is a bialgebra that becomes a coalgebra in $({}^H %
\mathfrak{M}^H,\square_H,H)$. We conclude by remarking that $R\#H$
always becomes an algebra in
$({}_{H}^{H}\mathfrak{M}_{H}^{H},\otimes_H,H)$, see Example
\ref{Gold claim}.

In conclusion the object $R\#H$ is an algebra in $({_{H}^{H}} %
\mathfrak{M}{_{H}^{H},\otimes }_{H},H)$ . As $\varepsilon $ is an
algebra
map in ${_{H}^{H}\mathcal{YD}},$ by Lemma \ref{Lemma 0.7}, $\overline{%
\varepsilon }:R\#H\rightarrow H$ is a left $H$\lin{}{}colinear map
so that it is a map in $\hm^H$ (it is always right
$H$\lin{}{}colinear). Finally it is\ easy to check that the image
of $\Delta _{R\#H}$ is included in $R\#H{\square }_{H}R\#H.$

\textit{(b) $\Rightarrow$ (a)} \ Since $\Delta _{R\#H}$ is
morphism of algebras, by Lemma \ref{delta multi}, it follows that
(\ref{eq:YD4}), (\ref
{eq:YD5}) and (\ref{eq:YD6}) hold true. As $\Delta _{R\#H}$ is left $H$%
\lin{}{}colinear and multiplicative, by Lemma \ref{Lemma 03} and
Lemma \ref {le:colin}, $\delta $ and $\omega $ are
$H$\lin{}{}colinear, so that (\ref{eq:YD2}) and (\ref{eq:YD3})
hold.

Since we have already proved that $\delta $ is left
$H$\lin{}{}colinear, we can
apply Lemma \ref{le:coass} and Lemma \ref{Lemma 0.4} to deduce (\ref{eq:YD7}%
) and (\ref{eq:YD8}) from the fact that $\Delta_{R\#H}$ is
coassociative.

Since $R\#H$ becomes a coalgebra in $({}^H
\mathfrak{M}^H,\square_H,H)$ it results that the canonical map
$\overline{\varepsilon}_{R\#H}$ that is associated to
$\varepsilon_{R\#H}$ is left $H$\lin{}{}colinear (see \ref
{cl:BecomesCoalg}). This property and the fact that $\varepsilon
_{R\#H}$ is an algebra map imply (\ref{eq:YD0}) and
(\ref{eq:YD1}), in view of Lemma \ref{Lemma 0.7}.

Since $\varepsilon _{R\#H}$ is a counit for $\Delta _{R\#H}$, and since $%
\Delta _{R\#H}$ is multiplicative, by Lemma \ref{Lemma 0.8} and
Lemma \ref
{le:counit} we conclude that (\ref{eq:YD9}) and (\ref{eq:YD10}) hold. Thus $%
(R,\varepsilon_{},\delta,\omega)$ is an Yetter--Drinfeld
quadruple.
\end{proof}

\begin{definition}
Let $(R,\varepsilon ,\delta ,\omega )$ be a Yetter--Drinfeld
quadruple. The smash algebra $R\#H$ endowed with the bialgebra
structure described in
Theorem \ref{pro Boss} will be called the \emph{bosonization} of $%
(R,\varepsilon ,\delta ,\omega )$ and will be denoted by
$R{\#}_bH$.
\end{definition}

\begin{proposition}
\label{pr:section} Let $H$ be a Hopf algebra and let $(R,\varepsilon_{},%
\delta,\omega)$ be a Yetter--Drinfeld quadruple.

(a) The map $\pi :R\#_bH\rightarrow H,\ \pi (r\#h):=\varepsilon
\left( r\right) h$, is a map in ${}_{H}^{H}\mathfrak{M}_{H}^{H}$
and a bialgebra morphism that induces the
$(H,H)$\lin{}{}bicomodule structure of $R\#H.$

(b) The map $\sigma :H\rightarrow R\#_bH, \sigma (h):=1_R\# h$, is an $(H,H)$%
\lin{}{}bilinear section of $\pi$ and an algebra morphism that
induces the $(H,H)$\lin{}{}bimodule structure of $R\#_bH.$
\end{proposition}

\begin{proof}
(a) Since $\Delta_{R\#H}$ is defined by (\ref{eq:dsm}) we have
\begin{equation}  \label{eqxxx}
\left((\pi\otimes\pi)\Delta_{R\#H}\right)(r\#h)=\sum\varepsilon_{}\left(%
\delta^1(r\#h{}_{\left( 1\right)})\right) \delta^2(r\#h{}_{\left(
1\right)}){}_{\langle-1\rangle}h{}_{\left(
2\right)}\otimes\varepsilon_{}\left(\delta^2(r\#h{}_{\left(
1\right)}){}_{\langle 0\rangle}\right) h{}_{\left( 3\right)}.
\end{equation}
Hence by (\ref{eqdel4}) and the second relation of (\ref{eq:YD0})
we get that
\begin{equation*}
\left((\pi\otimes\pi)\Delta_{R\#H}\right)(r\#h)=\sum\varepsilon_{}(r)h{}_{%
\left( 1\right)}\otimes h{}_{\left( 2\right)}.
\end{equation*}
Clearly $\varepsilon_{H}\pi=\varepsilon_{R\#H}$, so $\pi$ is
morphism of coalgebras. The first equality in (\ref{eq:YD0})
implies easily that $\pi$
is a morphism of algebras. As the right $H$\lin{}{}module and the right $H$%
\lin{}{}comodule structures are induced from the corresponding
ones of $H$ obviously $\pi$ is right $H$\lin{}{}linear and right
$H$\lin{}{}colinear. The fact that $ \pi$ is a morphism of left
modules follows by the first relation of (\ref {eq:YD0}). To prove
that $\pi$ is a morphism of left comodules we use the second
equality of (\ref{eq:YD0}) again.

By a straightforward computation, similar to that on that we
performed to prove (\ref{eqxxx}), we get:
\begin{equation*}
\left((R\#H\otimes\pi)\Delta_{R\#H}\right)(r\#h)=
\sum\delta^1(r\#h{}_{\left(
1\right)})\varepsilon_{}\left(\delta^2(r\#h{}_{\left(
1\right)})\right)\# h{}_{\left( 2\right)}\otimes h{}_{\left(
3\right)}=\sum r\#h{}_{\left( 1\right)}\otimes h{}_{\left(
2\right)}.
\end{equation*}
This relation means that $\pi$ induces the usual right
$H$\lin{}{}comodule structure on $R\#H$. Analogously one can prove
that $\pi$ induces the diagonal left coaction on $R\#H$.

(b) Very easy, left to the reader.
\end{proof}

The following Theorem was presented by the second author during
her talk at the "2003 Spring Eastern Sectional AMS Meeting"
(Special Session on Hopf Algebras and Quantum Groups), New York,
NY (U.S.A.), 12-13 April, 2003. There we were informed that the
dual form of the equivalence $(b)\Leftrightarrow (c)$ below, as
stated in Theorem \ref{th:dualBialgInBimod},
 has already been proved by P. Schauenburg (see 6.1 and Theorem 5.1 in
 \cite{Sch2}). Nevertheless, for sake of completeness, we decided
 to keep our proof.

\begin{theorem}
\label{th:BialgInBimod2}Let $A$ be a bialgebra and let $H$ be a
Hopf algebra. The following assertions are equivalent:

(a) \ $A$ is an object in ${}_{H}^{H}\mathfrak{M}_{H}^{H}$, the map $%
\varepsilon_{A}:A\rightarrow K$ is right $H$\lin{}{}linear and $A$
becomes an algebra in $({_{H}^{H}} \mathfrak{M}{_{H}^{H},\otimes }
_{H},H)$ and a coalgebra in $({^{H}} \mathfrak{M}{^{H}},{\square
}_{H},H)$.

(b)\ There is an algebra $R$ in ${_{H}^{H}\mathcal{YD}}$ and there
are maps $\varepsilon _{R}:R\rightarrow k,\delta :R\rightarrow
R\otimes R,\omega :H\rightarrow R\otimes R$ such that
$(R,\varepsilon _{R},\delta ,\omega )$ is a Yetter-Drienfeld
quadruple and $A$ is isomorphic as a bialgebra to the bosonization
$R{\#}_bH$ of this Yetter Drienfeld quadruple.

(c) \ There are a bialgebra map $\pi :A\rightarrow H$ and an $(H,H)$%
\lin{}{}bicolinear algebra map $\sigma :H\rightarrow A$ such that
$\pi \sigma = \mathrm{Id} _{H}$.

In this case, we can choose $R=A^{Co\left( H\right) }$ the diagram
of $A$.
\end{theorem}

\begin{proof}
(a) $\Rightarrow$ (b) \ By (\ref{cl:Settings}) the canonical map $%
\phi_A:R\#H \rightarrow A$ in ${}_{H}^{H}\mathfrak{M}_{H}^{H}$ is
an isomorphism of bialgebras, where the coalgebra structure on
$R\#H$ is
defined by $\Delta_{R\#H}:=(\phi_A^{-1}\otimes\phi_A^{-1})\Delta\phi_A$ and $%
\varepsilon_{R\#H}:=\varepsilon_{A}\phi_A.$ Clearly $R\#H$ becomes
an algebra in $({}_{H}^{H}\mathfrak{M}_{H}^{H},\otimes_H,H)$ and a
coalgebra in $({}^H \mathfrak{M}^H,\square_H,H)$ because $A$ does.
Let $\varepsilon_{}$ be the restriction of $\varepsilon_{A}$ to
$R$. Furthermore let us define the $K$\lin{}{}linear maps
$\widetilde{\delta}$, $\delta$ and $\omega$ as in (\ref
{eq:delta}) and (\ref{eq:DeltOmega}). By Theorem \ref{pro Boss} $%
(R,\varepsilon_{},\delta,\omega)$ is an Yetter--Drinfeld
quadruple. Of course the bosonization of this Yetter--Drinfeld
quadruple\ is the bialgebra $R\#H$ constructed above.

(b) $\Rightarrow$ (c) \ Apply Proposition \ref{pr:section} to get
a bialgebra map $\pi' :R\#H\rightarrow H$ and an
$(H,H)$\lin{}{}bicolinear algebra map $\sigma' :H\rightarrow R\#H$
such that $\pi' \sigma'=\mathrm{Id} _{H}$ and $H$ acts on $R\#H$
via $\sigma'$
and coacts on $R\#H$ via $\pi'$. Suppose that the isomorphism between $%
R\#H$ and $A$ is given by $f:A\rightarrow R\#H$. Then $\pi=\pi' f$ and $%
\sigma=f^{-1} \sigma'$ are the required morphisms.

(c) $\Rightarrow$ (a) \ Since $\pi$ is a morphism of coalgebras, $A$ is an $%
(H,H)$\lin{}{}bicomodule with the structures induced by $\pi $.
Similarly $\sigma$ defines an $(H,H)$\lin{}{}bimodule structure on
$A$. One can check easily that these structures define a structure
of Hopf bimodule on $A$. Also with
respect to these structures $A$ becomes an algebra in $({}_{H}^{H}\mathfrak{M}%
_{H}^{H},\otimes_H,H)$ and a coalgebra in $({}^H \mathfrak{M}^H,\square_H,H)$%
. The only thing that we have to prove is that $\varepsilon_{A}$ is right $H$%
\lin{}{}linear.

Since $\sigma$ is right $H$\lin{}{}colinear and the right coaction
on $A$ is iduced by $\pi$ we have:
\begin{equation*}
\sum\sigma(h{}_{\left( 1\right)})\otimes h{}_{\left(
2\right)}=\sum\sigma(h){}_{\left( 1\right)}\otimes
\pi\left(\sigma(h){}_{\left( 2\right)}\right).
\end{equation*}
By applying $\varepsilon_{A}\otimes H$ to this equation we obtain:
\begin{equation*}
\sum\varepsilon_{A}\left(\sigma(h{}_{\left(
1\right)})\right)h{}_{\left( 2\right)}=(\pi\sigma)(h)=h.
\end{equation*}
By applying $\varepsilon_{A}$ again we get $\varepsilon_{A}\left(\sigma(h)%
\right)=\varepsilon_{H}(h)$. We conclude by remarking that
$ah=a\sigma(h)$,
since the right action of $H$ on $A$ is induced by $\sigma$. Thus $%
\varepsilon_{A}(ah)=\varepsilon_{H}(h)\varepsilon_{A}(a)$, that is $%
\varepsilon_{A}$ is right $H$\lin{}{}linear.
\end{proof}

\begin{remark}
\label{boso}Let $\left( R,\varepsilon,\delta ,\omega\right)$ be an
Yetter--Drinfeld quadruple\ such that $\omega$ is trivial. Recall
that this means that:
\begin{equation*}
\omega (h)=\varepsilon _{H}(h)1_{R}\otimes 1_{R},\ \text{for all
}h\in H.
\end{equation*}
Then it is easy to check that relations
(\ref{eq:YD0})\lin{}{}(\ref{eq:YD10}) are
equivalent to the fact $(R,\delta ,\varepsilon)$ is a bialgebra in $({}_H^H%
\mathcal{YD},\otimes,K)$. Conversely, starting with a bialgebra $%
(R,\delta,\varepsilon_{})$ in the monoidal category
${}_H^H\mathcal{YD}$, we can regard $R$ as an Yetter--Drinfeld
quadruple\ with respect to the trivial cocycle $\omega$.
Furthermore, the bosonization of this Yetter--Drinfeld quadruple\
is the usual bosonization of the bialgebra $R$, i.e. as an algebra
is the smash product $R\#H$ and as a coalgebra is the cosmash
product i.e.
\begin{eqnarray*}
\Delta _{R\#H}\left( r\#h\right) &=&\sum r^{(1)}\otimes
r^{(2)}{}_{\langle -1\rangle}h_{(1)}\otimes r^{(2)}{}_{\langle
0\rangle}\otimes h_{\left( 2\right) }, \\ \varepsilon
_{R\#H}\left( r\#h\right) &=&\varepsilon \left( r\right)
\varepsilon \left( h\right),
\end{eqnarray*}
where, by notation, $\delta(r)=\sum r^{(1)}\otimes r^{(2)}$.
\end{remark}

\begin{corollary}
\ (D. Radford)\ \label{aradf} Let $H$ be a Hopf algebra and let
$A$ be a  bialgebra. Then the following statements are equivalent:

(a) $A$ is an object in ${}_{H}^{H}\mathfrak{M}_{H}^{H}$, the
counit of $A$ is right $H$\lin{}{}linear, $A$ becomes an algebra
in $({_{H}^{H}}\mathfrak{M}{_{H}^{H},\otimes }_{H}, H)$ and a
coalgebra in $({\ _{H}^{H}}\mathfrak{M}{_{H}^{H},{\square
}_{H}},H)$.

(b) \ The diagram $R$ of $A$ is a bialgebra in
$({_{H}^{H}\mathcal{YD}},\otimes,K )$ such that $A$ is isomorphic
as a bialgebra to the usual bosonization of $R$ by $H$.

(c) \ There are two bialgebra morphisms $\pi :A\rightarrow H$,
$\sigma :H\rightarrow A$ such that $\pi \sigma =$
$\mathrm{Id}_{H}$.
\end{corollary}

\begin{proof}
(a) $\Rightarrow$ (b) \ By Theorem \ref{th:BialgInBimod2} and by
the
previous remark it is enough to show that $\omega$ is trivial. Let $%
\Delta_{R\#H}:=(\phi_A^{-1}\otimes \phi_A^{-1})\Delta_A\phi_A^{-1}$, where $%
\phi_A$ is the canonical isomorphism between $R\#H$ and $A$. Since
$\Delta_A$ is a morphism in ${}_{H}^{H}\mathfrak{M}_{H}^{H}$ from
$A$ to $A\square_H A$ we get immediately:
\begin{equation}  \label{eq:triv}
\Delta_{R\#H} (1_{R}\# h) =\left[ \left( 1_{R}\#1_{H}\right)
\otimes \left( 1_{R}\#1_{H}\right) \right] h=\sum\left(
1_{R}\#h{}_{\left( 1\right)}\right) \otimes \left(
1_{R}\#h{}_{\left( 2\right)}\right)
\end{equation}
so that
\begin{equation*}
\omega \left( h\right) =(R\otimes \varepsilon _{H}\otimes R\otimes
\varepsilon _{H})\Delta (1_{R}\# h)=\varepsilon
_{H}(h)1_{R}\otimes 1_{R}.
\end{equation*}
\indent (b) $\Rightarrow$ (c) In view of Theorem
\ref{th:BialgInBimod2} we have only to prove that $\sigma$ is a
morphism of coalgebras. For an arbitrary Yetter--Drinfeld
quadruple\ $(R,\varepsilon_{},\delta,\omega)$ we showed that
$\varepsilon_{A}(\sigma(h))=\varepsilon_{A}(h)$, see the proof of
Theorem \ref{th:BialgInBimod2}. Assume that $\theta :A\rightarrow
B$ is a bialgebra isomorphism where $\left( B,m,u,\Delta
,\varepsilon \right) $ is
the usual bosonization of a bialgebra $R$ in $({_{H}^{H}\mathcal{YD}}$ $%
,\otimes )$. Also, by the proof of the above mentioned theorem we have $%
\sigma =\phi_A\sigma'$, where $\sigma'(h)=1\# h$. By (\ref
{eq:triv}) it follows that $\sigma'$ is a morphism of coalgebras,
hence $\sigma$ is so.

(c) $\Rightarrow$ (a) \ It is trivial, as $\sigma$ is a  morphism
of coalgebras.
\end{proof}

\begin{lemma}\label{matsu}
Let $A$ be a bialgebra over a field $K$ and let $I$ be a nilpotent
ideal and coideal of $A$. If the quotient bialgebra $A/I$ has an
antipode, then $A$ is a Hopf algebra.
\end{lemma}

\begin{proof}
Let us point out that an element $x$ in a ring $R$ is invertible
if it is invertible modulo a nil ideal $L$ of $R$. We apply this
to the ring $R=\Hom_K(A,A)$ endowed with the convolution product,
to the nil ideal $L=\Hom_K(A,I)$ and to $x=\mathrm{Id}_{A}$. The
quotient $R/L$ is isomorphic to the algebra $\Hom_K(A,A/I)$ and
through this identification the class of $\mathrm{Id}_{A}$
corresponds to the canonical projection of $p:A\rightarrow A/I$.
We conclude by remarking that the inverse of $p$ in
$\Hom_K(A,A/I)$ is $p\circ S$, where $S$ is the antipode of $A/I$.
\end{proof}

\begin{theorem}\label{te:A=YDq}
Let $A$ be a bialgebra over a field $K$. If the Jacobson radical
$J$ of $A$ is a nilpotent coideal such that $H:=A/J$ is a Hopf
algebra which has an $ad$\lin coinvariant integral and that every
canonical map $A/J^{n+1}\ra A/J^n$ splits in $\hm^H$, then $A$ is
isomorphic as a bialgebra with the bosonization $R\#_bH$ of a
certain Yetter--Drinfeld quadruple\
$(R,\varepsilon_{},\delta,\omega)$. In fact $A$ and $R\#_bH$ are
isomorphic Hopf algebras.
\end{theorem}

\begin{proof}
By Theorem \ref{te:section} there is an $(H,H)$\lin{}{}bicolinear
algebra
section $\sigma:H\rightarrow A$ of the canonical projection $%
\pi:A\rightarrow H$. We conclude by applying Theorem
\ref{th:BialgInBimod2} and Lemma \ref{matsu}.
\end{proof}

\begin{theorem}\label{coro:A=YDq} Let $A$ be a bialgebra over a field $K$. If the Jacobson radical
$J$ of $A$ is a nilpotent coideal such that $H:=A/J$ is a Hopf
algebra which is both semisimple and cosemisimple (e.g. when $H$
is semisimple over a field of characteristic $0$), then $A$ is
isomorphic as a bialgebra with the bosonization of a certain
Yetter--Drinfeld quadruple\ $(R,\varepsilon_{},\delta,\omega)$. In
fact $A$ and $R\#_bH$ are isomorphic Hopf algebras.
\end{theorem}

\begin{proof}
Apply Corollary \ref{co:SectonBicolin}, Theorem
\ref{th:BialgInBimod2} and Lemma \ref{matsu}.
\end{proof}

\section{Dual results and applications}
Of course all result of the previous sections can be dualized.
Because this process is based only on some elaborate computation
and does not require new ideas we shall just state the main
results that we shall use in this part of the paper.
\begin{claim}
We start by defining the Hochschild cohomology of a coalgebra in a
monoidal category $(\calM,\otimes,\mathbf{1})$. A triple
$(C,\Delta,\ep{})$ such that $C$ is an object in $\calM$ and
$\Delta:C\ra C\ot C$ and $\ep{}:C\ra \mathbf{1}$ are morphisms in
$\calM$ is a coalgebra in $\calM$ if it is an algebra in the dual
monoidal category $\calM^\circ$ of $\calM$. Recall that
$\calM^\circ$ and $\calM$ have same objects but
$\calM^\circ(X,Y)=\calM(Y,X)$. Similarly, for any coalgebra in
$\calM$, we define a $(C,C)$\lin bicomodule in $\calM$ to be an
$(C,C)$\lin bimodule in $\calM^\circ$. The category of all
$(C,C)$\lin bicomodules will be denoted by $\cmc$. It is an
abelian category (if $\calM$ is so). Furthermore, the class
$\mathcal I$ of all monomorphism in $\cmc$ that have a retraction
in $\calM$ is an injective class of monomorphisms. Note that if we
regard $C$ as an algebra in $\calM^\circ$ then $\mathcal I$ is the
projective class associated to this algebra, as in
(\ref{cl:E-Proj}).

We fix a coalgebra in a monoidal category $\cmct$. Now, for any
$(C,C)$\lin bicomodule $M\in\cmc$, we define the Hochschild
cohomology of $C$ with coefficients in $M$ by:
\[
\Hu(M,C)=\mathbf{Ext}^\bullet_{\mathcal{I}}(M,C) ,
\]
where $\mathbf{Ext}^\bullet_{\mathcal{I}}(M,-)$ are the relative
left derived functors of $\cmc(M,-)$. Note that $\Hu(M,C)$ is the
Hochschild cohomology of the algebra $C$ with the coefficients in
$M$ (regarded as objects in $\calM^\circ$).
\end{claim}

\begin{definition}A coalgebra in $\calM$ is called \emph{coseparable} if
and only if the comultiplication $\Delta:C\ra C\ot C$ has a
retraction in $\cmc$.
\end{definition}
\begin{proposition}
Let $C$ be a coalgebra in $\calM$. The following assertions are
equivalent:

(a) $C$ is coseparable.

(b) $C$ is $\mathcal I$-injective in $\cmc$.

(c) $\Hb{H}{1}(M,C)=0$, for all $M\in\cmc$.

(d) $\Hb{H}{n}(M,C)=0$, for all $M\in\cmc$, for all $n>0$.

\end{proposition}
\begin{proof}
Regard $C$ and $M$ as objects in $\calM^\circ$ and apply
Proposition \ref{prop:sep} and Theorem \ref{te:sep}.
\end{proof}

\begin{claim}
Also by working in $\calM^\circ$ we obtain the natural definition
of a formally smooth coalgebra $C$ in a monoidal category $\cmct$.

Let us consider a morphism of coalgebras $i:D\ra E$ which has a
retraction in $\calM$. We define $D\wedge
D:=\ker(\pi\ot\pi)\Delta_E$, where $(M,\pi)$ is the cokernel of
$i$. Note that $D\wedge D=E$ if and only if
$(\mathrm{Ker}\,i)^2=0$, where now $i$ is regarded as a morphism
of algebras in $\calM^\circ$ from $E$ to $D$. Hence by applying
Theorem \ref{X coro UniE} we have proved the following theorem.
\end{claim}

\begin{theorem}
\label{te:SmoothCoalg}Let $C$ be a coalgebra in
$(\mathcal{M},\otimes ,\mathbf{1} )$. Then the following
conditions are equivalent:

(a) $\Hb{H}{2}( M,C) =0$, for all $M\in \cmc$.

(b) If $i :D\rightarrow E$ is a coalgebra homomorphism that has a
retraction in $ \mathcal{M}$ and $D\wedge D=E$, then any coalgebra
homomorphism $f:D\rightarrow C$ can be extended to a coalgebra
homomorphism $E\rightarrow C.$
\end{theorem}

\begin{definition} A coalgebra $C$ will be called
\emph{formally smooth} if it satisfies one of the above equivalent
conditions.
\end{definition}
\begin{corollary}
Any coseparable coalgebra in a monoidal category is formally
smooth.
\end{corollary}

\begin{lemma}
\label{nescio}
Let $H$ be  Hopf algebra.

a) $H$ is coseparable as a coalgebra in $\mm_H$  if and only if
$H$ is cosemisimple.

b) $H$ is coseparable as a coalgebra in $_H\mm_H$  if and only if
there is an $ad$\lin invariant integral $\lambda\in H^\ast$ (see
Definition \ref{ad-invariant}). In particular, if $H$ is
semisimple and cosemisimple, then $H$ is coseparable in $_H\mm_H$.
\end{lemma}

\begin{proof}
Dual to Lemma \ref{le:SepInM^H}.
\end{proof}

\begin{theorem}
\label{te:SectionBicomod} \ Let $H$ be a Hopf algebra.

a) Let $ C$ be a coalgebra in $\mm_H$. If the coradical  $C_0$ of
$C$ is $H$ then there is a coalgebra map $\pi_C:C\rightarrow H$
which is a morphism in $\mm_H$ such that
$\pi_C\rest{H}=\mathrm{Id}_H$.

b) Let $ C$ be a coalgebra in $_H\mm_H$.  If $C_0=H$, $H$ has an
$ad$\lin invariant integral and every $C_n$ is a direct summand in
$C_{n+1}$ as an object in $_H\mm_H$, then there is a coalgebra map
$\pi_C:C\rightarrow H$ which is a morphism in $_H\mm_H$ such that
$\pi_C\rest{H}=\mathrm{Id}_H$.

c) Let $ C$ be a coalgebra in $_H\mm_H$. If $C_0=H$ is semisimple
then a morphism $\pi$ as in (b) exists.
\end{theorem}

\begin{proof}
a) \ Let us consider the coradical filtration
$(C_n)_{n\in\mathbb{N}}$. Obviously, every $C_n$ is a coalgebra in
$\mm_H$ and $C_{n+1}=C_n\wedge C_n$ (the wedge product is
performed in the coalgebra $C_{n+1}$). Indeed, by the definition
of the coradical filtration we have:
\begin{equation}\label{eq:corad}
C_{n+1}=\{x\in C\mid \Delta(x)\in C\ot C_n+C_0\ot C\},
\end{equation}
for every $n\geq 0$. Let $f_0:=\mathrm{Id}_H=\mathrm{Id}_{C_0}$.
Let us assume that we have constructed $f_0,\cdots,f_n$ morphisms
of coalgebras in $\mm_H$ such that $f_{i+1}\!\rest{{C_i}}=f_i:C_i
\ra H$, for all $i\in\{0,\cdots,n-1\}$. On the other hand, by
(\ref{eq:corad}), $C_{n+1}/C_n$ becomes a right $H=C_0$\lin
comodule with the structure induced by $\Delta$. Hence
$C_{n+1}/C_n$ is an object in $\mm^H_H$, so it is free as a right
$H$\lin module (by the fundamental theorem for Hopf modules). In
conclusion the inclusion $C_n\subseteq C_{n+1}$ has a retraction
in $\mm_H$. Since $H=C_0$ is cosemisimple, by Lemma
\ref{nescio}(a), it is coseparable in $\mm_H$ so that we can apply
Theorem \ref{te:SmoothCoalg} to find a morphism of coalgebras
$f_{n+1}:C_{n+1}\ra H$ such that $f_{n+1}\!\rest{{C_n}}=f_n$.
Hence there is a unique morphism of coalgebras $\pi_C:C\ra H$ in
$\mm_H$ such that $\pi_C\!\rest{{C_n}}=f_n$ for all natural
numbers $n$.

b) By the Lemma \ref{nescio}(b), $H$ is coseparable in $_H\mm_H$
and moreover by assumption $C_n$ is a direct summand of $C_{n+1}$
as an object in $_H\mm_H$. Hence we can apply Theorem
\ref{te:SmoothCoalg}, in the case when $\mathcal{M}={_H\mm_H}$, to
find a morphism of coalgebras $f_{n+1}:C_{n+1}\ra H$ in $_H\mm_H$
such that $f_{n+1}\!\rest{{C_n}}=f_n$.

c) Since $H$ is semisimple, it is separable so that the category
$_H\mm_H$ is semisimple. Hence any $C_n$ is a direct summand of
$C_{n+1}$ as a $H$-bimodule. Since $H=C_0$ is also cosemisimple,
by Lemma \ref{le:semicose} $H$ has an $ad$\lin invariant integral
and so that conclude by b).
\end{proof}

\begin{corollary}\label{co:Mas1}
\ Let $A$ be a Hopf algebra such that $A_0$, the coradical of $A$,
is a Hopf subalgebra.

a) There is a coalgebra map $\pi:A\rightarrow A_0$  which is a
morphism in $\mm_{A_0}$ such that
$\pi\rest{{A_0}}=\mathrm{Id}_{A_0}$.

b) Suppose that $A_0$ has an $ad$\lin invariant  integral and
every $A_n$ is a direct summand of $A_{n+1}$ as an $(A_0,A_0)$\lin
bimodule. Then there is a coalgebra map $\pi_A:A\rightarrow A_0$
which is a morphism in $_{A_0}\mm_{A_0}$ such that
$\pi_A\rest{{A_0}}=\mathrm{Id}_{A_0}$.
\end{corollary}
\begin{remarks}
a) A. Masuoka informed us that the first statement of Theorem
\ref{te:SectionBicomod} follows easily from \cite[Theorem
4.1]{Mas}.

 b) Statement (a) in Corollary \ref{co:Mas1} has already
been proved by Masuoka, see \cite[Theorem 3.1]{Mas}.
\end{remarks}
\begin{claim}\label{cl:SmashCoalgebra}
Let $H$ be a cosemisimple Hopf algebra. Suppose that $C$ is a
coalgebra in $\mm_H$ such that the coradical of $C$ is $H$. Then
by the above theorem there is a coalgebra map $\pi_C:C\ra H$ which
is right $H$\lin linear and $\pi_C(h)=h$, for any $h\in H$. Since
$\pi_C$ is a morphism of coalgebras it follows that $C$ is an
$(H,H)$\lin bicomodule. In fact, as $\pi_C$ is a morphism of right
$H$\lin modules, we can prove easily that $C$ is an object in
$\hm^H_H$, the comultiplication $\Delta_C$ is a morphism in
$\hm^H_H$ and $\mathrm{Im}\,\Delta_C\subseteq C\cotH C$. Let $R$
be the subspace of right coinvariant elements in $C$. By the
fundamental theorem of Hopf modules $\varphi_C:R\ot H\ra C$, $r\ot
h\mapsto rh$ is an isomorphism in $\mm^H_H$. It is also left
$H$\lin colinear, since $C$ belongs to $\hm^H_H$. Thus the
comultiplication $\Delta_C:C\ra C\cotH C$ can be identified with a
morphism $\Delta_{R\ot H}:R\ot H\ra (R\ot H)\cotH (R\ot H)$ in
$\hm^H_H$. By Lemma \ref{gamma} and (\ref{gammac}) there is
$\td:R\ot H\ra R\ot R$ such that relations (\ref{eq:delta}) and
(\ref{eqdd}) hold. Since $\Delta_{R\ot H}$ is right $H$\lin linear
we have:
\[
\Delta_{R\ot H}(r\ot h)=\Delta_{R\ot H}(r\ot 1)h= \sum
\left(\tilde{\delta}^{1}(r\ot 1)\ot \tilde{\delta} ^{2}(r\ot
1){}_{\langle-1\rangle}h_{(1)}\right)\otimes \left(\tilde{
\delta}^{2}(r\ot 1){}_{\langle 0\rangle}\ot h_{2}\right).
\]
If $\de(r)=\td(r\ot 1)$, and we use the notation $\de(r)=\sum
r^{(1)}\ot r^{(2)}$ then the above relation becomes: $$\Delta
(r\#h) =\sum \left(r^{( 1) } \ot r^{(
2)}_{\langle-1\rangle}h_{(1)}\right)\ot \left(r^{(
2)}_{\langle0\rangle}\ot h_{(2)}\right).$$ An easy computation
shows that $\de:R\ra R\ot R$ is coassociative and left $H$\lin
colinear since $\Delta_{R\ot H}$ is so. Moreover, if $\ep{}$ is
the restriction of $\ep{C}$ to $R$, then $\ep{}$ is a counit for
$\de$. In conclusion we have proved the following proposition.
\end{claim}

\begin{proposition}\label{pr:cosmash}
Let $H$ be a cosemisimple Hopf algebra. Suppose that $C$ is a
coalgebra in $\mm_H$ such that the coradical $C_0$ of $C$ is $H$.
Then $C$ is an object in $\hm^H_H$ such that $R$, the space of
right coinvariant elements of $C$, is an $H$\lin comodule
coalgebra and $C$ is isomorphic as a coalgebra, via a morphism in
$\hm^H_H$, with the smash product coalgebra $R\# H$ of $R$ by $H$.
\end{proposition}

\begin{corollary}\label{co:ProjOnR}
Keep the notation and assumptions from the preceding proposition.
Then there is a right $H$\lin linear coalgebra morphism
$\pi_R:C\ra R$ such that $\pi_R\rest{R}=\mathrm{Id}_R$, where $R$
is regarded as a right module with trivial action.
\end{corollary}

\begin{proof}
Obviously $\pi':R\# H\ra R$, given by $\pi'(r\# h)=\ep{H}(h)r$, is
a morphism of coalgebras, it is right $H$\lin linear and
${\pi'}\rest{R}=\mathrm{Id}_R$. Hence $\pi_R:=\pi'\varphi^{-1}_C$
has the same properties, as the canonical map $\varphi_C:R\# H\ra
C$ is an isomorphism of coalgebras in $^H\mm^H_H$ and
$\varphi_C(r\# 1)=r$, for every $r\in R$.
\end{proof}

\begin{lemma}\label{le:Filtration}
Let $C$ be a coalgebra. Suppose that there is a group$\,$\lin like
element $c_0\in C$ such that $C_0=Kc_0$, i.e. $C$ is connected.
Let $(C'_n)_{n\in\mathbb{N}}$ be a coalgebra filtration in $C$
such that $C'_0=C_0$. Then, for every $c\in C'_n$, we have:
\begin{equation}\label{eq:Filtration}
\Delta(c)-c\ot c_0-c_0\ot c\in C'_{n-1}\ot C'_{n-1}.
\end{equation}
In particular, if $c\in C'_1$ then $\Delta(c)=c\ot c_0+c_0\ot
c-\ep{}(c)c_0\ot c_0.$
\end{lemma}
\begin{proof}
Since $(C'_n)_{n\in\mathbb{N}}$ is a coalgebra filtration we have
$\Delta(C'_n)\subseteq\sum_{i+j=n} C'_i\ot C'_j.$ Hence there are
$c'$, $c''$ and $x\in C'_{n-1}\ot C'_{n-1}$ such that:
\begin{equation}\label{eq:filt}
\Delta(c)=c'\ot c_0+c_0\ot c''+x.
\end{equation}
By applying $\ep{}\ot C$ and $C\ot\ep{}$ to this relation we
deduce that:
\begin{eqnarray*}
c&=&\ep{}(c')c_0+c''+x_1,\\ c&=&\ep{}(c'')c_0+c'+x_2,
\end{eqnarray*}
where $x_1$, $x_2$ are in $C'_{n-1}$ since $x\in C'_{n-1}\ot
C'_{n-1}$. We conclude the first part of the lemma by substituting
$c'$ and $c''$ in (\ref{eq:filt}). Now, if $c\in C'_1$ then
$\Delta(c)=c\ot c_0+c_0\ot c+\alpha c_0\ot c_0,$ for a certain
$\alpha$ in $K$. By applying $\ep{}\ot\ep{}$ we deduce that
$\alpha=-\ep{}(c)$.
\end{proof}

\begin{claim}
Let $H$ be a cosemisimple Hopf algebra. We shall denote by
$\widehat{H}$ the set of isomorphism classes of simple left
$H$\lin comodules. It is well$\,$\lin known that for every
$\tau\in\widehat{H}$ there is a simple subcoalgebra $C(\tau)$ of
$H$ such that $\rho_V(V)\subseteq C(\tau)\ot V$, where
$(V,\rho_V)$ is an arbitrary comodule in $\tau$. Moreover, we have
$H=\oplus_{\tau\in\widehat{H}}C(\tau)$.
\end{claim}

\begin{theorem}\label{te:Corad_Filt}
Keep the notation and assumptions from the statement and the proof
of Proposition \ref{pr:cosmash}. Let $(C_n)_{n\in\mathbb{N}}$ be
the coradical filtration of $C$.

a) \ For every natural number $n$ we have $C_n\simeq R_n\# H$
(isomorphism in $^H\mm^H_H$). In particular $C_n$ is freely
generated as an $H$\lin module by elements $r\in C$ satisfying the
relation:
\begin{equation}\label{eq:C_n}
\Delta(r)=\sum r\un{-1}\ot r\un{0}+r\ot 1_H + C_{n-1}\ot C_{n-1}.
\end{equation}
\indent b) \ $C_1$ verifies the following equation:
\begin{equation}\label{eq:C_1}
C_1=C_0+\sum_{\tau\in\widehat{H}}\left(C(\tau)\wedge K
1_H\right)H,
\end{equation}
\end{theorem}

\begin{proof}
a) \ By Proposition \ref{pr:cosmash}, $C_n$ is the smash product
coalgebra $R'_n\# H$. By the construction of $R'_n$ we have
$R'_n=R\bigcap C_n$. Since $C_n$ is isomorphic in $^H\mm^H_H$ with
$R'_n\# H$ it results that $C_n$ is free as a right $H$\lin
module.

Note that $(R'_n)_{n\in\mathbb{N}}$ is not \emph{a priori} a
coalgebra filtration in $R$, since $R$ is not a subcoalgebra of
$C$ (its comultiplication is $\de$, see (\ref{cl:SmashCoalgebra})
for its definition). Recall also that we use the notation
$\de(r)=\sum r^{(1)}\ot r^{(2)}$.

Let us prove that $(R'_n)_{n\in\mathbb{N}}$ is indeed a coalgebra
filtration. Let $\pi_R$ be the coalgebra morphism from Corollary
\ref{co:ProjOnR}. Then $R'_n=\pi_R(C_n)$, so
$(R'_n)_{n\in\mathbb{N}}$ is a coalgebra filtration of $R$, as
$\pi_R$ is surjective. By \cite[Corollary 5.3.5]{Mo} the coradical
of $R$ is included in $\pi_R(H)=K1_H$, hence $R$ is connected and
$R'_0=R_0$. By Lemma \ref{le:Filtration}, applied to the
filtration $(R'_n)_{n\in\mathbb{N}}$ we deduce that:
\[
\de(r)\in r\ot 1_H+ 1_H\ot r+ R'_{n}\ot R'_{n},
\]
for any $r\in R'_{n+1}$. By induction it results that
$R'_n\subseteq R_n$, for every $n$. On the other hand, for $r\in
R_{n+1}$ we have:
\[
\de(r)\in r\ot 1_H+ 1_H\ot r+ R_{n}\ot R_{n}.
\]
Since $C$ is isomorphic to the smash product coalgebra via the
canonical map $\phi_C$ we get:
\begin{equation}\label{eq:CalcDelta}
\Delta(r)\in \sum r^{(1)}r^{(2)}\un{-1}\ot r^{(2)}\un{0}+R_{n}H\ot
R_{n}H=\sum r\un{-1}\ot r\un{0}+r\ot 1_H + R_{n}H\ot R_{n}H.
\end{equation}
If we assume, by induction, that $R_n=R'_n$ then $\Delta(r)\in
C\ot C_n+H\ot C$, that is $r\in C_{n+1}$. Thus $r\in
C_{n+1}\bigcap R=R'_{n+1}$. In conclusion the filtrations
$(R'_n)_{n\in\mathbb{N}}$ and $(R_n)_{n\in\mathbb{N}}$ are equal,
and $C_n\simeq R_n\# H$. Note also that, by (\ref{eq:CalcDelta}),
every element in $R_{n}$ satisfies (\ref{eq:C_n}), so (a) is
proved.

b) \ By the proof of the first part it follows that every $R_n$ is
a subobject in $^H\mm^H_H$ of $R$. Let us decompose $R_1$ as a
direct sum of left $H$\lin comodules:
\begin{equation}\label{eq:Decomp}
R_1=K1_H\oplus R'_1=K1_H\oplus\left(\oplus_{i=1}^n V_i\right),
\end{equation}
where each $V_i$ is simple. Let $\tau_i$ be the isomorphism class
of $V_i$. Take $i\in \{1,\dots,n\}$ and $r\in V_i$. As in the
proof of (\ref{eq:CalcDelta}), by using the second equality in
Lemma \ref{le:Filtration}, one can show that:
\[
\Delta(r)=\sum r\un{-1}\ot r\un{0}+r\ot 1_H -\ep{}(r)1_H\ot
1_H=\sum r\un{-1}\ot r\un{0}+\left(r-\ep{}(r)1_H\right)\ot 1_H.
\]
Hence $\Delta(r)\in C(\tau_i)\ot C+C\ot K1_H$ which proves that
$r\in C(\tau_i)\wedge K1_H$. Thus, in view of the decomposition
(\ref{eq:Decomp}), we have proved the inclusion
``$\;\subseteq\;$'' of (\ref{eq:C_1}), as $C$ is generated as a
right $H$\lin module by $R$. The other inclusion is trivial since,
for $\tau \in \widehat{H}$ and $c\in C(\tau)\wedge K1_H$, we have:
\[
\Delta(c)\in C(\tau)\ot C+C\ot K1_H\subseteq H\ot C+C\ot H.
\]
Thus $c\in H\wedge H=C_1$, so we deduce $(C(\tau)\wedge
K1_H)H\subseteq C_1$, as $C_1$ is a right submodule of $C$.
\end{proof}
\begin{remark}
Let $A$ be a Hopf algebra such that $A_0$, the coradical of $A$,
is a subalgebra. In \cite[Lemma 4.2]{AS5} it is shown that
equation (\ref{eq:C_1}) holds true for $C:=grA$. In \cite[Remark
3.2]{CDMM} it is pointed out that the proof of (\ref{eq:C_1}),
given in \cite{AS5} for $grA$, also works  in the case $C:=A$,
since $A$ is a cosmash by \cite[Theorem 3.1]{Mas}.
\end{remark}
\begin{definition}
Let $H$ be a Hopf algebra and let $(R,\delta,\varepsilon)$ be a
coalgebra in the category $({_{H}^{H}\mathcal{YD}},\otimes,K)$.
Assume that $ m:R\otimes R\rightarrow R,$ and $\xi :R\otimes
R\rightarrow H$ are $K$\lin{}{}linear maps and fix one element
$1\in R$. The quadruple $ (R,1,m ,\xi)$ will be called \emph{dual
Yetter--Drinfeld quadruple} if and only if, for all $r,s,t\in R$
and $h\in H$, the following relations hold:
\begin{eqnarray}
&^{h}1=\varepsilon _{H}(h)1\qquad \text{and}\qquad \rho
_{R}(1)=1_{H}\otimes 1;  \label{eq:YD0'} \\ &\delta (1)=1\otimes
1\qquad \text{and}\qquad \varepsilon (1)=1_K; \label{eq:YD1'} \\
&^{h}m(r\otimes s)=\sum m(\,{}^{h_{(1)}}r \otimes{}^{h_{(2)}}s);
\label{eq:YD2'} \\ &\sum \xi ({}^{h_{(1)}}r\otimes
{}^{h_{(2)}}s)=\sum h_{(1)}\xi (r\otimes s)Sh_{(2)};
\label{eq:YD3'} \\ &\delta m=(m\otimes m)\delta_{R\otimes R}
\qquad \text{and}\qquad \varepsilon m=m_K(\varepsilon
\otimes\varepsilon);  \label{eq:YD4'} \\ &\Delta_H \xi=(m_H\otimes
H)(\xi\otimes H\otimes\xi)(R\otimes R\otimes\rho_{R\otimes
R})\delta_{R\otimes R} \quad\text{and} \quad\varepsilon _{H}\xi
=m_K(\varepsilon\otimes\varepsilon); \label{eq:YD5'} \\ & c_{R,H}
(m\otimes \xi)\delta_{R\otimes R}= (m_H\otimes R)(\xi\otimes
H\otimes m)(R\otimes R\otimes \rho_{R\otimes R})\delta_{R\otimes
R}; \label{eq:YD6'} \\ & m(R\otimes m)=m(m\otimes R)(R\ot R\ot
\mu_R)(R\otimes R\otimes \xi\otimes R)(\delta_{R\otimes R}\otimes
R);  \label{eq:YD7'} \\ &m_H(\xi\otimes H)(R\otimes
m\otimes\xi)(R\otimes\delta_{R\otimes R})= m_H(\xi\otimes
H)(R\otimes c_{H,R})(m\otimes\xi\otimes R)(\delta_{R\otimes
R}\otimes R);  \label{eq:YD8'} \\ &m(r\otimes 1)=r=m(1\otimes r);
\label{eq:YD9'} \\ &\xi (r\otimes 1)=\xi (1\otimes r)=\varepsilon
(r)1_{H}. \label{eq:YD10'}
\end{eqnarray}
\end{definition}

\begin{remark}
Note that these relations can be interpreted as follows:

(\ref{eq:YD0'}) \lin{}{} $1$ is left $H$\lin{}{}invariant and left
$H$\lin{}{}coinvariant;

(\ref{eq:YD1'}) \lin{}{} $1$ is a group\lin{}{}like element;

(\ref{eq:YD2'}) \lin{}{} $m$ is left $H$\lin{}{}linear;

(\ref{eq:YD3'}) \lin{}{} $\xi $ is left $H$\lin{}{}linear, where
$H$ is a module with the adjoint action;

(\ref{eq:YD4'}) \lin{}{} $m{}$ is a morphism of coalgebras, where
on $R\otimes R$ we consider the coalgebra structure that uses the
braiding $c$;

(\ref{eq:YD5'}) \lin{}{} $\xi $ is a normalized cocycle; more
generally, if $C$ is a left $H$\lin{}{}comodule coalgebra then a
map $\psi:C\ra H$ is called a non\lin{}{}commutative 1 cocycle if
\[\Delta_H(\psi(c))=\sum \psi\left(c\ro{1}\right)(c\ro{2})\un{-1}\ot
\psi\left((c\ro{2})\un{0}\right)
\]

(\ref{eq:YD6'}) \lin{}{} $\xi $ measures how far $m$ is to be a
morphism of left $H$\lin{}{}comodules (if $\xi $ is trivial, i.e.
for every $r,s\in R$ we have $\xi
(r\otimes s)=\varepsilon (r)\varepsilon (s)$, then $m$ is left $H$%
\lin{}{}colinear); we shall say that $m$ is a \emph{twisted
morphism} of left $H$\lin{}{}comodules; we shall use the notation
$m(r\otimes s)=rs,$ so equation (\ref {eq:YD6'}) can be rewritten
as follows:
\begin{equation*}
\sum (r^{(1)}{}^{r_{\left\langle-1\right\rangle
}^{(2)}}s^{(1)})_{\left\langle-1\right\rangle }\xi
(r_{\left\langle 0\right\rangle }^{(2)}\otimes s^{(2)})\otimes
(r^{(1)}{}^{r{}_{\langle-1\rangle }^{(2)}}s^{(1)})_{\left\langle
0\right\rangle }=\sum \xi (r^{(1)}\otimes {}^{r_{\left\langle
-2\right\rangle }^{(2)}}s^{(1)})r_{\left\langle-1\right\rangle
}^{(2)}s_{\left\langle-1\right\rangle }^{(2)}\otimes
r_{\left\langle 0\right\rangle }^{(2)}{}s{}_{\langle 0\rangle
}^{(2)}
\end{equation*}

(\ref{eq:YD7'}) \lin{}{} when $\xi$ is trivial then
(\ref{eq:YD7'}) is equivalent to the fact that $m$ is associative;
so, in general,
we shall say that $m$ is $\xi$\lin{}{}\emph{%
associative}; here $\mu_R$ denotes the $H$\lin{}{}action on $R$;

(\ref{eq:YD8'}) \lin{}{} we shall just say that $m$ and $\xi$ are
\emph{compatible}; it is equivalent to:
\begin{equation*}
\sum \xi (r\otimes s^{(1)}\;{}^{s_{\left\langle-1\right\rangle
}^{(2)}}t^{(1)})\xi (s_{\left\langle 0\right\rangle} ^{(2)}\otimes
t^{(2)})=\sum \xi( r^{(1)}\ {}^{r_{\left\langle-1\right\rangle
}^{(2)}}s^{(1)}\otimes \;{}^{\xi ( r_{\left\langle 0\right\rangle
}^{(2)}\otimes s^{(2)}) _{(1)}}t)\xi (r_{\left\langle
0\right\rangle }^{(2)}\otimes s^{(2)})_{(2)}
\end{equation*}

(\ref{eq:YD9'}) \lin{}{} $m$ is a \emph{unitary map} with respect to $%
1$;

(\ref{eq:YD10'}) \lin{}{} $\omega$ is a \emph{unitary map} with
respect to $1$;

\noindent Since $1$ satisfies the last two relation we shall call
it the \emph{unit} of the {dual Yetter--Drinfeld quadruple} $R$.
By analogy $m$ will be called the \emph{multiplication} of $R$.
Finally, we shall say that $\xi$ is the \emph{cocycle} of $R$.
\end{remark}

\begin{theorem}
\label{th:YDdual} Let $R$ be a coalgebra in ${_{H}^{H}\mathcal{YD}}$. If $%
1\in R$, $m :R\ot R\rightarrow R$ and $\xi :R\ot R\rightarrow H$
are linear maps then the following assertions are equivalent:

(a) $(R,1,m,\xi)$ is a dual Yetter--Drinfeld quadruple.

(b) The smash product coalgebra $R\#H$ of $R$ by $H$ is a
bialgebra with unit $1\# 1_H$ and multiplication:
\[m_{\sm}(r\#h\ot s\#k)=\sum
m\left(r^{(1)}\otimes\;{}^{r^{(2)}_{\left\langle-1\right
\rangle}h\ro{1}}s^{(1)}\right)\ot
\xi\left(r^{(2)}_{\left\langle0\right\rangle}
{\otimes}\;{}^{h\ro{2}}s^{(2)}\right)h\ro{3}k,
\]
and $R\#H$ becomes a coalgebra in $({_{H}^{H}}%
\mathfrak{M}_{H}^{H},\cotH,H)$ and an algebra in $({_{H}}%
\mathfrak{M}{_{H}},\ot_{H},H).$
\end{theorem}

\begin{definition}
Let $(R,1 ,m ,\xi)$ be a dual Yetter--Drinfeld quadruple. The
smash product coalgebra $R\#H$ endowed with the bialgebra
structure described in the preceding
Theorem will be called the \emph{bosonization} of $%
(R,1,m,\xi)$ and will be denoted by $R{\#}\,^bH$.
\end{definition}

As we already remarked, before Theorem \ref{th:BialgInBimod2}, the
equivalence $(b)\Leftrightarrow (c)$ below has already been proved
by P. Schauenburg (see 6.1 and Theorem 5.1 in
 \cite{Sch2}).

\begin{theorem}
\label{th:dualBialgInBimod}Let $A$ be a bialgebra and let $H$ be a
Hopf algebra. The following assertions are equivalent:

(a) \ $A$ is in ${}_{H}^{H}\mathfrak{M}_{H}^{H}$, $1$ is right
$H$\lin{}{}coinvariant and $A$ becomes a coalgebra in $({_{H}^{H}}
\mathfrak{M}_{H}^{H},\cotH,H)$ and an algebra in $({^{H}}
\mathfrak{M}{^{H}},{\ot }_{H},H)$.

(b)\  There is a coalgebra $R$ in ${_{H}^{H}\mathcal{YD}}$, there
is an element $1\in R$ and linear maps $m :R\ot R\rightarrow R$,
$\xi :R\ot R\rightarrow H$ such that $(R,1,m,\xi)$ is a dual
Yetter-Drinfeld quadruple and $A$ is isomorphic, as a bialgebra,
with the bosonization $R\#\,^bH$ of $(R,1,m,\xi)$.

(c) \ There are a bialgebra map $\sigma :H\rightarrow A$ and an $(H,H)$%
\lin{}{}bilinear coalgebra map $\pi :A\rightarrow H$ such that $\pi \sigma =%
\mathrm{Id} _{H}$.

In this case, we can choose $R=A^{Co\left( H\right) }$ the diagram
of $A$.
\end{theorem}

\begin{theorem}\label{te:Boson}
Let $A$ be bialgebra over a field $K$. Suppose that the coradical
$H$ of $A$ is a semisimple sub\lin bialgebra of $A$ with antipode.
Then $A$ is isomorphic as a bialgebra with the bosonization
$R\#\,^b H$ of a certain dual Yetter--Drinfeld quadruple\
$(R,1,m,\xi)$. In fact $A$ and $R\#\,^b H$ are isomorphic Hopf
algebras.
\end{theorem}

\proof The first assertion is dual to Theorem \ref{coro:A=YDq}. In
view of a famous Takeuchi's result (see \cite[Lemma 5.2.10]{Mo})
$A$ and $R\#\,^b H$ are Hopf algebras and hence they are also
isomorphic as Hopf algebras.
\endproof

\begin{example}
Let $p$ be an odd prime and let $K$ an infinite field containing a
primitive $p$--$\,$th root of the unit $\lambda .$ Let $C$ be a
cyclic group of order $p^{2}$ with generator $c$. For every $a\in
K,a\neq 0,$ let $\
A:=H\left( a\right) $ be the Hopf algebra constructed by Beatty, D\u{a}sc%
\u{a}lescu and Gr\"{u}nenfelder in \cite{BDG}. $A$ has dimension
$p^{4}$, with basis $\left\{ c^{i}x_{1}^{j}x_{2}^{r}\mid 0\leq
i\leq p^{2}-1,0\leq j,r\leq p-1\right\} $ where $c,x_{1},x_{2}$
are subject to:
\begin{gather*}
c^{p^{2}}=1,x_{1}^{p}=c^{p}-1,x_{2}^{p}=c^{p}-1, \\
x_{1}c=\lambda ^{-1}cx_{1},x_{2}c=\lambda
cx_{2},x_{2}x_{1}=\lambda
x_{1}x_{2}+a\left( c^{2}-1\right) , \\
\Delta \left( c\right) =c\otimes c,\Delta \left( x_{1}\right)
=c\otimes x_{1}+x_{1}\otimes 1,\Delta \left( x_{2}\right)
=c\otimes x_{2}+x_{2}\otimes 1.
\end{gather*}
$A$ is a pointed Hopf algebra with coradical $H:=KC$. Let $\sigma
:H\rightarrow A$ be the canonical injection and let $\pi
:A\rightarrow H$ be
the obvious projection. It is straightforward to show that$\ A,H,\pi $ and $%
\sigma $ fulfills the requirements of Theorem
\ref{th:dualBialgInBimod}(c). Let
\begin{equation*}
R=A^{coH}=\left\{ b\in A\mid \sum b_{\left( 1\right) }\otimes \pi
\left( b_{\left( 2\right) }\right) =b\otimes 1\right\}.
\end{equation*}
We have that $R$ is the $K$--$\,$subspace of $A$ spanned by the products $%
x_{1}^{j}x_{2}^{r},$ where $0\leq j,r\leq p-1.$ In view of Theorem
\ref{th:dualBialgInBimod}, one gets a dual Yetter--Drinfeld
quadruple $(R,1,m,\xi )$ such that
$A$ is isomorphic as a bialgebra, with the bosonization $%
R\#\,^{b}H$ of $R$ by $H.$ We point out that $\xi $ is not
trivial. In fact we have:
\begin{equation*}
\xi \left( x_{2}\otimes x_{1}\right) =a\left( c^{2}-1\right) .
\end{equation*}
Clearly, the dual Hopf algebra $A^{\ast }$ fulfills the
requirements of Theorem \ref{th:BialgInBimod2} with respect to
$H^{\ast },\sigma ^{\ast }$ and $\pi ^{\ast }$. Let $\iota
:R\rightarrow A$ be the canonical injection. Then we have that the
restriction $\Lambda $ of $\iota ^{\ast }$ to $\left( A^{\ast
}\right) ^{coH^{\ast }}$ \
\begin{equation*}
\Lambda :\left( A^{\ast }\right) ^{coH^{\ast }}\rightarrow R^{\ast
}
\end{equation*}
is an isomorphism. Let $\alpha :R^{\ast }\otimes R^{\ast
}\rightarrow \left( R\otimes R\right) ^{\ast }$ be the usual
isomorphism. Then we have the following commutative diagram:
\begin{equation*}
\begin{diagram}[h=2em] H^* & \rTo<{\omega} && (A^*)^{coH^*}\otimes (A^*)^{coH^*} \\
 \dTo<{\xi^* } & && \dTo>{
\Lambda \otimes \Lambda } \\ (R \otimes R)^* && \lTo<{\alpha} & R^* \otimes R^* \\
\end{diagram}
\end{equation*}
In fact we have:
\begin{eqnarray*}
\left[ \left( \alpha \left( \Lambda \otimes \Lambda \right) \omega
\right) \left( \chi \right) \right] \left( r\otimes s\right)
&=&\left( \varepsilon
_{R}\#\chi \right) m_{R\#H}\left( r\#1\otimes s\#1\right)  \\
&=&\sum {\varepsilon _{R}}\left[ m\left( r^{(1)}\otimes
\;{}^{r_{\left\langle -1\right\rangle }^{(2)}}s^{(1)}\right) \right] {\chi }%
\left[ \xi \left( r_{\left\langle 0\right\rangle }^{(2)}{\otimes }%
\;s^{(2)}\right) \right]  \\
&=&\sum {\varepsilon _{R}}\left( r^{(1)}\right) {\varepsilon
_{R}}\left( {}^{r_{\left\langle -1\right\rangle
}^{(2)}}s^{(1)}\right) {\chi }\left[ \xi
\left( r_{\left\langle 0\right\rangle }^{(2)}{\otimes }\;s^{(2)}\right) %
\right]  \\
&=&\sum {\varepsilon _{R}}\left( r^{(1)}\right) \varepsilon
_{H}\left( r_{\left\langle -1\right\rangle }^{(2)}\right)
{\varepsilon _{R}}\left( s^{(1)}\right) {\chi }\left[ \xi \left(
r_{\left\langle 0\right\rangle
}^{(2)}{\otimes }\;s^{(2)}\right) \right]  \\
&=&{\chi }\left[ \xi \left( r{\otimes }\;s\right) \right] =\left[
\xi ^{\ast }\left( \chi \right) \right] \left( r\otimes s\right) .
\end{eqnarray*}
It follows that we can identify the Yetter--Drinfeld quadruple
$(\left( A^{\ast}\right) ^{coH^{\ast }},\varepsilon ,\delta
,\omega )$ with the Yetter--Drinfeld quadruple $\left( R^{\ast
},\left( u_{R}\right) ^{\ast },m^{\ast },\xi ^{\ast }\right) ,$ where $%
\left( u_{R}\right) ^{\ast }:R^{\ast }\rightarrow K$ is the evaluation at $%
1\in R$. In particular we observe that we get a nontrivial
bosonization since $\omega$ is not trivial.
\end{example}

 \noindent \textbf{Acknowledgements.} We thank A.
Masuoka for his very interesting remarks that helped us to improve
an earlier version of this paper. We also would like to thank
Nicolas Andruskiewitsch and Christian Lomp for some useful
suggestions.

\begin{minipage}[t]{14.6cm}\vspace*{2mm}\sc \footnotesize
A. Ardizzoni\ --\ University of Ferrara, Department of
Mathematics, Via Machiavelli 35, Ferrara,  I-44100, Italy
\vspace*{1mm}. {\small\it email:} {\small \rm ardiz@dm.unife.it}%

C. Menini\ --\ University of Ferrara, Department of Mathematics,
Via Machiavelli 35, I-44100, Ferrara, Italy\vspace*{1mm}.
{\small\it email:} {\small\ttfamily \rm men@dns.unife.it}%

D. \c{S}tefan\ --\ University of Bucharest, Faculty of
Mathematics, Strada Academiei 14, Bucharest, RO-70109, Romania.
{\small\it email:} {\small\ttfamily \rm dstefan@al.math.unibuc.ro}
\end{minipage}

\begin{thebibliography}{100}
\bibitem[AD]{AD}N. Andruskiewitsch, J. Devoto,
Extensions of Hopf algebras, Algebra i Analiz {\bf 7} (1995),
22--61, and also as a preprint in MPI-Bonn series in June 1993.

\bibitem[AS1]{AS1}N. Andruskiewitsch, H.-J. Schneider, Hopf algebras of order
$p^2$ and braided Hopf algebras of order $p$, J. Algebra {\bf 199}
(1998), 430--454.

\bibitem[AS2]{AS2}N. Andruskiewitsch, H.-J. Schneider, Lifting of
quantum linear spaces and pointed Hopf algebras of order $p^3$, J.
Algebra {\bf 209} (1998), 658--691.

\bibitem[AS3]{AS3}N.\,Andruskiewitsch and H.\,-J.\,Schneider,
Finite quantum groups and Cartan matrices, Adv. Math. {\bf 154}
(2000), 1--45.

\bibitem[AS4]{AS4}N.\,Andruskiewitsch and H.\,-J.\,Schneider,
\emph{Pointed Hopf algebras}, in S. Montgomery and H.-J. Schneider
(eds.), "New Directions in Hopf Algebras", MSRI Publ. 43,
Cambridge Univ. Press, 2002, pp. 1-68.

\bibitem[AS5]{AS5}N.\,Andruskiewitsch and H.\,-J.\,Schneider,
\textit{On the coradical filtration of Hopf algebras whose
coradical is a Hopf subalgebra}, Bol. Acad. Nacional Cienc.
C\'ordoba {\bf 65} (2000), 45--50.

\bibitem[BDG]{BDG}M. Beattie, S. D\u{a}sc\u{a}lescu, L.
Gr\"{u}nenfelder, \emph{On the number of types of finite
dimensional Hopf algebras}, Invent. math. \textbf{136} (1999),
1--7.

\bibitem[BDK]{BDK}  B. Bakalov, A. D'Andrea, V.G. Kac, \emph{Theory of
finite pseudoalgebras}, Adv. Math. \textbf{162} (2001), 1--140.

\bibitem[CDMM]{CDMM} C. C\u{a}linescu, S. D\u{a}sc\u{a}lescu, C.
Menini and A. Masuoka, \emph{Quantum lines over non$\,$\lin
cocommutative Hopf algebras}, preprint 2002.

\bibitem[CQ]{CQ} J. Cuntz and D. Quillen, \emph{Algebra extensions
and nonsingularity}, J. of AMS, \textbf{8} (1995), 251--289.

\bibitem[DNC]{DNC}  S. D\u{a}c\u{a}lescu, C. N\u{a}st\u{a}sescu and \c{S}.
Raianu, \emph{Hopf Algebras}, Marcel Dekker, 2001.

\bibitem[Do]{Doi}  Y. Doi, \emph{Homological coalgebra}, J. Math. Soc.
Japan \textbf{33} (1981), 31--50.

\bibitem[EG]{EtGe}  P. Etingof, S. Gelaki, \emph{On finite-dimensional
semisimple and cosemisimple Hopf algebras in positive
characteristic.} Internat. Math. Res. Notices \textbf{16} (1998),
851--864.

\bibitem[Ka]{Ka}  C. Kassel, \emph{Quantum groups}, Graduate Text in
Mathematics \textbf{155}, Springer, 1995.

\bibitem[Maj1]{Maj1} S. Majid, \emph{Crossed products by braided groups and bosonization},
J. Algebra \textbf{163} (1994), 165--190.

\bibitem[Maj2]{Maj2}  S. Majid, \emph{Foundations of quantum group theory},
Cambridge University Press, 1995.

\bibitem[Mas]{Mas} A.\,Masuoka, \textit{Hopf cohomology vanishing via
approximation by Hochschild cohomology}, in press for the Banach
Center Publ. Vol. \textbf{61}.

\bibitem[Mo]{Mo}  S. Montgomery, \emph{Hopf Algebras and their actions on
rings, }CMBS Regional Conference Series in Mathematics
\textbf{82}, 1993.

\bibitem[HS]{HiSt}  P. J. Hilton and U. Stambach, \emph{A course in
Homological algebra}, Graduate Text in Mathematics \textbf{4},
Springer, New York, 1971.

\bibitem[Ra1]{Rad}  D. Radford, \emph{Hopf algebras with projection}, J.
Algebra \textbf{92} (1985), 322--347.

\bibitem[Ra2]{Radf}  D. Radford, \emph{The Trace Function and Hopf Algebras}%
, J. Algebra \textbf{163} (1994), 583--622.

\bibitem[Raf]{Raf}  M. D. Rafael, Comm. Aalgebra \textbf{18
}(1990), 1445-1459.

\bibitem[Sch1]{Sch1}  P. Schauenburg, \emph{Hopf Modules and
Yetter\lin{}{}Drinfeld Modules}, J. Algebra \textbf{169} (1994),
874--890.

\bibitem[Sch2]{Sch2}  P. Schauenburg, \emph{The structure of Hopf algebras with a weak projection}, Algebr. Represent. Theory \textbf{3} (1999),
187-211.

\bibitem[SvO]{SvO}  D. \c{S}tefan and F. van Oystaeyen,
\emph{The Wedderburn$\,$\lin Malcev theorem for comodule
algebras}, Comm.\ Algebra \textbf{27} (1999), 3569--3581.

\bibitem[Sw]{Sw}  M. Sweedler, ``Hopf Algebras'', Benjamin, New York,
1969.

\bibitem[TW]{TW} E.J. Taft and R.L. Wilson, \emph{On antipodes in pointed Hopf
algebras}, J. Algebra \textbf{29} (1974), 27--42.

\bibitem[We]{We}  C. Weibel, \emph{An introduction to homological algebra, }%
Cambridge Studies in Advanced Mathematics \textbf{38}, Cambridge
University Press, 1994.
\end{thebibliography}
\end{document}